\documentclass{amsart}
\usepackage{amssymb,latexsym,graphics}
\usepackage[mathscr]{eucal}

\newtheorem{thm}{Theorem}[section]

\newtheorem{lem}[thm]{Lemma}
\newtheorem{prop}[thm]{Proposition}

\theoremstyle{remark}
\newtheorem{rem}[thm]{Remark}

\theoremstyle{definition}

\newcommand{\lp}[2]{\Vert \, #1 \, \Vert_{#2}}

\newcommand{\ret}{\vspace{.3cm}}

\begin{document}

\title[Global Regularity for $(4+1)$ NLW]
{Global Regularity and Scattering for General Non-Linear Wave Equations 
\textbf{II}. $(4+1)$ Dimensional Yang--Mills Equations in the Lorentz 
Gauge}
\author{Jacob Sterbenz}
\address{Department of Mathematics, Princeton University,
Princeton, NJ 08540}
\email{sterbenz@math.princeton.edu}
\thanks{The author was supported in part by an NSF postdoctoral 
fellowship.}
\date{}
\dedicatory{}
\commby{}


\begin{abstract}
We continue here with previous investigations \cite{STR_generic1}
on the global behavior of general type non-linear wave equations
for a class of small, scale-invariant initial data. In particular,
we show that the $(4+1)$ dimensional Yang-Mills equations are
globally well posed with asymptotically free behavior for a wide
class of initial data sets which include general charges. The
method here is based on the use of a new set of Strichartz estimates
for the linear wave equation which incorporates extra weighted 
smoothness assumptions with respect to the angular variable, along
with the construction of appropriate micro-local function spaces
which take into account this type of additional regularity.
\end{abstract}

\maketitle

\section{Introduction}

The goal of this paper is to give a proof and description
of the global regularity properties of a wide class of
non-linear wave equations on $(4+1)$ dimensional Minkowski
space. This is a continuation of our previous work 
\cite{STR_generic1}. All the equations we shall consider here 
are semi-linear wave equations with derivative non-linearities.
The generic form for such an object can be written as follows:
\begin{equation}
	\Box \phi^I \ = \ \mathcal{N}(\phi^I , \nabla \phi^I) \ ,
	\label{generic_system}
\end{equation}
where $\Box = -\partial_t + \Delta_x$ is the usual D'Lambertian
with $\Delta_x = \partial_1^2 +\partial_2^2 + \partial_3^2 
+\partial_4^2$ the Laplacean on $\mathbb{R}^4$.
Here, the superscript notation in the $\phi^I$ denotes that
we may be considering a system of equations, where the $I$
can be thought of as an index. The non-linearity $\mathcal{N}$
on the  right hand side of \eqref{generic_system} is some function
of $\phi^I$ and its first partial derivatives, collectively
denoted by $\phi^I$. However, we do not allow $\mathcal{N}$
to contain second order derivatives of $\phi$. Also, we will
restrict ourselves here to the case where the non-linearity
$\mathcal{N}$ has constant coefficients and is a polynomial
of degree $2$ or higher in the vector $(\phi^I,\nabla\phi^I)$
with no interaction of the type $\phi^2$. This is not a severe
restriction because it still includes the class of gauge--field 
equations on Minkowski space, which is one of the main
motivations of this work: \\

Specifically, let $(G,\mathfrak{g})$ be a compact, semi--simple
Lie group with Lie algebra $\mathfrak{g}$. For a given set 
of $\mathfrak{g}$--valued functions on Minkowski space $\{A_\alpha\}$,
$\alpha = 0,\ldots, 4$, we form the curvature 2-form:
\begin{equation}
	F_{\alpha\beta} \ = \ \partial_\alpha A_\beta - \partial_\beta
	A_\alpha + [A_\alpha,A_\beta] \ . \label{YM_curv}
\end{equation} 
Then $F$ is said to satisfy the Yang--Mills equations if:
\begin{equation}
	D^\beta F_{\alpha\beta} \ = \ 0 \ , \label{YM_eqs}
\end{equation}
where $D_\alpha F = \partial_\alpha F + [A_\alpha,F]$ 
denotes the gauge covariant derivative of $F$ in the
direction of $\partial_\alpha$. Expanding out the equation
\eqref{YM_eqs} in terms of the gauge potentials $\{A_\alpha\}$,
we arrive at the following second order system of PDE:
\begin{equation}
	\Box A_\beta \ = \ \partial_\beta 
	\partial^\alpha A_\alpha
        - \Big([\partial^\alpha A_\alpha,A_\beta] + 
	[A_\alpha,\partial^\alpha
        A_\beta ]+ [A^\alpha,F_{\alpha\beta}]\Big) \ . 
	\label{2nd_order_A}
\end{equation}
If we now make the a-priori assumption that:
\begin{equation} 
	\partial^\alpha A_\alpha=0 \ , \label{lorentz_gauge}
\end{equation}
the so called \emph{Lorentz gauge} condition, then \eqref{2nd_order_A}
reduces to:
\begin{equation}
	\Box A_\beta \ = \  
	- [A_\alpha,\partial^\alpha
        A_\beta ]-  [A^\alpha,F_{\alpha\beta}]\ . 
	\label{2nd_order_A_lorentz}
\end{equation}
In fact, it turns out that if $\{A_\alpha\}$ is a solution to 
\eqref{2nd_order_A_lorentz} such that at $t=0$ one has the gauge
condition:
\begin{align}
	\partial^\alpha A_\alpha \, (0) \ &= \ 0 \
	&\partial^\alpha \partial_t A_\alpha \, (0) \ &= \ 0 \ ,
	\label{initial_lorentz}
\end{align}
then \eqref{lorentz_gauge} is satisfied for all times where
the solution in sufficiently smooth. That is, the Lorentz gauge
condition propagates. Notice that if the second condition
in \eqref{initial_lorentz} is to be satisfied, then from the 
equations \eqref{2nd_order_A_lorentz}, we must have that the 
temporal potential $A_0$ satisfies the following elliptic
constraint equation at $t=0$:
\begin{equation}
	\Delta_x A_0 (0) \ = \ 
	\partial^i A_i (0) + [A_\alpha,\partial^\alpha
        A_\beta ](0) +  [A^\alpha,F_{\alpha\beta}](0) \ . 
	\label{elliptic_constraint}
\end{equation}
This implies that in general, the initial data for the system
\eqref{2nd_order_A_lorentz} (and indeed for the 
system \eqref{2nd_order_A}) cannot decay better than
$r^{-2}$. This fact is known as the \emph{charge problem},
and causes certain difficulties in the global theory of
\eqref{YM_eqs}. We will discuss this in more detail shortly.\\

We now return to the more general discussion of equations
of type \eqref{generic_system}.
Our main concern will be the global in time regularity
properties for these kind of systems. This type of 
question has been considered 
by many authors for various spatial dimensions, 
and it is not possible to give here a complete 
account of all the progress that has been made to date. 
In the case of $4$ spatial dimensions, that is $(4+1)$
dimensional Minkowski space, the first general 
theory of the global behavior of non-linear systems of the form
\eqref{generic_system} was given by the breakthrough work
of S. Klainerman \cite{Kl_uniform}. Specifically, he showed that if the
non-linearity on the right hand side of 
\eqref{generic_system} is schematically of the form:
\begin{equation}
	\Box \phi^I = |\nabla \phi^I|^2 \ , \label{generic_wm}
\end{equation}
and the Cauchy data:
\begin{align}
	\phi^I(0) \ &= \ f^I \ , &\partial_T\phi^I\, (0)
	\ &= \ g^I \ , \label{cauchy_data}
\end{align}
is sufficiently smooth and decays sufficiently fast at
(space-like) infinity, then as long as the corresponding norms are
small enough a global solution to \eqref{generic_wm}
with this initial data exists. The method of that paper was
based on controlling the $L^\infty$ norm $\nabla\phi$ 
to such an extent that Duhamel's principle could be
used globally in time. That is, one makes crucial use of the 
\emph{uniform} decay of the solutions to the non-linear problem 
\eqref{generic_wm}. This decay is provided by certain weighted
energy estimates which can naturally be recover by commuting
the various weights (vector fields) with the linear equation
on the left hand side of \eqref{generic_wm}.
At the outset, the result of 
\cite{Kl_uniform} did not include the case where the non-linear system
\eqref{generic_system} takes the form (again schematically):
\begin{equation}
	\Box \phi^I = \phi^I\, \nabla \phi^I \ , \label{generic_ym}
\end{equation}
which is the general interaction type of the Lorentz Yang--Mills
equations \eqref{2nd_order_A_lorentz}.
A somewhat more involved argument is needed to get around the
fact that the natural quantity one can gain $L^\infty$ control of
via the energy method is $\nabla \phi^I$ instead of $\phi^I$.
This problem was handled by H\"ormander in \cite{Hor_4d}, who used certain
Riesz potential modifications of the usual energy to 
gain the needed $L^\infty$ control on $\phi^I$.  
Notice that in some sense, the problem \eqref{generic_ym}
is critical with respect to the decay of $\phi^I$. This can easily
be seen by integrating the naive asymptotic one would have
for $\nabla \phi^I$:
\begin{equation}
	|\nabla \phi^I| \ \sim  \ 
	\frac{const.}{(t+r)^\frac{3}{2}(|t-r| +1)^\frac{1}{2}} \ , \notag
\end{equation}
to obtain $|\phi^I| \sim t^{-1}$, which just fails to be integrable
globally in time.\footnote{Of course one can obtain the correct decay
directly for $\phi^I$ through the use of Morawetz type multipliers.
However, the price one pays for this is the presence of extra
weights in the energy integral. Once these are taken into account,
one will see again that (at least alone the forward light cone t=r)
the decay in $(t+r)$ is critical.} Heuristically, this means that 
one cannot close a boot-strapping argument for the system 
\eqref{generic_ym} without allowing either modifying the energy in some
way, allowing it to grow, 
or finding some family of exact space--time integrals to control
error estimates which come from the differentiation of the
non-linearity. A major drawback of the results
\cite{Kl_uniform}--\cite{Hor_4d} is that they assume the initial
data \eqref{cauchy_data} is either compactly supported or
decays in such a way that it is
$L^2$. That is, the data is assumed to decay like $r^{-2 - \epsilon}$
as $r\to \infty$. As we have mentioned in our previous discussion
of the Yang-Mills equations \eqref{YM_eqs}, this kind of decay rate
is not quite attainable.\footnote{However, note that one is only
off by logarithmic divergence. It is likely that this problem
can be overcome by using an appropriate fractionally weighted modification
of the usual vector-field method. We will not pursue these ideas here,
as our approach is much more general and includes initial data
sets which decay at a rate that would be highly singular
to any straight forward modification of the vector-field technique.} 
For these reasons, as well as its intrinsic interest, we introduce
here a completely different method for studying the global behavior
of \eqref{generic_ym} which is  based on recent advances in the
low regularity theory of general non-linear wave equations. 
As a point of comparison, this method allows us to handle initial data
which only decays like $r^{-1-\epsilon}$ at infinity.\\

The method we employ here is not based in any way on the 
uniform properties of solutions to the system \eqref{generic_system}.
Instead, our point of departure will be the following simple
observation: Let $\phi^I$ be a given solution to the system
\eqref{generic_ym}. Then is is easy to see that if one performs the
scale transformation:
\begin{equation}
	\phi^I_\lambda(t,x) \ = \ \lambda\, \phi^I(\lambda t,\lambda x)
	\ , \label{ym_scaling}
\end{equation}
the resulting function $\phi^I_\lambda$ is also a solution
to the set of equations \eqref{generic_ym}. This is just a reflection
of the fact that the equations \eqref{generic_ym} are homogeneous.
Suppose now that one could produce a Banach space $B$ which is
dimensionless with respect to the scale transformation 
\eqref{ym_scaling} at time $t=0$. That is, one has the identity:
\begin{equation}
	\lp{\big(\phi^I(0),\partial_t\phi^I\, (0)\big)}{B} \ = \ 
	\lp{\big(\phi^I_\lambda(0),\lambda (\partial_t\phi^I)_\lambda
	\, (0)\big)}{B} \ . \label{initial_phi_norm}
\end{equation}
Suppose that furthermore, one had an existence theorem which said
that any set of initial data for which \eqref{initial_phi_norm}
is small enough, there is a local in time solution to 
\eqref{generic_ym} with this initial data. 
Then, by simply re-scaling, such a 
local existence theorem would necessarily be global in time. Because 
the equations we are considering are hyperbolic, it is natural
to look for a $B$ which is an energy type space. A simple calculation
shows that on $\mathbb{R}^4$, the Sobolev space which is
scale invariant with respect to \eqref{ym_scaling} (at $t=0$)
is the energy space $\dot{H}^1$. However, it is not at all 
unreasonable to expect that 
such a space is far too weak to control solutions to
\eqref{generic_ym} locally in time. For example, $\dot{H}^1$
in $4$ spatial dimensions is a whole
derivative (and then some) away from controlling $L^\infty$.
In fact, one can see immediately from looking at the first
non-trivial Picard iterate to \eqref{generic_ym} 
that one starts to loose regularity as soon as the
initial data is rougher than $H^{1 + \frac{1}{4}}$ (see \cite{KSel_bilin}).
Furthermore, by an adaptation of the $(3+1)$ dimensional 
counterexamples of Lindblad \cite{Lind_counter}, one should be able to show 
that certain instances of 
the equations \eqref{generic_ym} are ill-posed in the
Sobolev spaces $H^s$ when $s < 1 + \frac{1}{4}$. This is 
in stark contrast to the situation in $(5+1)$ and higher dimensions,
where one can come arbitrarily close to the scale invariant
Sobolev space $H^\frac{n-2}{2}$ \cite{Tataru}, and can in fact recover 
local existence in the scale invariant Besov space 
$\dot{B}^{\frac{n-2}{2},1}$ in $(6+1)$ and higher dimensions 
\cite{STR_generic1}.\\

The reason why the low dimensional setting
is more difficult to control than the higher dimensional
regime it that ``parallel'' interactions in the non-linearity
on the right hand side of \eqref{generic_ym} become stronger
and stronger as the dimension decreases. Closely related to this
is the range of validity of the so called Strichartz estimates. 
Specifically, in $(4+1)$ dimensions, one looses the $L^2(L^4)$
Strichartz estimate which clearly plays a major role 
via Duhamel's principle in the
well posedness theory equations with quadratic type interactions (that
is, one looks to put the non-linearity in $L^1(L^2)$). For an important
class of equations with special structure in the non-linearity,
this interaction of parallel waves is largely destroyed, and one
can gain the needed improvement over the $H^{1 + \frac{1}{4}}$
barrier to come arbitrarily close to the scaling. For
example, this was accomplished by Klainerman-Tataru in
\cite{KTat} for the Yang-Mills
equations \eqref{2nd_order_A} with the Coulomb gauge enforced. 
Going even further in this direction, it
should be possible to combine the Besov space technique of 
\cite{STR_generic1} with the compound null structure\footnote{
A close inspection of the proof in \cite{KTat} will show that this
is needed to get around the failure of certain end-point bilinear
$L^1(L^\infty)$ estimates. Specifically, compounding the
non-linearity of \eqref{2nd_order_A} in the Coulomb gauge, and taking into 
account various cancellations due to the null structures
present, one arrives at a set of equations that morally looks like
$\Box \phi = \Delta^{-1}(\phi\nabla\phi)\cdot\nabla\phi$. This
equation looks a lot like wave-maps, except that the weights
are distributed in a more unfavorable fashion. In particular,
while it is true that one can get $\Delta^{-1}(\phi\nabla\phi)
\hookrightarrow L^1(L^\infty)$ at fixed frequency, there is no
room left to add over the low frequencies in a 
$High\times High \Rightarrow Low$
interaction. Therefore, even if one assumes a Besov structure 
for the $\phi$, there is not enough room to close. However, this is 
exactly the bad frequency interaction which is eliminated by the 
tri-linear null structure of \cite{MacStr_MKG}.} discussed
in \cite{MacStr_MKG} to push the global well-posedness theory
of these (Coulomb gauge) equations to the scale invariant
Besov space $\dot{B}^{1,1}$. Finally it is conjectured, and a major
open problem of this subject, that by either working
with the curvature \eqref{YM_curv} directly or by making use of the
Coulomb gauge restriction of the equations \eqref{2nd_order_A},
the equations \eqref{YM_eqs} are well posed in the
scale invariant Sobolev space $\dot{H}^1$. \\

However, our interest here is in the Lorentz gauge equations
\eqref{2nd_order_A_lorentz}, and more generally equations which
are generically of the type \eqref{generic_ym}. An inspection
of the non-linearity in  \eqref{2nd_order_A_lorentz} reveals
that it does not seem to contain the special ``null structure''
of the non-linearity of \eqref{2nd_order_A} in the Coulomb gauge
(at least at the bilinear level), and it is
a tentative conjecture that these specific 
equations are in fact ill-posed
for regularities less than $H^{1+\frac{1}{4}}$. This brings
into question whether one can prove scale invariant global existence
in the spirit of \cite{STR_generic1}. It is clear from the above
discussion that any modification to that theory will need to
go away from translation invariant spaces. That is, one is led
to look for a theory which includes the low regularity micro-local
techniques of \cite{STR_generic1}, but somehow makes crucial
use of the weighted vector-field from 
\cite{Kl_uniform}. One idea is to understand how 
the presence of homogeneous weighted derivatives effects the range 
of validity of the Strichartz estimates. Because the main
obstacle to improved estimates of this type is the presence
of waves which are highly concentrated along a given null direction,
it is natural to expect that the rotation generators:
\begin{equation}
	\Omega_{ij} \ = \ x_i\partial_j - x_j \partial_i 
	\ , \label{rotation_gen}
\end{equation}
play a distinguished role because they penalize such objects.
This indeed turns out to be the case, and one gains a significant
improvement at the level of both linear and bilinear estimates
as was discussed in 
\cite{St_str}. This observation will form the basis for the
first main ingredient of the approach we take here, which
is to better control the linear theory. 
At the non-linear level, one would expect that 
the rotations \eqref{rotation_gen} also play a major role 
because they would help to eliminate parallel interactions
coming in the right hand side of \eqref{generic_ym}. In other words,
one would hope that in some sense the rotations \eqref{rotation_gen}
could substitute for the null-structures
one makes use of in the Coulomb gauge. Again, this turns out to be 
the case and will form the second main pillar of the approach we take
here which is to build function spaces that take into account ``angular
concentration'' phenomena. What we will do is prove
prove the following theorem:\\

\begin{thm}[Global well posedness for the system \eqref{generic_ym}]
\label{GWP_theorem}
        For the generic system of non--linear wave equations 
	\eqref{generic_ym} on $(4+1)$ dimensional Minkowski space,
        there exists constants $0 < \epsilon_0,C$ such that if
        \begin{equation}\label{initial_smallness}
                \lp{(f^I,g^I)}{\dot B_\Omega^{1,1}\times\dot B_\Omega^{0,1}}
                \ \ \leqslant \ \ \epsilon_0 \ ,
        \end{equation}
	where $\dot B_\Omega^{1,1}$ is the Banach space with norm:
	\begin{equation}
	        \lp{h}{\dot B^{1,1}_\Omega} \ = \ 
		\lp{h}{\dot{B}^{1,1}} + \sum_{i < j}
		\lp{\Omega_{ij} h}{\dot{B}^{1,1}} \ , \label{sph_besov_norm}
	\end{equation}
	and likewise for $B_\Omega^{0,1}$,
        then there exits a global solution $\psi^I$ to the system
	\eqref{generic_ym} with initial data $(f^I,g^I)$
        which satisfies the stability condition:
        \begin{equation}\label{cont_cond}
                \lp{\psi^I}{C(\dot B^{1,1}_\Omega)\cap
                C^{(1)}(\dot B^{0,1}_\Omega)} \ \
                \leqslant \ \ C\lp{(f^I,g^I)}{\dot B^{1,1}_\Omega\times
                \dot B^{0,1}_\Omega} \ .
        \end{equation}
	In particular, there is no energy growth of the solution
	to \eqref{generic_ym}.
        The solution $\psi^I$ is unique and depends smoothly on the
	initial data in the following sense: There
        exists a sequence of smooth functions $(f^I_N,g^I_N)$ such that:
        \begin{equation}
                \lim_{N\to \infty}
                \lp{(f^I,g^I) - (f^I_N,g^I_N)}{\dot B^{1,1}_\Omega
		\times\dot B^{0,1}_\Omega}
                 \ =  \ 0 \ . \notag
        \end{equation}
        For this sequence of functions, there exists a sequence of unique
        smooth global solutions $\psi^I_N$ of \eqref{generic_ym}
        with this initial data. Furthermore, the $\psi^I_N$ converge to
        $\psi^I$ as follows:
        \begin{equation}
                \lim_{N\to \infty}
                \lp{\psi^I - \psi^I_N}{ C(\dot{B}^{1,1}_\Omega)\cap
                C^{(1)}(\dot{B}^{0,1}_\Omega) }  \ = \ 0 \ . \notag
        \end{equation}
        Also, $\psi^I$ is the only solution which may be obtained as a
        limit (in the above sense) of solutions to \eqref{generic_system}
        with regularizations of ${(f^I,g^I)}$ as initial data.
        Finally, $\psi^I$ retains any extra smoothness inherent in the initial
        data. That is, if $(f^I,g^I)$ also has finite 
	$\dot{H}^{s}_\Omega\times
        \dot{H}^{s-1}_\Omega$ norm, for $1 < s$, then so does
        $\psi^I$ at fixed time and one has the following estimate:
        \begin{equation}\label{smooth_cont_cond}
                \lp{\psi^I}{C(\dot{H}^s_\Omega)\cap C^{(1)}
		(\dot{H}^{s-1}_\Omega)} \ \
                \leqslant \ \ C\lp{(f^I,g^I)}{\dot{H}^s_\Omega
		\times \dot{H}^{s-1}_\Omega} \ .
        \end{equation}
\end{thm}\ret\ret

\noindent In a straightforward way, our estimates also 
address the issue of the asymptotic freedom of the system 
\eqref{generic_ym}. As an immediate corollary of our approach,
we have that:\\

\begin{thm}\label{scattering_result}
        Using the same notation as above, for our solution $\psi^I$ to the
	system \eqref{generic_system} with initial data $(f^I,g^I)$, 
	there exists
        data sets $({f^I}^\pm,{g^I}^\pm)$, such that if
        ${\psi^I}^\pm$ is the solution to the homogeneous wave equation,
	$\Box{\psi^I}^\pm = 0$,  with this initial data, 
	the following asymptotics hold:
        \begin{align}
                \lim_{t \rightarrow \infty}
                \lp{{\psi^I}^+ - \psi^I}{\dot B_\Omega^{1,1}\cap\partial_t
                \dot B_\Omega^{0,1}}
                \ &= \ 0 \ , \label{scattering_result1} \\
                 \lim_{t \rightarrow -\infty}
                \lp{{\psi^I}^- - \psi^I}{\dot B_\Omega^{1,1}\cap\partial_t
                \dot B_\Omega^{0,1}}
                \ &= \ 0 \ . \label{scattering_result2}
        \end{align}
        Furthermore, the scattering operator retains any additional regularity
        inherent in the initial data. That is, if $(f^I,g^I)$ has finite
        $\dot H_\Omega^s\times\dot H_\Omega^{s-1}$ norm, then so does 
	$({f^I}^\pm,{g^I}^\pm)$, and the following asymptotics hold:
        \begin{align}
                \lim_{t \rightarrow \infty}
                \lp{{\psi^I}^+ - \psi^I}{\dot H_\Omega^s\cap\partial_t 
		\dot H_\Omega^{s-1}}
                \ &= \ 0 \ , \label{smooth_scattering_result1} \\
                 \lim_{t \rightarrow -\infty}
                \lp{{\psi^I}^- - \psi^I}{\dot H_\Omega^s\cap\partial_t 
		\dot H_\Omega^{s-1}}
                \ &= \ 0 \ . \label{amooth_scattering_result2}
        \end{align}
\end{thm}\ret

\begin{rem}
In the statement of the generic system \eqref{generic_ym} and
in proof of Theorem \ref{GWP_theorem} we have ignored the cubic
type interactions $(\phi^I)^3$ which appear on the left hand
side of \eqref{2nd_order_A_lorentz}. Notice that these terms
respect the scaling \eqref{ym_scaling}. It turns out that they
are trivial to treat
in the spaces we use here by taking a product of the
$L^3(L^6)$ Strichartz estimate which is available
in $(4+1)$ dimension. The only real issue is to make sure that
one recovers the Besov structure for $High\times High \Rightarrow Low$ 
frequency
interactions, but this is again a triviality due to the room in the 
$L^3(L^6)$ estimate (one does not even have to use bilinear
estimates to do this). 
\end{rem}\ret

\begin{rem}
For convenience, we have chosen to work here with the spaces
\eqref{sph_besov_norm} which involve a whole angular (momentum)
derivative. As the reader will see shortly, there is much room
in the dyadic estimates of our proof. Specifically, it should be
possible to prove our theorem with the use of only a little
more then $\frac{1}{2}$ an angular derivative. However, this would
force one to work out $L^\infty$ paraproducts in the angular variable,
which would bring another layer of technical complications that we
have chosen to avoid. However, this still 
leaves an interesting gap because
based on the local theory one would expect that, say for compactly
supported initial data, there is global regularity for
small $H^{1 + \frac{1}{4} + \epsilon}$ norm. Therefore, in some sense,
our estimates seem to fall $\frac{1}{4}$ a derivative short of
the optimal level. Perhaps this gap can be eliminated by somehow
incorporating (fractional powers of) the other invariant 
vector-fields. In particular, the boosts 
$\Omega_{0i} = t\partial_i + x_i \partial_t$. We will say no more
of this here.
\end{rem}\ret

\begin{rem}
We have not included here a specific discussion of the first
set of model equations \eqref{generic_wm}. It turns out that these
are a bit easier to treat than the equations \eqref{generic_ym}.
In other words, the difference between \eqref{generic_wm}
and \eqref{generic_ym} which can be seen at the level of decay
can also be seen at the micro-local level. Specifically, for the
equations \eqref{generic_wm}, the estimate \eqref{HH_L1Linf_with_Omega} 
below would
be much easier to prove because it would not need the bilinear
estimates \eqref{impr_ang_str_est}. 
It should be noted however, that the somewhat
more involved version of the estimate \eqref{HH_L1Linf_with_Omega} 
which we use here
can also directly be used in the proof of the well-posedness
of equations of type \eqref{generic_wm}.
\end{rem}

\ret\ret

\section{Notation and preliminary setup}

For quantities $A$ and $B$, we denote by $A \lesssim B$ to mean
that $A \leqslant C\cdot B$ for some large constant $C$. The constant
$C$ may change from line to line, but will always remain
fixed for any given instance where this notation appears.
Likewise we use the notation $A\sim B$ to mean that
$\frac{1}{C}\cdot B \leqslant A \leqslant C\cdot B$. We also
use the notation $A \ll B$ to mean that $A \leqslant
\frac{1}{C}\cdot B$ for some
large constant $C$. This is the notation we will use throughout
the paper to break down quantities into the standard cases:
$A\sim B$, or $A \ll B$, or $B \ll A$; and
$A \lesssim B$, or $B \ll A$, without ever discussing
which constants we are using. We will also employ the following
notation to indicate arbitrarily small adjustments to a given
numerical value:\ For a given constant $A$, we write $A+$
(resp. $A-$) to mean that for any sufficiently small $0 < \epsilon$,
we may replace $A+$ by $A+\epsilon$ (resp. $A-\epsilon$) on the
line where it occurs and still have a true estimate. However,
we \emph{do not} assume any uniformity in this notation. That is, any
implicit constants which appear in conjunction with $A\pm$ may 
depend on $\epsilon$. An example of this is the $L^\infty$  
Sobolev estimate:
\begin{equation}
	\lp{f}{L^\infty(\mathbb{R}^4)} \ \lesssim \
	\lp{f}{H^{2+}(\mathbb{R}^4)} \ . \notag
\end{equation}
Also, if two separate occurrences of the $A\pm$ notation appear on the 
same line, we will not assume that the same $\epsilon$ is being
used for each separate occurrence.\\

For a given function of two variables
$(t,x)\in \bf{R}\times\bf{R}^4$ we write the spatial
and space--time Fourier transform as:
\begin{align}
        \widehat{u}(t,\xi) &= \int e^{-2\pi i \xi\cdot x} \, u(t,x) \ dx
        \ , \notag \\
        \widetilde{u}(\tau,\xi) &= \int e^{-2\pi i (\tau t + \xi\cdot x)}
        \ u(t,x) \ dt  dx \ . \notag
\end{align}
respectively. Because we are not assuming any extra structure in the 
non-linearity of \eqref{generic_ym}, we will work almost exclusively
with the space-time Fourier transform.\\

For a given function $f$ of the spatial variable only, we denote by:
\begin{align}
	e^{it\sqrt{-\Delta}}f\, (x) \ &= \ 
	\int e^{\pi i(t|\xi| +  \xi\cdot x)} \, \hat{f}(\xi) \ d\xi
        \ , \notag \\
	e^{-it\sqrt{-\Delta}}f\, (x) \ &= \ 
	\int e^{\pi i(-t|\xi| +  \xi\cdot x)} \, \hat{f}(\xi) \ d\xi
        \ , \notag
\end{align}
the forward and backward wave propagation of $f$.\\

Let $E$ denote any fundamental solution to the homogeneous wave equation:
i.e., one has the formula $\Box E = \delta$. We define the standard
Cauchy parametrix for the wave equation via the rule:
\begin{equation}
        \frac{1}{\Box} F \ = \ E*F - W (E*F) \ . \notag
\end{equation}
Here and in the sequel, for any test function $H$ 
we use the notation $W(H)$ to
denote the solution to the homogeneous wave equation with initial data
$\left(H(0),\partial_t H\, (0)\right)$.
Explicitly, one has the identity:
\begin{equation}\label{Duhamel_integral}
        \widehat{\frac{1}{\Box}F}\, (t,\xi) \
        = \ - \int_0^t \frac{\sin\left(2\pi
        |\xi|(t-s)\right)}{2\pi |\xi|} \widehat{F}(s,\xi) \ ds \ .
\end{equation}\\
 
For any function $F$ which is supported away from the light cone in
Fourier space, we shall use the following notation for division by the
symbol of the wave equation:
\begin{equation}
        \frac{1}{\varXi} F \ = \
        E*F \ . \notag
\end{equation}
Of course, the definition of $\frac{1}{\varXi}$ does not depend on $E$
so long as for $F$ is supported away from the light cone;
and for us that will always be the case when we use this notation.
Explicitly, one has the formula:
\begin{equation}
        \mathcal{F}\left[ \frac{1}{\varXi} F
        \right](\tau,\xi) \
        = \ \frac{1}{4\pi^2(\tau^2 - |\xi|^2)}
        \widetilde{F}(\tau,\xi) \ . \notag
\end{equation}\\

Next, we record here some basic results from spherical harmonic
analysis. For more details on this material, see the companion
paper to this work \cite{St_str}. The first order of business
concerns defining fractional powers of the spherical Laplacean:
\begin{equation}
	\Delta_{sph} \ = \ \sum_{i<j} \Omega_{ij}^2 \ . \notag
\end{equation}
As is well known, this can be done via spectral resolution, and we
write:
\begin{equation}
	|\Omega|^s \ =\ (-\Delta_{sph})^\frac{s}{2} \ . \label{sph_lap_s}
\end{equation}
The operator \eqref{sph_lap_s} kills off the spherically symmetric
part of any function it is applied to. Because of this, we will 
employ the following ``inhomogeneous'' version of this operator:
\begin{equation}
	\langle \Omega \rangle^s f \ =\ F_0 + |\Omega|^s f \ , \notag
\end{equation}
where $f$ is a function of the spatial variable, 
and $f_0$ denotes the spherically symmetric
part of $f$. That is:
\begin{equation}
        f_0(r) \ = \ \frac{1}{|\mathbb{S}^3|} \int_{\mathbb{S}^3}
	\ f(r\omega)\ d\omega \ . \notag
\end{equation}
A key property of the operators $\langle\Omega\rangle^s$
is that they commute with the spatial (and thus space-time) Fourier
transform:
\begin{equation}
	\widehat{\langle\Omega\rangle^s f} \ = \ \langle\Omega\rangle^s
	\widehat{f} \ . \notag
\end{equation}
Also, we have the following equivalence of Sobolev type norms\footnote{
The $\Omega$ subscript in conjunction with numerical superscripts,
e.g. the $s$ in $H^s_\Omega$, will always denote angular derivatives
in this section.
In other places in the paper the $\Omega$ subscript will always mean
\emph{one} angular derivative, while the superscripts will denote
translation invariant derivatives. An example of this is the notation 
$\dot{B}^{1,1}_\Omega$ introduced on line \eqref{sph_besov_norm}.}
involving the unit power $\langle\Omega\rangle$:
\begin{equation}
	\lp{f}{H^1_\Omega(\mathbb{R}^4)}^2 \ = \ 
	\lp{\langle \Omega \rangle f}{L^2(\mathbb{R}^4)}^2 \ = \
	\lp{f_0}{L^2(\mathbb{R}^4)}^2 + 
	\sum_{i<j}  \lp{\Omega_{ij} f}{L^2(\mathbb{R}^4)}^2 \ . 
	\label{norm_equiv}
\end{equation}
Because all of the norms we build here will be based on the unit powers
$\langle\Omega\rangle$,
we will by abuse of notation replace any instance of a single $\Omega_{ij}$
with the operator $\langle \Omega \rangle$, and although it is not
strictly true, we will assume that there is the \emph{point-wise} Leibniz
rule:
\begin{equation}
	\langle\Omega \rangle (fg) \ = \ \langle\Omega \rangle f\cdot g +
	f\cdot\langle\Omega \rangle g \ . \notag 
\end{equation}
This will be a great convenience to us because some of the norms
we define below are $L^\infty$ based (in particular 
\eqref{fixed_Z_norm} and \eqref{ang_Z_norm}), 
where it would be difficult to define
paraproducts for fractional powers $\langle \Omega \rangle^s$.\\

Finally, we record here two basic results which follow from the 
Littlewood-Paley theorem for the sphere, in conjunction with 
interpolation in weighted
spaces of the type $L^p(\ell^2_s)$ (again, see \cite{St_str} for details):\\

\begin{prop}[Sobolev embedding on the sphere]\label{sph_sob_prop}
If $F$ is a test function on the unit sphere $\mathbb{S}^3\subset
\mathbb{R}^4$, then the following
estimate holds for $2 \leqslant p< \infty$:
\begin{equation}
	\lp{F}{L^p(\mathbb{S}^3)} \ \lesssim \ 
	\lp{ \langle\Omega\rangle^{3(\frac{1}{2} - \frac{1}{p})}
	F}{L^2(\mathbb{S}^3)} \ , \notag
\end{equation}
where the implicit constants depend only on $p$. In particular, if $f$ is
a test function on $\mathbb{R}^4$ and $s < \frac{3}{2}$, then one has:
\begin{equation}
	\int_0^\infty \ \lp{f(r)}{L^\frac{6}{3 - 2s}( \mathbb{S}^3)}^2\ r^3 \, dr
	 \ \lesssim \ 
	\int_0^\infty \ \lp{\langle\Omega\rangle^s f(r)}
	{L^2( \mathbb{S}^3)}^2\ r^3 \, dr \ = \ 
	\lp{f}{H^s_\Omega
	(\mathbb{R}^4)} \ . \label{sph_sob_embed}
\end{equation}
\end{prop}\ret

\begin{prop}[Interpolation of spherical Sobolev spaces]\label{sph_interp_prop}
Let $W_\Omega^{s,p}$ denote the norm:
\begin{equation}
	\lp{f}{W_\Omega^{s,p}(\mathbb{R}^4) } \ = \ \lp{\langle\Omega\rangle^s
	f}{L^p(\mathbb{R}^4) } \ , \notag
\end{equation}
for functions $f$ of the spatial variable only. 
Then for $1 < p_1 , p_2 < \infty$ 
one has the following interpolation spaces:
\begin{equation}
	\left(  W_\Omega^{s_1,p_1} , W_\Omega^{s_2,p_2}    \right)_t
	\ = \ W_\Omega^{s,p} \ , \label{sph_interp_result}
\end{equation}
where $s = (1-t)s_1 + ts_2$ and $\frac{1}{p} = \frac{(1-t)}{p_1}
+ \frac{t}{p_2}$.
\end{prop}

\ret\ret

\section{Strichartz Estimates}

We list here the space--time estimates for the homogeneous wave
equations which form the foundation for our proof of Theorem
\ref{GWP_theorem}. As we have mentioned before, all of these are
of ``Strichartz type''. The first group of estimates we will use
are just the classical Strichartz estimates for the wave equation
which we state for the case of $(4+1)$ dimensions:\\

\begin{prop}[Frequency localized
``classical'' Strichartz estimates on $\mathbb{R}^{(4+1)}$
including endpoints (see \cite{KT_str})]\label{classical_str_th}
Let $n=4$ be the number of spatial dimensions, and let
$\sigma=\frac{3}{2}$ be the corresponding Strichartz admissible
exponent. If $f$ is any function of the spatial
variable only, denote by $f_1 = P_1 f$ its unit frequency
projection (see \eqref{s_cutoff}).  Then one has the following
family estimates for $2\leqslant q$:
\begin{equation}
        \lp{e^{it\sqrt{-\Delta}}f_1}{L^q(L^r)} \ \lesssim \
        \lp{f_1}{L^2} \ , \label{classical_str_est}
\end{equation}
where $\frac{1}{q} + \frac{\sigma}{r} \leqslant \frac{\sigma}{2}$. 
\end{prop}\ret

As was discussed in the introduction, the estimates 
\eqref{classical_str_est} alone are not strong enough to close a global
iteration argument for non-linear wave equations of the form 
\eqref{generic_system}. What is needed is an improvement 
of the range of admissible $(q,r)$ indices on the left hand
side of \eqref{classical_str_est}. It is well known that this cannot
be accomplished within the context of translation invariant
smoothness assumptions on the initial data (see \cite{St_str}). 
However, incorporating extra weighted smoothness assumptions
for the angular variable provides the needed mechanism to overcome this 
obstacle. The corresponding estimates are:\\

\begin{prop}[Frequency localized
Strichartz estimates for angularly regular initial
data (see \cite{St_str})]\label{ang_str_th}
Let $n=4$ be the number of spatial dimensions, and let
$\sigma_\Omega=3$ denote the four dimensional angular Strichartz
admissible exponent. Let $f_1$ be a unit frequency function of the
spacial variable only (as above). Then for indices $(q,r)$ such that
$\frac{1}{q} + \frac{\sigma}{r} \geqslant \frac{\sigma}{2}$ and 
$\frac{1}{q} + \frac{\sigma_\Omega}{r} <  \frac{\sigma_\Omega}{2}$, 
and for every $0 < \epsilon$, there is a $C_\epsilon$ which depends
only on $\epsilon$ such that the following estimates hold:
\begin{equation}
        \lp{e^{it\sqrt{-\Delta}}f_1}{L^q(L^r)} \ \lesssim \
        C_\epsilon\ \lp{\langle \Omega \rangle^s f_1}{L^2} 
	\ , \label{ang_str_est}
\end{equation}
where $s= (1 + \epsilon)(\frac{n-1}{r} + \frac{2}{q}- \frac{n-1}{2})$.
\end{prop}\ret

In practice, only a small subset of the indices $(q,r)$ in the
two propositions listed above will  be of use to us. 
These are $(\infty,2)$, $(2,\infty)$,
$(2,6)$, and $(2,3+)$. To highlight this fact, we list out
the corresponding instances of \eqref{classical_str_est} and 
\eqref{ang_str_est}:\\

\begin{align}\hline \hline \notag \\
	\lp{e^{it\sqrt{-\Delta}}f_1
	}{L^\infty(L^2)} \ &\lesssim \ \lp{f_1}{L^2} \ , 
	\label{st_energy_str_est} \\
	\lp{e^{it\sqrt{-\Delta}}f_1
	}{L^2(L^\infty)} \ &\lesssim \  \lp{f_1}{L^2}\ , 
	\label{st_L2Linf_str_est} \\
	\lp{e^{it\sqrt{-\Delta}}f_1
	}{L^2(L^6)} \ &\lesssim \  \lp{f_1}{L^2} \ , 
	\label{st_L2L6_str_est} \\
	\lp{\langle\Omega\rangle^\frac{1}{2}\,
	e^{it\sqrt{-\Delta}}f_1
	}{L^2(L^{3+})} \ &\lesssim \  \lp{\langle \Omega\rangle f_1}{L^2}\ , 
	\label{st_L2L3_str_est}
	\\ \notag \\ \hline \hline \notag 
\end{align}\\

In our proof of Theorem \ref{GWP_theorem}, we will need more than
just the linear estimates 
\eqref{classical_str_est}--\eqref{ang_str_est}. This is a common feature
of lower dimensional problems, and is necessitated by the presence
of certain bad $High\times High \Rightarrow Low$
frequency interactions. The standard device for dealing with this 
problem is the use of bilinear Strichartz estimates. 
We will use here T. Tao's fine--course scale idea for dealing with 
these (see \cite{TKR} and \cite{Tao_semi}). The basic idea is to fix
a scale, say $\frac{1}{\mu}$ for $\mu\ll 1$, and then decompose the domain
of the spatial variable into cubes with side lengths $\sim\frac{1}{\mu}$. 
Then, one replaces the usual $L^r$ norm in the spatial variable with
$\ell^r(L^2)$, where the $L^2$ norm is taken on the ``fine'' scale
of each individual cube, while the $\ell^r$ norm represents the ``coarse''
scale which is summation over all cubes. One reason this method
is so powerful, is that it allows one to use the bilinear construction process
directly in an iteration procedure, where  resorting to the canned estimates
that this method ultimately provides may be unduly burdensome. This
is crucial when dealing with eccentric multipliers as we do here.
Therefore, we will only state the two scale estimates themselves,
without mentioning the various bilinear estimates which follow as a
corollary. For point of reference, we point out here that these 
estimates will only be used in the proof of estimate 
\eqref{HH_L1Linf_with_Omega} below. We begin by stating the
classical two--scale estimates:\\

\begin{prop}[Frequency localized two--scale 
Strichartz estimates (\cite{TKR} ,\cite{Tao_semi})]\label{imp_str_th}
Let $n=4$ be the number of spatial dimensions. Let
$0 < \mu \lesssim 1$ be a given parameter. Let $\{Q_\alpha\}$ be a partition
of $\mathbb{R}^n$
into cubes of side length $\sim\frac{1}{\mu}$.
Then if $f_1$ is a unit frequency function of the spatial
variable only, the following estimates hold:
\begin{equation}
        \lp{\big(\sum_\alpha \lp{
	e^{it\sqrt{-\Delta}}f_1}
        {L^2(Q_\alpha)}^r\big)^\frac{1}{r}}{L_t^2}
        \ \lesssim \ \mu^{-1} \,
        \lp{f_1}{L^2}
        \ , \label{impr_str_est}
\end{equation}
where $6 \leqslant r$.
\end{prop}\ret

\noindent Next, we state the improvement to
\eqref{impr_str_est} which incorporates angular regularity:\\

\begin{prop}[Frequency localized two--scale Strichartz estimates for
angularly regular data; endpoint case (\cite{St_str})]\label{imp_ang_str_th}
Let $n=4$ be the number of spatial dimensions
and let $f_1$ be a unit frequency function of the
spatial variable only. Let
$0 < \mu \lesssim 1$ be given, and let $\{Q_\alpha\}$ be a partition
of $\mathbb{R}^n$
into cubes of side length $\sim\frac{1}{\mu}$. Then for
every $0 < \epsilon$, there is a $C_\epsilon$ and 
$ 3 < r_\epsilon$ depending on $\epsilon$,
such that $r_\epsilon \to 3$ as $\epsilon\to 0$
such that the following estimate holds:
\begin{equation}
        \lp{\big(\sum_\alpha \lp{
	e^{it\sqrt{-\Delta}}f_1 }
        {L^2(Q_\alpha)}^{r_\epsilon}\big)^\frac{1}{r_\epsilon}}{L_t^2}
        \ \lesssim \ C_\epsilon \, \mu^{-\frac{1}{2} -2\epsilon} \,
        \lp{\langle \Omega \rangle^{\frac{1}{2} + \epsilon} \, f_1 }{L^2}
        \ . \label{impr_ang_str_est}
\end{equation}
\end{prop}\ret

\noindent In practice, we will only need a single 
instance of both \eqref{impr_str_est}
and \eqref{impr_ang_str_est}. We state these in encapsulated form here
for the convenience of the reader:\\

\begin{align}\hline \hline \notag \\
	\lp{\big(\sum_\alpha \lp{
	e^{it\sqrt{-\Delta}}f_1 }
	{L^2(Q_\alpha)}^6\big)^\frac{1}{6}}{L_t^2}
	\ &\lesssim \ \mu^{-1}\, \lp{f_1}{L^2} \ , 
	\label{st_imp_L2L6_str_est} \\
	\lp{\big(\sum_\alpha \lp{
	e^{it\sqrt{-\Delta}}f_1 }
	{L^2(Q_\alpha)}^{3+}\big)^\frac{1}{3+}}{L_t^2}
	\ &\lesssim \ \mu^{-(\frac{1}{2}+)}\, \lp{\langle \Omega \rangle
	f_1}{L^2} \ . 
	\label{st_imp_L2L3_str_est}
	\\ \notag \\ \hline \hline \notag 
\end{align}\\

\noindent Notice that we have added an extra $\frac{1}{2}$ an angular 
derivative to the right hand side of \eqref{st_imp_L2L3_str_est} above.
This is how we will use this estimate in this paper, and is a reflection
of the fact that we have elected to work with integral powers 
of the angular momentum operator $|\Omega|$ in this work.

\ret\ret

\section{Multipliers, Functions Spaces, and Scattering}

In this section, we will set up much of notation to be used in the proof
of Theorem \ref{GWP_theorem} and we will construct the function spaces 
used to iterate the problem \eqref{generic_ym}.
For the most part, the approach taken here is similar to
that of \cite{STR_generic1}, with the simple addition of 
angular derivatives. A notable
exception occurs in the definition of the special $L^1(L^\infty)$ 
``outer block'' norms \eqref{Z_norm} and \eqref{ang_Z_norm}. Because of the 
need to capture extra savings in our
estimates based on angular regularity, these are a bit more involved than
their cousins used in \cite{STR_generic1}. 
We strongly recommend that the reader
first read that paper as a warm up to the present work because it represents
a simplified version of the type of decompositions and
estimate combinations used here. \\

\ret

\subsection{Multipliers and angular restrictions}

Let $\varphi$ be  a smooth bump function (i.e.
supported on the set $|s| \leqslant 2$ such
that $\varphi = 1$ for $|s| \leqslant 1$). In what follows, it will be
a great convenience for us to assume that $\varphi$ may change its exact
form for two separate instances of the symbol $\varphi$ (even if they
occur on the same line). In this way, we may assume without loss of generality
that in addition to being smooth, we also have the idempotence identity
$\varphi^2= \varphi$. We shall use this convention for all the cutoff functions
we introduce in the sequel.\\

For $\lambda\in 2^\mathbb{Z}$,
we denote the dyadic scaling of $\varphi$ by
$\varphi_\lambda(s) = \varphi(\frac{s}{\lambda})$. The most basic Fourier
localizations we shall use here are with respect to the spatial and
space-time Fourier variable and the distance from the cone in Fourier space. 
Accordingly, we form the Littlewood-Paley type cutoff functions:
\begin{align}
        p_\lambda(\xi) &= \varphi_{2\lambda}(|\xi|)
        - \varphi_{\frac{1}{2}\lambda}(|\xi|) \ , \label{s_cutoff} \\
        s_{\lambda}(\tau,\xi) &= \varphi_{2\lambda}(|(\tau,\xi)|)
        - \varphi_{\frac{1}{2}\lambda}(|(\tau,\xi)|) \ , \label{st_cutoff} \\
        c_d(\tau,\xi) &= \varphi_{2d}(|\tau| - |\xi|) -
        \varphi_{\frac{1}{2} d}(|\tau| - |\xi|) \ . \label{cone_cutoff}
\end{align}
We now denote the corresponding Fourier multiplier operator via the formulas
$\ \ \widetilde{S_\lambda u} = s_\lambda \widetilde{u}\ \ $ and
$\ \ \widetilde{C_d u} = c_d \widetilde{u}\ \ $ respectively. We also use a
multi-subscript notation to denote products of the above operators, e.g.
$\ \ S_{\lambda,d} = S_\lambda C_d\ \ $. We shall use the notation:
\begin{equation}
        S_{\lambda,\bullet \leqslant d} = \sum_{\delta \leqslant d}
        S_{\lambda,\delta} \ , \label{d_cone_nbd}
\end{equation}
to denote cutoff in an $O(d)$ neighborhood of the light cone in Fourier
space. At times it will also be convenient to write
$\ \ S_{\lambda,d \leqslant \bullet} = S_\lambda -
S_{\lambda,\bullet < d}\ \ $. We shall also use the notation
$S^\pm_{\lambda,d}$ etc. to denote the multiplier $S_{\lambda,d}$
cutoff in the half space $\ \ \pm \tau > 0\ \ $.\\

The other type of Fourier localization which will be central to our analysis
will be the restriction of the spatial
angular variable $\omega$, where $\xi = |\xi|\,\omega$. We accomplish this
as follows: For each small parameter $\eta\lesssim 1$, we decompose the
unit sphere in $\mathbb{R}^4$ into angular sectors of size $\sim\eta$ with
bounded overlap independent of $\eta$. We label the corresponding partition
of unity by their angles and write them as $b^\omega_{\eta}$. It is
clear that this construction can be done in such a way that all the
$b^\omega_{\eta}$ are (essentially) rotations of each other. Note
that the multipliers $b^\omega_{\eta}p_\lambda$ essentially
cutoff on parallelepipeds of size 
$\lambda\times(\lambda\eta)\times(\lambda\eta)\times(\lambda\eta)$. Furthermore,
after rotating each of these multipliers onto the $\xi_1$ axis, one has the 
following bounds:
\begin{align}
        |\partial_1^N \, b^\omega_{\eta}p_\lambda
	| &\leqslant C_N
        \lambda^{-N} \ ,
        & |\partial_i^N  \, b^\omega_{\eta}p_\lambda 
	| &\leqslant C_N
        (\lambda\eta)^{-N} \ , \label{block_smoothness}
\end{align}
for $i=2,3,4$. In particular we see that each operator 
$B^\omega_{\eta}P_\lambda$ is given by convolution with an $L^1$ kernel,\\

A major defect of the $S_{\lambda,d}$ multipliers is that they are not uniformly
bounded on most Lebesgue spaces. However, if we first localize them further
onto a block that is directed along the light cone in Fourier space of
dimensions $\lambda\times\sqrt{\lambda d}\times\sqrt{\lambda d}\times
\sqrt{\lambda d}\times d$, then the resulting kernels will be uniformly in 
$L^1$. In the sequel we shall write these special block localizations as:
\begin{align}
        S^\omega_{\lambda,d} &= B^\omega_{(\frac{d}{\lambda})^\frac{1}{2}}
	P_\lambda
        S_{\lambda,d} \ ,
        & S^\omega_{\lambda,\bullet \leqslant d} &=
        B^\omega_{(\frac{d}{\lambda})^\frac{1}{2}} P_\lambda
        S_{\lambda,\bullet \leqslant d} \ . \notag
\end{align}
It is important to note that the above multipliers are cutoffs in the
region of Fourier space where $|\tau| \lesssim |\xi|$.
We now record some useful multiplier bounds, the proofs of which can be found
in \cite{STR_generic1}:\\

\begin{lem}[Multiplier boundedness on Lebesgue 
spaces]\label{generic_mult_lemma}
\ \ \ 
\begin{enumerate}
        \item The following multipliers are given by $L^1$ kernels: \ \ 
        $\lambda^{-1}\nabla S_\lambda$, 
	$B^\omega_{(\frac{d}{\lambda})^\frac{1}{2}}P_\lambda$,  
	$S^\omega_{\lambda,d}$,
        and $(\lambda d)\varXi^{-1}S^\omega_{\lambda,d}$ . In particular,
        all of these are bounded on every mixed Lebesgue space $L^q(L^r)$. 
        \item The following multipliers are bounded on the spaces
        $L^q(L^2)$, for $1 \leqslant q \leqslant \infty$: \ \ 
        $S_{\lambda,d}$, and $S_{\lambda,\bullet \leqslant d}$. 
\end{enumerate}
\end{lem}\ret

In the sequel, we shall also need the following somewhat stronger version of the
boundedness of the multipliers $S^\omega_{\lambda,d}$ and 
$(\lambda d)\varXi^{-1}S^\omega_{\lambda,d}$:

\begin{lem}[Multiplier boundedness on special 
$L^1(L^\infty)$ spaces]\label{spec_mult_lemma}
Let $u$ be a function of space and time, then the following estimates
hold:
\begin{align}
	\int \sup_\omega \lp{\frac{1}{\varXi}S^\omega_{\lambda,d}u\, (t)}
	{L^\infty}\, dt \ &\lesssim \ \frac{1}{\lambda d}
	\int \sup_\omega \lp{S^\omega_{\lambda,d}u\, (t)}
	{L^\infty}\, dt \ , \label{sup_L1Linf_bound} \\
	\int \left( \sum_\omega 
	\lp{\frac{1}{\varXi}S^\omega_{\lambda,d}u\, (t)}{L^\infty}^2
	\right)^\frac{1}{2}\, dt \ &\lesssim \
	\frac{1}{\lambda d} \int \left( \sum_\omega 
	\lp{S^\omega_{\lambda,d}u\, (t)}{L^\infty}^2
	\right)^\frac{1}{2}\, dt \ . \label{sqsum_L1Linf_bound}\\
	\int \sup_\omega \lp{S^\omega_{\lambda,d}u\, (t)}
	{L^\infty}\, dt \ &\lesssim \ 
	\int \sup_\omega \lp{ B^\omega_{(\frac{d}{\lambda})^\frac{1}{2}}
	P_\lambda  S_\lambda u\, (t)}
	{L^\infty}\, dt \ , \label{sup_L1Linf_bound_no_box} \\
	\int \left( \sum_\omega 
	\lp{S^\omega_{\lambda,d}u\, (t)}{L^\infty}^2
	\right)^\frac{1}{2}\, dt \ &\lesssim \
	\int \left( \sum_\omega 
	\lp{ B^\omega_{(\frac{d}{\lambda})^\frac{1}{2}}
	P_\lambda  S_\lambda u\, (t)}{L^\infty}^2
	\right)^\frac{1}{2}\, dt \ . \label{sqsum_L1Linf_bound_no_box}
\end{align}
\end{lem}\ret

\begin{proof}[proof of estimates 
\eqref{sup_L1Linf_bound}--\eqref{sqsum_L1Linf_bound_no_box}]
It suffices to prove the implications \eqref{sup_L1Linf_bound}
and \eqref{sqsum_L1Linf_bound}, as the proofs of
\eqref{sup_L1Linf_bound_no_box} and \eqref{sqsum_L1Linf_bound_no_box}
follow from virtually identical reasoning. Using the idempotence 
relation $S^\omega_{\lambda,d} = S^\omega_{\lambda,d}S^\omega_{\lambda,d}$,
and writing $K^\omega$ for the convolution kernel of the operator
$(\lambda d) \varXi^{-1}S^\omega_{\lambda,d}$, for the estimate
\eqref{sup_L1Linf_bound} we can bound:
\begin{align}
        &(\lambda d) \, \int
        \sup_\omega \lp{\frac{1}{\varXi}S^\omega_{\lambda,d}u\, (t)}
	{L^\infty_x}\, dt \notag \\
	= \
	&\int \sup_\omega \lp{K^\omega * S^\omega_{\lambda,d}u\, (t)}
	{L^\infty_x}\, dt \ , \notag \\
	\lesssim \  
	&\int \left( \sup_\omega \int \lp{K^\omega(t-s)}{L^1_x} 
	\lp{ S^\omega_{\lambda,d}u\, (s)}{L^\infty_x}\,ds\right) 
	dt \ , \notag \\
	\lesssim \
	&\int\int \ \sup_\omega \lp{K^\omega(t-s)}{L^1_x} 
	\ \sup_\omega\lp{ S^\omega_{\lambda,d}u\, (s)}{L^\infty_x}\,ds\, 
	dt \ . \notag 
\end{align}
We now use the fact that all of the $K^\omega$ are essentially spatial
rotations of each other to show that for each fixed $(t-s)$:
\begin{equation}
        \lp{K^\omega(t-s)}{L^1_x} \ \lesssim \ \lp{K^{\omega_0}(t-s)}{L^1_x}
	\ , \notag
\end{equation}
where $\omega_0$ is, say, the angular sector in the direction of the
$\xi_1$ axis in Fourier space. Thus, we have that:
\begin{align}
        &(\lambda d) \, \int
        \sup_\omega \lp{\frac{1}{\varXi}S^\omega_{\lambda,d}u\, (t)}
	{L^\infty_x}\, dt \notag \\
	\lesssim \ &\int\int \ \lp{K^{\omega_0}(t-s)}{L^1_x} 
	\ \sup_\omega\lp{ S^\omega_{\lambda,d}u\, (s)}{L^\infty_x}\,ds\, 
	dt \ . \notag \\
	\lesssim \ &\lp{K^{\omega_0}}{L^1(L^1)}\, \int
	\sup_\omega\lp{ S^\omega_{\lambda,d}u\, (s)}{L^\infty_x}\,ds
	\ . \notag \\
	\lesssim \ &\int
	\sup_\omega\lp{ S^\omega_{\lambda,d}u\, (s)}{L^\infty_x}\,ds \ . \notag
\end{align}
Similarly, for the estimate \eqref{sqsum_L1Linf_bound}, we compute:
\begin{align}
        &(\lambda d) \, \int\left(
        \sum_\omega \lp{\frac{1}{\varXi}S^\omega_{\lambda,d}u\, (t)}
	{L^\infty_x}^2\right)^\frac{1}{2}\, dt \notag \\
	= \
	&\int \left( \sum_\omega \lp{K^\omega * S^\omega_{\lambda,d}u\, (t)}
	{L^\infty_x}^2\right)^\frac{1}{2}\, dt \ , \notag \\
	\lesssim \  
	&\int \left( \sum_\omega \left( \int \lp{K^\omega(t-s)}{L^1_x} 
	\lp{ S^\omega_{\lambda,d}u\, (s)}{L^\infty_x}\, ds \right)^2  \right)
	^\frac{1}{2}\, dt \ , \notag \\
	\lesssim \
	&\int\int \left( \sum_\omega \lp{K^\omega(t-s)}{L^1_x}^2\, 
	\lp{ S^\omega_{\lambda,d}u\, (s)}{L^\infty_x}^2\right)^\frac{1}{2}
	\,ds\, dt \ . \notag \\
	\lesssim \ &\int\int \ \sup_\omega \lp{K^\omega(t-s)}{L^1_x}\left( 
	\sum_\omega
	\lp{ S^\omega_{\lambda,d}u\, (s)}{L^\infty_x}^2\right)^\frac{1}{2}
	\,ds\, dt \ . \notag \\
	\lesssim \ &\lp{K^{\omega_0}}{L^1(L^1)}\, \int\left(
	\sum_\omega\lp{ S^\omega_{\lambda,d}u\, (s)}
	{L^\infty_x}^2\right)^\frac{1}{2}\,ds \ . \notag \\
	\lesssim \ &\int\left(
	\sum_\omega\lp{ S^\omega_{\lambda,d}u\, (s)}
	{L^\infty_x}^2\right)^\frac{1}{2}\,ds\ . \notag
\end{align}
This completes the proof of estimates 
\eqref{sup_L1Linf_bound}--\eqref{sqsum_L1Linf_bound_no_box}.
\end{proof}\ret

\ret

\subsection{Function spaces}

Our next step will be to use the above multipliers to define dyadic versions
of the function spaces we will Picard iterate in. We first 
define those dyadic norms which guarantee that the various
Strichartz estimates of Theorems \ref{classical_str_th}, 
\ref{ang_str_th}, \ref{imp_str_th}, and \ref{imp_ang_str_th} hold. These are:
\begin{align}
        \lp{u}{X^{\frac{1}{2}}_{\lambda,p}}^p &=
        \sum_{ d\in 2^\mathbb{Z}} d^{\frac{p}{2}}
        \lp{S_{\lambda,d} u}{L^2}^p \ , &\hbox{(``classical'' $H^{s,\delta}$)}
        \label{Xsd_norm} \\ 
        \lp{u}{Y_\lambda} &=
        \lambda^{-1}\lp{\Box S_\lambda u}{L^1(L^2)} \ .
        &\hbox{(Duhamel)} \label{duhamel_norm} 
\end{align}
Because these function spaces give no weight to solutions
to the homogeneous wave equation, we will always need to use them in
concert with the fixed frequency energy space 
$S_\lambda\left(L^\infty(L^2)\right)$.
Grouping all of these together, we form our first main dyadic function
space:
\begin{equation}
        \lp{u}{F_\lambda} = \left( X^\frac{1}{2}_{\lambda,1}
        + Y_\lambda \right) \cap S_\lambda
        \left( L^\infty(L^2) \right) \ . \label{fixedF_norm}
\end{equation}
We also define the corresponding norms with angular derivatives added 
as follows:
\begin{align}
        \lp{u}{F_{\Omega,\lambda}} \ = \ 
	\lp{\langle\Omega\rangle u}{F_\lambda} \ . \label{ang_fixedF_norm} 
\end{align}
At times it will also be convenient for us to write:
\begin{align}
        \lp{u}{X^\frac{1}{2}_{\Omega,\lambda,1}} \ &= \
	\lp{\langle\Omega\rangle u}{X^\frac{1}{2}_{\lambda,1}} \ , \notag \\
	\lp{u}{Y_{\Omega,\lambda}} \ &= \ 
	\lp{\langle\Omega\rangle u}{Y_{\lambda}} \ . \notag
\end{align}
A key property of the spaces $F_\lambda$ and $F_{\Omega,\lambda}$
is that their elements can be written as integrals over solutions
to the homogeneous wave equation where the integration involves only
the temporal variable (see \cite{STR_generic1} for details).
This is sometimes referred to as foliation, or the \emph{trace method}. 
Using this technique, the Strichartz estimates from
Theorems \ref{classical_str_th}, \ref{ang_str_th}, \ref{imp_str_th}, and 
\ref{imp_ang_str_th} can be transferred to the $F$ spaces
in a straightforward way, even though some of the estimates involve weighted
angular derivatives. This is one of the main reasons why the 
use of Strichartz estimates
which only involve the rotation vector fields is so crucial 
to the approach taken in this paper. Indeed, if the Strichartz estimates
\eqref{ang_str_est} and \eqref{impr_ang_str_est} 
contained any of the other invariant vector fields from \cite{Kl_uniform}, there
would be no way to transfer them to the $F_{\Omega,\lambda}$ via usual
temporal foliation.
We now list the various foliated instances of estimates
\eqref{classical_str_est}, \eqref{ang_str_est}, \eqref{impr_str_est}, and
\eqref{impr_ang_str_est} which will be used in this paper: \\

\begin{align}\hline \hline \notag \\
	\lp{S_1 u}{L^\infty(L^2)} \ &\lesssim \ \lp{u}{F_1} \ , 
	\label{F_energy_str_est} \\
	\lp{S_1 u}{L^2(L^\infty)} \ &\lesssim \ \lp{u}{F_1} \ , 
	\label{F_L2Linf_str_est} \\
	\lp{S_1 u}{L^2(L^6)} \ &\lesssim \ \lp{u}{F_1} \ , 
	\label{F_L2L6_str_est} \\
	\lp{S_1 \langle\Omega\rangle^\frac{1}{2}\,
	u}{L^2(L^{3+})} \ &\lesssim \ \lp{u}{F_{\Omega,1}} \ , 
	\label{F_L2L3_str_est}\\
	\lp{S_\lambda u}{L^\infty(L^2)} \ &\lesssim \ \lp{u}{F_\lambda} \ , 
	\label{high_F_energy_str_est} \\
	\lp{S_\lambda u}{L^2(L^6)} \ &\lesssim \ \lambda^\frac{5}{6}
	\lp{u}{F_\lambda} \ , 
	\label{high_F_L2L6_str_est} \\
	\lp{S_\lambda \langle\Omega\rangle^\frac{1}{2}\,
	u}{L^2(L^{3+})} \ &\lesssim \ \lambda^{\frac{1}{6}+}
	\lp{u}{F_{\Omega,\lambda}} \ , 
	\label{high_F_L2L3_str_est} \\
	\lp{\big(\sum_\alpha \lp{S_1 u\, (t)}
	{L^2(Q_\alpha)}^6\big)^\frac{1}{6}}{L_t^2}
	\ &\lesssim \ \mu^{-1}\, \lp{u}{F_1} \ , 
	\label{imp_L2L6_str_est} \\
	\lp{\big(\sum_\alpha \lp{S_1 u\, (t)}
	{L^2(Q_\alpha)}^{3+}\big)^\frac{1}{3+}}{L_t^2}
	\ &\lesssim \ \mu^{-(\frac{1}{2}+)}\, \lp{u}{F_{\Omega,1}} \ . 
	\label{F_imp_L2L3_str_est}
	\\ \notag \\ \hline \hline \notag 
\end{align}\\

An important property of the $Y_\lambda$ spaces is that they are nearly
contained in the $X^\frac{1}{2}_{\lambda,1}$ spaces. More precisely, this is 
true at a fixed dyadic distance from the light cone in Fourier space. To see 
this, note that by duality and the estimate \eqref{F_energy_str_est},
we have the inclusions:
\begin{equation}
        \lambda \varXi^{-1} L^1(L^2)
        \ \subseteq \ \lambda \varXi^{-1}
        X^{-\frac{1}{2}}_{\lambda,\infty} \ \subseteq \
        X^\frac{1}{2}_{\lambda,\infty} \ . \notag
\end{equation}
Using the above embedding at fixed distance from the cone, and by dyadic 
summing, shows that we have the following estimates:\\

\begin{align}\hline \hline \notag \\
        d^\frac{1}{2}\ \lp{S_{\lambda,d}u}{L^2(L^2)} \ \lesssim \
	\lp{u}{F_\lambda} \ , \label{L2_F_est}\\
	d^\frac{1}{2}\ \lp{S_{\lambda,d \leqslant 
	\bullet}u}{L^2(L^2)} \ \lesssim \
	\lp{u}{F_\lambda} \ . \label{sum_L2_F_est}
        \\ \notag \\ \hline \hline \notag 
\end{align}\\

It is well known, the norms \eqref{Xsd_norm} and \eqref{duhamel_norm}
are not strong enough to iterate wave equations which contain
derivatives in their nonlinearities like \eqref{generic_ym}. This
is due to the presence of a very specific $Low\times High$ frequency 
interaction in the term $\phi\,\nabla\phi$. What is needed to circumvent 
this problem
is to add some extra $L^1(L^\infty)$ norms to the $F$ spaces. This idea 
originally goes back to the work of Klainerman--Machedon \cite{KMQij}, and
was later used to its full extent in Tataru \cite{Tataru}. 
In our previous work \cite{STR_generic1},  we used a slight innovation
on the norms in \cite{Tataru} which allowed one to work in a scale
invariant setting. Here it will be necessary for us 
to use somewhat more technical versions of those norms, in part
because some of the estimates we need to prove here are 
essentially tri-linear in nature.
To define these norms for a given test function $u$ at frequency $\lambda$, 
we first consider all ways that one may write:
\begin{equation}
        u \ = \ \sum_\alpha \ u^\alpha \ . \label{u_sum}
\end{equation} 
For each $u^\alpha$, we consider a set of solid
angles, $\{\theta_{\alpha,d}\}$, with the property that each
$\sqrt{\frac{d}{\lambda}}\leqslant |\theta_{\alpha,d}|$. We then measure:
\begin{equation}
        \lp{u^\alpha}{Z^{ \{\theta_{\alpha,d}\} }_\lambda} \ = \
	\sum_d\ |\theta_{\alpha,d}|\ \int \left(
	\sum_{ \theta_{\alpha,d} } \ \sup_{ \omega \subseteq \theta_{\alpha,d}}
	\lp{S^\omega_{\lambda,d}u^\alpha\, (t)}{L^\infty_x}^2\right)^\frac{1}{2}
	dt \ . \label{fixed_Z_norm}
\end{equation}
In the above expression, the inclusion $\omega \subseteq \theta_{\alpha,d}$
indicates that the solid angle $\omega$, when considered as a spherical cap
of dimension $\sqrt{\frac{d}{\lambda}}\times\sqrt{\frac{d}{\lambda}}\times
\sqrt{\frac{d}{\lambda}}$ (i.e. one has $|\omega| \lesssim 
\sqrt{\frac{d}{\lambda}}$), is contained in the spherical cap defined by
$\theta_{\alpha,d}$. It is important to note that the multipliers
$S^\omega_{\lambda,d}$ and $(\lambda d)\varXi^{-1}S^\omega_{\lambda,d}$
are bounded on the norm \eqref{fixed_Z_norm}. The proof of this is 
similar to the proofs of \eqref{sqsum_L1Linf_bound} and 
\eqref{sqsum_L1Linf_bound_no_box} above. We record these facts as:\\

\begin{align}\hline \hline \notag \\
        \int \left(
	\sum_{ \theta } \ \sup_{ \omega \subseteq \theta}
	\lp{\frac{1}{\varXi}
	S^\omega_{\lambda,d}u\, (t)}{L^\infty_x}^2\right)^\frac{1}{2}
	dt \ &\lesssim \ \frac{1}{\lambda d}\int \left(
	\sum_{ \theta } \ \sup_{ \omega \subseteq \theta}
	\lp{ S^\omega_{\lambda,d}u\, (t)}{L^\infty_x}^2\right)^\frac{1}{2}
	dt  \ , \label{spec_sqsum_L1Linf_bound} \\
	\int \left(
	\sum_{ \theta } \ \sup_{ \omega \subseteq \theta}
	\lp{ S^\omega_{\lambda,d}u\, (t)}{L^\infty_x}^2\right)^\frac{1}{2}
	dt \ &\lesssim \ \int \left(
	\sum_{ \theta } \ \sup_{ \omega \subseteq \theta}
	\lp{ B^\omega_{(\frac{d}{\lambda})^\frac{1}{2}}P_\lambda
	u\, (t)}{L^\infty_x}^2\right)^\frac{1}{2} dt 
	\ . \label{spec_sqsum_L1Linf_bound_no_box}
        \\ \notag \\ \hline \hline \notag 
\end{align}\\

\noindent We now define the $Z_\lambda$ norms to be the infimum
over all possible choices of the sum \eqref{u_sum} and the angle sets
$\{\theta_{\alpha,d}\}$:
\begin{equation}
        \lp{u}{Z_\lambda} \ = \ \inf_{u = \sum_\alpha u^\alpha}
	\left(\inf_{ \{\theta_{\alpha,d}\} } \ \sum_\alpha\
	\lp{u^\alpha}{Z^{ \{\theta_{\alpha,d}\} }_\lambda}\right) 
	\ . \label{Z_norm}
\end{equation}
It is important to note here that in the vast majority of instances,
we will only need to estimate \eqref{fixed_Z_norm} for $u^\alpha=u$ and
$|\theta_{\alpha,d}| \sim \sqrt{\frac{d}{\lambda}}$. That is, for the most
part we will be dealing with the norm:
\begin{equation}
        \lp{u}{Z'_\lambda} \ = \ \sum_d\left(\frac{d}{\lambda}
	\right)^\frac{1}{2}\int \left(\sum_\omega\lp{S^\omega_{\lambda,d}u\,
	(t)}{L^\infty_x}^2\right)^\frac{1}{2} dt \ . \label{Z'_norm}
\end{equation}
Notice that one has the inclusion $Z'_\lambda \subseteq Z_\lambda$, so it is
sufficient to be able to bound the right hand side of \eqref{Z'_norm}. The
only instance where it is more convenient to work with the larger angles in
the norm \eqref{Z_norm}, is in estimate \eqref{LH_L1Linf_Omega_on_low} below. 
See the comments in the proof of that estimate for more information as to why
\eqref{Z_norm} is necessary.\\

We will also need an analog of the $Z_\lambda$ norm for situations where
$u$ cannot absorb an extra angular derivative. This is given by the following:
\begin{equation}
        \lp{u}{Z_{\Omega,\lambda}} \ = \ \sum_d\ \int
	\sup_\omega\lp{S^\omega_{\lambda,d}u\, (t)}{L^\infty_x}\, dt \ .
	\label{ang_Z_norm}
\end{equation}
The reason we have defined $Z_{\Omega,\lambda}$ separately instead of
defining it as $Z_\lambda$ is mainly for notational convenience, as it will 
give a uniform look to the statements of various estimates which follow.
We now add the norms \eqref{Z_norm} and
\eqref{ang_Z_norm} together with the norm \eqref{ang_fixedF_norm} to form
the fixed frequency version of the main function space which we will iterate in:
\begin{equation}
        G_{\Omega,\lambda} \ = \ F_{\Omega,\lambda}\ \cap\ \lambda
	|\Omega|^{-1}Z_\lambda\ \cap \
	\lambda Z_{\Omega,\lambda} \ , \label{fixed_G_norm}
\end{equation}
where $\lambda |\Omega|^{-1}Z_\lambda$ is the space with the norm
$\lambda^{-1}\lp{|\Omega| u}{Z_\lambda}$, while 
$\lambda Z_{\Omega,\lambda}$
has the norm $\lambda^{-1}\lp{u}{Z_{\Omega,\lambda}}$. The overall spaces
we will use are the analogs of the Besov space $\dot{B}^{1,1}$ and the
Sobolev spaces $\dot{H}^s$ when $1 < s$:
\begin{align}
        \lp{u}{G^1_\Omega} \ &= \ \sum_\lambda  \ \lambda\, \lp{u}
	{G_{\Omega,\lambda}} \ , \label{G_norm} \\
	\lp{u}{F_\Omega^s}^2 \ &= \ \sum_\lambda \ \lambda^{2s}\lp{u}
	{F_{\Omega,\lambda}}^2 \ . \label{F_norm} 
\end{align}
The space $G^1_\Omega$ will be the foundation of our iteration procedure, 
while the $F^s_\Omega$ space is axillary and will be used to show
that solutions to \eqref{generic_system} retain any extra regularity 
inherent in the initial data.\\

\ret

\subsection{Scattering in the $F$ spaces}
As it turns out, our scattering result, Theorem \ref{scattering_result},
is contained for
free in the structure of the $F$ spaces in the sense that any element of
those spaces can be approximated by a free wave (solution to the homogeneous
wave equation) at temporal infinity. Using a straightforward approximation
argument where one truncates the very high and very low frequencies, 
we can reduce things to proving scattering for 
functions truncated at a fixed dyadic frequency:\\

\begin{lem}[Scattering in the space $F_{\Omega,\lambda}$]
\label{fixed_freq_scattering}
For any function $u_\lambda \in F_{\Omega,\lambda}$, 
there exists a set of initial
data \ $(f_\lambda^\pm,g_\lambda^\pm)\
\in \ P_\lambda(L_\Omega^2)\times\lambda P_\lambda (L_\Omega^2)$
\ such that the following asymptotic holds:
\begin{align}
        \lim_{t\to \infty}\lp{u_\lambda(t)- W(f_\lambda^+,g_\lambda^+)(t)}
        {\dot{H}_\Omega^1\cap\partial_t(L_\Omega^2)} \ 
	&= \ 0 \ , \label{plus_scat_dyadic}\\
        \lim_{t\to -\infty}\lp{u_\lambda(t)- W(f_\lambda^-,g_\lambda^-)(t)}
        {\dot{H}_\Omega^1\cap\partial_t(L_\Omega^2)} \ 
	&= \ 0 \ . \label{minus_scat_dyadic}
\end{align}
\end{lem}\ret

\begin{proof}[proof of Lemma \ref{fixed_freq_scattering}]
The proof depends on the fact that one may write:
\begin{equation}
        u_\lambda \ = \ u_{\mathring X_{\Omega,\lambda}} 
	+ u^+_{X^{1/2}_{\Omega,\lambda,1}} +
        u^-_{X^{1/2}_{\Omega,\lambda,1}} + u_{Y_{\Omega,\lambda}} \ , \notag
\end{equation}
where $u_{\mathring X_{\Omega,\lambda}}$ is a solution to the homogeneous
wave equation with $L^2_\Omega$ data, $u^\pm_{X^{1/2}_{\Omega,\lambda,1}}$
are functions in $X^{1/2}_{\Omega,\lambda,1}$ whose Fourier transforms
are also functions and are cut off in the upper (resp. lower) half plane (in
Fourier space), and $u_{Y_{\Omega,\lambda}}$ is in the space
$Y_{\Omega,\lambda}$. For a discussion of this, see \cite{STR_generic1}.
We now define the scattering data implicitly by the relations:
\begin{align}
        W(f_\lambda^+,g_\lambda^+)(t)
        \ &= \ u_{\mathring X_{\Omega,\lambda}} + \int_{0}^\infty |D_x|^{-1}
        \sin\left(|D_x|(t-s)\right)\Box u_{Y_{\Omega,\lambda}}
	(s)\, ds \ , \notag \\
        W(f_\lambda^-,g_\lambda^-)(t)
        \ &= \ u_{\mathring X_{\Omega,\lambda}} + \int_{-\infty}^0 |D_x|^{-1}
        \sin\left(|D_x|(t-s)\right)\Box u_{Y_{\Omega,\lambda}}(s)\, ds \ . \notag
\end{align}
Using the fact that $\Box u_{Y_{\Omega,\lambda}}$ 
has finite $L^1(L_\Omega^2)$ norm, we are reduced to showing the limits:
\begin{equation}
        \lim_{t\to \pm\infty}
        \lp{u^+_{X^{1/2}_{\Omega,\lambda,1}}(t) + 
	u^-_{X^{1/2}_{\Omega,\lambda,1}}(t)}
        {\dot{H_\Omega}^1\cap\partial_t(L_\Omega^2)} \ = \ 0 \ . \notag
\end{equation}
This is a straightforward exercise in Plancherel's theorem and Dominated
Convergence for sequences of integrals. The key is to use the foliation
of the $X^{1/2}_{\Omega,\lambda,1}$ spaces alluded to above and the fact
that all weighted derivatives involve the spatial variable only. We refer the
interested reader to the work \cite{STR_generic1} for a full account.
\end{proof}

\ret\ret

\section{Some Preliminary Estimates of Sobolev Type}

In this section we provide some basic estimates of Sobolev type
which will be needed in the sequel, as well as some inclusions which 
result from these for the function spaces we introduced in the last 
section. We begin with the basic local Sobolev estimate, also
know as Bernstein's inequality:\\

\begin{lem}[Local Sobolev estimates]\label{local_sobolev_lem}
Let $u$ be a test function on $\mathbb{R}^4$, then 
one has the following frequency localized estimates:
\begin{equation}
	\lp{B^\omega_{\eta}P_1u}{L^p} \ \lesssim \
	\eta^{3(\frac{1}{r}-\frac{1}{p})}\, 
	\lp{u}{L^r} \ . \label{local_sobolev_est} 
\end{equation}
\end{lem}\ret

\noindent
For the convenience of the reader, we 
highlight here some specific instances of \eqref{local_sobolev_est},
some of which have been rescaled, that
will be used in the sequel. In all of the estimates below, the integration
is taken over $\mathbb{R}^4$:\\

\begin{align}\hline \hline \notag \\
	\lp{ B^\omega_{d^\frac{1}{2}}
	P_1 u}{L^\infty} \ &\lesssim \
	d^{\frac{1}{4}}\, \lp{u}{L^{6}} \ ,
	\label{LinftL6_loc_sob}\\
	\lp{ B^\omega_{(\frac{d}{\mu})^\frac{1}{2}}
	P_\mu u}{L^\infty} \ &\lesssim \
	\mu^\frac{5}{12} d^{\frac{1}{4}}\, \lp{u}{L^{6}} \ ,
	\label{low_LinftL6_loc_sob}\\
	\lp{ B^\omega_{d^\frac{1}{2}}
	P_1 u}{L^\infty} \ &\lesssim \
	d^{\frac{1}{2}-}\, \lp{u}{L^{3+}} \ ,
	\label{LinftL3_loc_sob}\\
	\lp{ B^\omega_{(\frac{d}{\mu})^\frac{1}{2}}
	P_\mu u}{L^\infty} \ &\lesssim \ \mu^{\frac{5}{6}-}
	d^{\frac{1}{2}-}\, \lp{u}{L^{3+}} \ ,
	\label{low_LinftL3_loc_sob}\\
	\lp{ B^\omega_{(\frac{d}{\mu})^\frac{1}{2}}
	P_\mu u}{L^\infty} \ &\lesssim \
	\mu^{\frac{5}{4}-}d^{\frac{3}{4}-}\, \lp{u}{L^{2+}} \ ,
	\label{LinftL2_loc_sob} \\
	\lp{ B^\omega_{(\frac{d}{\mu})^\frac{1}{2}}
	P_\mu u}{L^\infty} \ &\lesssim \
	\mu^{\frac{5}{4}}d^{\frac{3}{4}}\, \lp{u}{L^{2}} \ ,
	\label{strict_LinftL2_loc_sob} \\
	\lp{ B^\omega_{(\frac{d}{\mu})^\frac{1}{2}}
	P_\mu u}{L^\infty} \ &\lesssim \
	\mu^{\frac{5}{3}-}d^{1-}\, \lp{u}{L^{\frac{3}{2}+}} \ ,
	\label{LinftL3/2_loc_sob}\\
	\lp{ B^\omega_{(\frac{d}{\mu})^\frac{1}{2}}
	P_\mu u}{L^2} \ &\lesssim \ \mu^\frac{5}{12} d^\frac{1}{4}
	\, \lp{u}{L^\frac{3}{2}} \ , \label{L3/2_loc_sob} \\
	\lp{ B^\omega_{(\frac{d}{\mu})^\frac{1}{2}}
	P_\mu u}{L^2} \ &\lesssim \ \mu^{\frac{5}{6}-} d^{\frac{1}{2}-}
	\, \lp{u}{ L^{\frac{6}{5}+}} \ , \label{L6/5_loc_sob} \\
	\lp{ P_\mu u}{L^2} \ &\lesssim \ \mu^{2} 
	\, \lp{u}{ L^1 } \ . \label{L1_loc_sob}
	\\ \notag \\ \hline \hline \notag 
\end{align}\\

It will be of crucial importance to us in the resolution
of Theorem \ref{GWP_theorem} to know that our solutions to 
\eqref{generic_ym}
are not concentrating on small angular regions in Fourier space. 
This control of the solutions we construct will allow us to
use the improved Strichartz estimates \eqref{ang_str_est} and
\eqref{impr_ang_str_est}. However, because we must test our solutions
for angular regularity to gain this information, 
it will be necessary for us to deal with
situations where we need to squeeze some extra savings out of a term
that has no extra translation invariant
derivatives to give and therefore cannot
be put in any other space than $L^\infty(L^2)$. 
In other words, situations where there is extra angular regularity
present in the absence of dispersion. An example of this is a
$Low\times High$ frequency interaction where the angular derivative
falls on the low frequency term. In some sense, this
is the major technical difficulty which needs to be overcome when 
applying the estimates \eqref{ang_str_est} and \eqref{impr_ang_str_est} 
to nonlinear problems which contain derivatives in the nonlinearity.
To deal with this problem, we shall employ the
following Lemma which we state for arbitrary spatial dimension.
This is essentially an uncertainty principle for
how a function can be localized in the angular variable in Fourier
space:\\

\begin{lem}[Angular concentration estimates]\label{ang_const_lem}
Let $2\leqslant n$ be a given integer. Then for any test function 
u on $\mathbb{R}^n$, and any
$2 \leqslant p < \infty$ one has the following estimate:
\begin{equation}
	\sup_\omega\lp{B^\omega_{\eta}u}{L^p} \ \lesssim \
	\eta^s\, \lp{\langle\Omega\rangle^s u}{L^p} \ , \label{ang_const_est}
\end{equation}
where $0 \leqslant s < \frac{n-1}{p}$.
\end{lem}\ret


\begin{proof}[proof of estimate \eqref{ang_const_est}]
Upon rotation onto the positive $\xi$ axis in Fourier space, the multiplier
$b^\omega_\eta(\xi)$ satisfies the following differential bounds:
\begin{align}
        |   \partial_{\xi_{i_1}}\ldots\partial_{\xi_{i_k}} b^\omega_\eta\ (\xi)| \ \lesssim \
	2^{-jk} \ ,& &\xi_1 \in [2^j,2^{j+1}] \ . \notag
\end{align}
Furthermore, for each fixed fixed $\xi_1 \in [2^j,2^{j+1}]$, one has that the support
of $b^\omega_\eta(\xi_1,\xi')$ lies in (perhaps some thickening of)
the region $\xi' \in [-2^{j+1},2^{j+1}]\times \ldots \times 
[-2^{j+1},2^{j+1}]$. Therefore, for every permutation of the variables
$(\xi_1,\ldots,\xi_n)$, one has the following integral bound:
\begin{equation}
        \sup_{0 < k \leqslant n} \ \sup_{\xi_{k+1} , \ldots , \xi_{n}} \
	\int_\mathcal{D}\ \big| \partial_{\xi_1}\ldots\partial_{\xi_k}
	\, b^\omega_\eta(\xi_1,\ldots,\xi_{k};\xi_{k+1},\ldots,\xi_n)\big|
	\ d\xi_1\ldots d\xi_{k} \ \lesssim \ 1 \ , \notag
\end{equation}
where $\mathcal{D}$ is any dyadic rectangle\footnote{That is, one generated
by the usual dyadic partition of the coordinate axis.}
on $\mathbb{R}^k$. Therefore, by the Marcinkiewicz Multiplier
Theorem (see for instance p. 109 of \cite{Ssi}), one has that for
$2\leqslant p <\infty$:
\begin{equation}
        \lp{B^\omega_{\eta}u}{L^p} \ \lesssim \ \lp{u}{L^p} \ . \notag
\end{equation}
Therefore, by using the $n$--dimensional version of
Proposition \ref{sph_interp_prop}, it suffices to show the following estimate in $L^2$:
\begin{align}
        \lp{B^\omega_{\eta}u}{L^2} \ \lesssim \ \eta^s\ \lp{\langle \Omega
	\rangle^s\, u}{L^2} \ ,& &s \ < \ \frac{n-1}{2} \ . \notag
\end{align}
By the Plancherel Theorem and the fact that $\langle \Omega \rangle^s$
commutes with the Fourier transform, this is equivalent to showing that:
\begin{align}
        \lp{b^\omega_{\eta}\widehat u}{L^2} \ 
	\lesssim \ \eta^s\ \lp{\langle \Omega
	\rangle^s\, \widehat u}{L^2} \ ,& &s \ < \ \frac{n-1}{2} \ . \notag
\end{align}
This can now be accomplished via a simple use of
H\"olders inequality followed by the $n$ dimensional
version of the angular Sobolev embedding Proposition \eqref{sph_sob_embed}
as follows: For
$s  <  \frac{n-1}{2}$ we compute that:
\begin{align}
        \lp{b^\omega_{\eta}\widehat u}{L^2}^2 \ &= \ 
	\int_0^\infty \ \lp{b^\omega_{\eta}\widehat u\, (r)}
	{L^2(\mathbb{S}^{n-1})}^2 \ r^{n-1} \ dr \ , \notag \\
	&\lesssim \ \int_0^\infty \ \lp{b^\omega_\eta}
	{L^\frac{n-1}{s}(\mathbb{S}^{n-1})       }^2 \ 
	\lp{\widehat u\, (r)}{L^\frac{2(n-1)}{n-1-2s}(\mathbb{S}^{n-1})}^2
	\ r^{n-1} \ dr \ , \notag \\
	&\lesssim \ \eta^{2s} \ \int_0^\infty \ 
	\lp{\langle \Omega\rangle^s\, \widehat u\, (r)}
	{L^2(\mathbb{S}^{n-1})}^2 \ r^{n-1} \ dr \ , \notag \\
	&= \ \eta^{2s} \ \lp{\langle \Omega\rangle^s \widehat u}{L^2}^2
	\ . \notag
\end{align}
This ends the proof \eqref{ang_const_est}.
\end{proof}\ret

\noindent
We list here some specific incarnations of \eqref{ang_const_est} which
will appear in the sequel. These norms are taken over $\mathbb{R}^4$ for
the Lebesgue spaces, and $\mathbb{R}^{(4+1)}$ for the $F$ spaces:\\

\begin{align}\hline \hline \notag \\
	\sup_\omega \lp{ B^\omega_{(\frac{d}{\mu})^\frac{1}{2}}
	u}{L^{3+}} \ &\lesssim \ \left(\frac{d}{\mu}\right)^\frac{1}{4}
	\, \lp{\langle\Omega\rangle
	^\frac{1}{2}\, u}{L^{3+}} \ , \label{L3_ang_const} \\
	\sup_\omega\lp{ B^\omega_{(\frac{d}{\mu})^\frac{1}{2}}
	u}{L^2} \ &\lesssim \ \left(\frac{d}{\mu}\right)^\frac{1}{2}
	\, \lp{\langle\Omega\rangle\, u}{ L^2} \ , \label{L2_ang_const} \\
	\sup_\omega \lp{B^\omega_{(\frac{d}{\mu})^\frac{1}{2}}
	u}{F_\lambda} \ &\lesssim \ \left(\frac{d}{\mu}\right)^\frac{1}{2}
	\, \lp{u}{F_{\Omega,\lambda}} \ , \label{F_ang_const} 
	\\ \notag \\ \hline \hline \notag 
\end{align}\\

We conclude this section by using the Lemmas 
\ref{ang_const_lem}--\ref{local_sobolev_lem} to 
show that the $Z$ norms need only be
recovered for those pieces of $u\in F_{\Omega,\lambda}$ which are 
in the $X$ space portion of things:\\

\begin{lem}[$Z$ norm recovery for functions in the $Y$ spaces] 
On $\mathbb{R}^{4+1}$ one has the following uniform inclusions:
\begin{align}
        Y_\lambda \ &\subseteq \lambda Z'_\lambda \ , \label{YZ_incl} \\
	Y_{\Omega,\lambda} \ &\subseteq \ \lambda Z_{\Omega,\lambda} \ . \
	\label{ang_YZ_incl}
\end{align}
\end{lem}\ret

\begin{proof}[proof of the inclusions \eqref{YZ_incl} and \eqref{ang_YZ_incl}]
We'll start with the inclusion \eqref{YZ_incl}. For a fixed $d$, we can use
the special multiplier bound \eqref{sqsum_L1Linf_bound} to compute that:
\begin{align}
        &\int \left( \sum_\omega\lp{S^\omega_{\lambda,d}u\, (t)}{L^\infty_x}^2
	\right)^\frac{1}{2}dt \ , \notag \\
	\lesssim \ &\int \left( 
	\sum_\omega\lp{\frac{1}{\varXi}S^\omega_{\lambda,d}\Box u\, 
	(t)}{L^\infty_x}^2
	\right)^\frac{1}{2}dt \ , \notag \\
	\lesssim \ &\frac{1}{\lambda d}\int \left( 
	\sum_\omega\lp{S^\omega_{\lambda,d}\Box u\, 
	(t)}{L^\infty_x}^2
	\right)^\frac{1}{2}dt \ . \notag
\end{align}
We now use the local Sobolev estimate \eqref{strict_LinftL2_loc_sob} 
and the multiplier
boundedness Lemma \ref{generic_mult_lemma} to conclude that:
\begin{align}
        &\int \left( 
	\sum_\omega\lp{S^\omega_{\lambda,d}\Box u\, 
	(t)}{L^\infty_x}^2
	\right)^\frac{1}{2}dt \ , \notag \\
	\lesssim \ &\lambda^\frac{5}{4} d^\frac{3}{4}\ \int \left( 
	\sum_\omega\lp{S^\omega_{\lambda,d}\Box u\, 
	(t)}{L^2_x}^2
	\right)^\frac{1}{2}dt \ , \notag \\
	\lesssim \ &\lambda^\frac{5}{4} d^\frac{3}{4}\ \int  
	\lp{S_{\lambda,d}\Box u\, 
	(t)}{L^2_x}\, dt \ , \notag \\
	\lesssim \ &\lambda^\frac{5}{4} d^\frac{3}{4}
	\lp{\Box S_\lambda u}{L^1(L^2)} \ . \notag
\end{align}
Multiplying through by the factor $(\lambda d)^{-1}(d/\lambda)^\frac{1}{2}$
summing the last line of the above estimate over $d$ yields:
\begin{align}
        \lp{u}{Z'_\lambda} \ &\lesssim \ \sum_d 
	\left(\frac{d}{\lambda}\right)^\frac{1}{4}\ 
	\lp{\Box S_\lambda u}{L^1(L^2)} \ , \notag \\
	&\lesssim \ \lambda\, \lp{u}{Y_\lambda} \ . \notag
\end{align}\\

To prove the inclusion \eqref{ang_YZ_incl}, we first use the special
multiplier bound \eqref{sup_L1Linf_bound}, the local Sobolev
estimate \eqref{strict_LinftL2_loc_sob}, and the concentration estimate
\eqref{L2_ang_const} to compute that for fixed $d$ one has:
\begin{align}
        &\int\ \sup_\omega \lp{S^\omega_{\lambda,d}u\, (t)}{L^\infty_x}
	\, dt \ , \notag \\
	\lesssim \ &\frac{1}{\lambda d}\int 
	\sup_\omega \lp{S^\omega_{\lambda,d}\Box u\, (t)}{L^\infty_x}
	\, dt \ , \notag \\
	\lesssim \ &\left(\frac{\lambda}{d}\right)^\frac{1}{4} \ \int
	\ \sup_\omega \lp{S^\omega_{\lambda,d}\Box u\, (t)}{L^2_x}
	\, dt \ , \notag \\
	\lesssim \ &\left(\frac{d}{\lambda}\right)^\frac{1}{4} \ \int
	\ \lp{S_{\lambda,d}\Box \langle\Omega\rangle u\, (t)}{L^2_x}
	\, dt \ . \notag 
\end{align}
Using the multiplier bound \eqref{generic_mult_lemma} and summing 
over $d$ now yields:
\begin{align}
        \lp{u}{Z_{\Omega,\lambda}} \ &\lesssim \
	\sum_d \left(\frac{d}{\lambda}\right)^\frac{1}{4}
	\lp{S_{\lambda,d}\langle\Omega\rangle u}{L^1(L^2)} \ , \notag \\
	&\lesssim \ \lambda\ \lp{u}{Y_{\Omega,\lambda}} \ . \notag
\end{align}
\end{proof}

\ret\ret

\section{Bilinear Decompositions for Small Angles}\label{bilin_decomp_sect}

In this section, we will list various bilinear decompositions
for frequency localized products of the form:
\begin{equation}
        S_\mu(S_{\lambda_1}u\cdot S_{\lambda_2}v) \ . \label{basic_product}
\end{equation}
The discussion here will be in the same style as
Section 8 of \cite{STR_generic1}, and we will leave 
some of the details to that paper.
Our first task will be to decompose certain instances of the 
product \eqref{basic_product}, where $\lambda_1 \sim \lambda_2$ and
$\mu \lesssim \lambda_1$. This, of course, will be used in the
$High \times High$ type frequency interactions in the sequel. Here we will
pay special attention to the localized product:
\begin{equation}
        S_{\mu,d}( 
	S_{\lambda,\bullet \leqslant \min\{d,c\mu\}}u  \cdot 
	S_{\lambda,\bullet \leqslant \min\{d,c\mu\}}v ) \ , 
	\label{important_high_product}
\end{equation} 
where $c\ll 1$ is some small constant which we shall fixed in the proof.
The most important feature of \eqref{important_high_product} is that
the two terms in the product are restricted to be closer to the light--cone
in Fourier space than their output. By a simple computation of the 
convolution variables, this in turn implies that the \emph{angle} of
interaction between these two terms is restricted to within
$O(\frac{d}{\mu})$ of the angle of the output (in Fourier space of course).
To see this, notice that all we are talking about here is a matter of
the support of the cutoff function associated with the various 
multipliers in \eqref{important_high_product}. Therefore, it suffices to study
the convolution:
\begin{equation}
        s_{\mu,d}( 
	s_{\lambda,\bullet \leqslant \min\{d,c\mu\}} \ * \ 
	s_{\lambda,\bullet \leqslant \min\{d,c\mu\}} ) \ , 
	\label{important_high_conv}
\end{equation}
Because of the restrictions involved, we can assume without loss of generality
that the two terms in the convolution \eqref{important_high_conv} are
supported in the lower resp. upper half plane. Then for any
$(\tau,\xi)\in supp\{ s^-_{\lambda,\bullet \leqslant \min\{d,c\mu\}} \}$
and 
$(\tau',\xi')\in supp\{ s^+_{\lambda,\bullet \leqslant \min\{d,c\mu\}} \}$
we have that:
\begin{align}
        O(d) \ &= \ \Big| \ |\tau + \tau'| \ - \  
	|\xi + \xi'| \ \Big| \ , \notag \\
	&= \ \Big| \ \big| |\xi| - |\xi'| + O(d)\big| 
	\ - \ |\xi + \xi'| \ \Big| \ , \notag \\
	&= \ \Big| \ O(d) \ + \ \big| |\xi| - |\xi'| \big| 
	\ - \ |\xi + \xi'| \ \Big| \ . \notag 
\end{align}
Therefore we are able to conclude that:
\begin{equation}
        \lambda \Theta_{\xi,-\xi'}^2 \ \lesssim 
	\ \Big| \ \big| |\xi| - |\xi'| \big| 
	\ - \ |\xi + \xi'| \ \Big| \ = \ O(d) \  , \notag
\end{equation}
showing that we in fact have $|\Theta_{\xi,-\xi'}| \lesssim 
\sqrt{\frac{d}{\lambda}}$.
However, this extra precision will not concern us here, and it will
suffice to know that $|\Theta_{\xi,-\xi'}| \lesssim \sqrt{\frac{d}{\mu}}$.
It will also be important for us to know that the angle between
$\xi + \xi'$ and $\pm\xi$ does not exceed $O(\sqrt{\frac{d}{\mu}})$. 
This information will be used
to decompose products of the form \eqref{important_high_conv} by only
localizing one of the factors in Fourier space. This will be crucial
to us in the proof of estimate \eqref{HH_L1Linf_with_Omega} below, 
where the use of eccentric
multipliers in conjunction with the coarse scale decomposition of 
physical space needed for Proposition \ref{imp_ang_str_th} 
will not be possible 
due to the uncertainty principle. To obtain this type of decomposition, we 
compute:
\begin{align}
        O(d) \ &= \ \Big| \ |\tau'| \ - \  
	|\xi'| \ \Big| \ , \notag \\
	&= \ \Big| \ |(\tau + \tau') - \tau | \ - \  
	|(\xi + \xi') - \xi | \ \Big| \ , \notag \\
	&= \ \Big| \ \big| \pm |\xi + \xi'| - |\xi| + O(d)\big| 
	\ - \ |(\xi + \xi') - \xi| \ \Big| \ , \notag \\
	&= \ \Big| \ O(d) \ + \ \big| \pm |\xi + \xi'| - |\xi| \big| 
	\ - \ |(\xi + \xi') - \xi| \ \Big| \ , \notag \\
	&= \ \Big| \ O(d) \ + O(\mu \Theta_{\xi+\xi', \pm\xi}^2)
	\ \Big| \ . \notag
\end{align}
Therefore we must have that $O(\mu \Theta_{\xi+\xi', \pm\xi}^2) = O(d)$,
or $\Theta_{\xi+\xi', \pm\xi} = O(\sqrt{\frac{d}{\mu}})$.\\

The other piece of information we will need to know here is  that the
restricted convolution \eqref{important_high_conv} is supported
in the region where $|\tau + \tau'| \lesssim |\xi + \xi'|$. 
This will be used to deconstruct the $s_{\mu,d}$
multiplier that appears there into a sum over the angular
pieces $s^\omega_{\mu,d}$, each of which is supported in that
region. We calculate:
\begin{align}
        |\tau + \tau'| \ &= \ \big| |\xi| - |\xi'| + 
	O(\min\{c\mu,d\}) \big| \ , \notag \\
	&\lesssim \ \big| |\xi| - |\xi'| \big| + O(\min\{c\mu,d\})
	\ , \notag \\
	&\lesssim \ |\xi + \xi'| + O(\min\{c\mu,d\})
	\ . \notag 
\end{align}
A quick look at the support of $s_{\mu,d}$ shows that if one had
$|\xi + \xi'| \ll |\tau - \tau'|$, one must have $|\xi + \xi'| \ll \mu$
and $|\tau + \tau'|\sim \mu$. But, by the last line above this would imply 
that $|\tau + \tau'| \lesssim O(\min\{c\mu,d\})$. Thus, as long as $c$ is
chosen to be sufficiently small, we would have $|\tau + \tau'|\ll\mu$,
a contradiction. Therefore, throughout  this paper, we will assume that
$c$ has been chosen so small as to guarantee 
$|\tau + \tau'| \lesssim |\xi + \xi'|$ for convolutions of the type
\eqref{important_high_conv}.\\

What the above calculations taken together show, is
that one may replace the cutoff $s_{\mu,d}$ in the restriction
\eqref{important_high_conv} with a sum over the cutoffs
$s^\omega_{\mu,d}$, and for each term in this sum one has 
that the two factors in the convolution are supported on antipodal blocks
which differ by an angle at most $O(\sqrt{\frac{d}{\mu}})$. Furthermore,
the sum over $\omega$ is essentially diagonal with respect to these antipodal
block pairs in that there is only one (essentially) pair of antipodal
blocks for each $\omega$. The following diagram is useful for visualizing 
this:

\begin{figure}[h]
        \scalebox{.75}{\hspace{.4in}\includegraphics{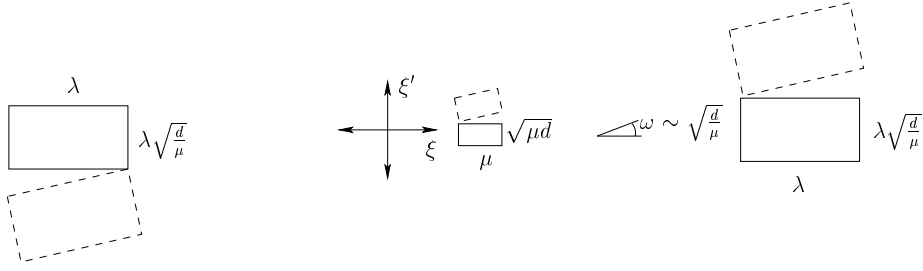}}
        \caption{Spatial supports of multipliers in an angular decomposition.}
        \label{HH_angle_decomp_fig}
\end{figure}\ret

\noindent We note here that by using the same computations as above,
one can provide similar decompositions for expressions of the form
$S_{\mu,\bullet \leqslant d} (S_{\lambda,\bullet \leqslant d}u
\cdot \nabla S_{\lambda,d}v)$ and  
$S_{\mu,\bullet \leqslant d} (S_{\lambda,d}u\cdot\nabla
S_{\lambda,\bullet < d}v)$, where $d$ is in the range $d < c\mu$ with
$c\ll 1$ the fixed small number defined above.
We recored all of these decompositions in the following:\\

\begin{lem}[$High\times High$ angular 
decomposition]\label{HH_angle_decomp_lem}
For the expression:
\begin{equation}
	S_{\mu,d}( 
	S_{\lambda,\bullet \leqslant \min\{d,c\mu\}}u  \cdot \nabla
	S_{\lambda,\bullet \leqslant \min\{d,c\mu\}}v ) \ , \notag 
\end{equation}
one has the following angular decomposition:

\begin{align}
        &s^\pm_{\mu,d}( 
	s^-_{\lambda,\bullet \leqslant \min\{d,c\mu\}}  * 
        s^+_{\lambda,\bullet \leqslant \min\{d,c\mu\}} ) \ , 
	\notag \\ \notag \\
	= \ \sum_{\substack{ \omega_1,\omega_3 \ : \\
        |\omega_1\mp\omega_3|\sim (\frac{d}{\mu})^\frac{1}{2} } }
        &{s^{\omega_1}_{\mu,d}}^\pm \left(
        s^-_{ \lambda,\bullet \leqslant \min\{d,c\mu\} }  \ * \ 
        b^{\omega_3}_{(\frac{d}{\mu})^\frac{1}{2} }
        s^+_{\lambda, \bullet \leqslant \min\{d,c\mu\} } \right) \ , 
	\notag \\ \notag \\
        = \ \sum_{\substack{ \omega_1,\omega_2,\omega_3 \ : \\
        |\omega_1\mp\omega_3|\sim (\frac{d}{\mu})^\frac{1}{2} \\
        |\omega_2 + \omega_3| \sim (\frac{d}{\mu})^\frac{1}{2} } }
        &{s^{\omega_1}_{\mu,d}}^\pm \left(
	b^{\omega_2}_{ (\frac{d}{\mu})^\frac{1}{2} }
        s^-_{ \lambda,\bullet \leqslant \min\{d,c\mu\} }  \ * \ 
        b^{\omega_3}_{ (\frac{d}{\mu})^\frac{1}{2} }
        s^+_{\lambda, \bullet \leqslant \min\{d,c\mu\} } \right) \ , \notag
\end{align}\\

\noindent
for the convolution of the associated cutoff functions in Fourier space.
There is a similar decomposition for the terms
$S_{\mu,\bullet \leqslant d} (S_{\lambda,\bullet \leqslant d}u
\cdot \nabla S_{\lambda,d}v)$ and  \\
$S_{\mu,\bullet \leqslant d} (S_{\lambda,d}u\cdot\nabla
S_{\lambda,\bullet < d}v)$, where $d$ is in the range $d < c\mu$ and
$c\ll 1$ is the small number fixed above.
\end{lem}\ret

\noindent
Because the sum on the right hand side of the expression in the above
lemma is essentially diagonal, we shall save notation in the sequel by 
abusively writing:\\

\begin{align}\hline \hline \notag \\
	&S_{\mu,d}( 
	S_{\lambda,\bullet \leqslant \min\{d,c\mu\}}u\cdot \nabla 
        S_{\lambda,\bullet \leqslant \min\{d,c\mu\}}v ) \ , \notag \\
	= \ \sum_\omega &S^{\pm\omega}_{\mu,d} \left(
        S_{ \lambda,\bullet \leqslant \min\{d,c\mu\} }u\cdot\nabla   
        B^{\omega}_{ (\frac{d}{\mu})^\frac{1}{2} }
        S_{\lambda, \bullet \leqslant \min\{d,c\mu\} }v \right) \ , 
	\label{HH_angle_partial_decomp} \\
	= \ \sum_\omega &S^{\pm\omega}_{\mu,d} \left(
	B^{-\omega}_{ (\frac{d}{\mu})^\frac{1}{2} }
        S_{ \lambda,\bullet \leqslant \min\{d,c\mu\} }u\cdot\nabla
        B^{\omega}_{ (\frac{d}{\mu})^\frac{1}{2} }
        S_{\lambda, \bullet \leqslant \min\{d,c\mu\} }v \right) \ . 
	\label{HH_angle_decomp}
	\\ \notag \\ \hline \hline \notag 
\end{align}\\

\noindent 
We shall also use this shorthand for the same decomposition applied to the 
terms $S_{\mu,\bullet \leqslant d} (S_{\lambda,\bullet \leqslant d}u\cdot\nabla
S_{\lambda,d}v)$
and $S_{\mu,\bullet \leqslant d} (S_{\lambda,d}\cdot\nabla
S_{\lambda,\bullet < d}v)$.\\

Next, we move on to several decompositions which are
dual to Lemma \eqref{HH_angle_decomp_lem}. These take place in the 
presence of a $Low\times High$ frequency interaction. The first such
decomposition will be used when the low frequency term controls the
angles. The validity of this decomposition follows from essentially
the same calculations as used for \eqref{HH_angle_partial_decomp} and
\eqref{HH_angle_decomp} above. For a proof, see \cite{STR_generic1}.
We record it here as:\\

\begin{lem}[$Low\times High$ wide angle 
decomposition]\label{LH_wide_angle_decomp_lem}
For the expression:
\begin{equation}
	S_{\lambda , \bullet < \min\{d,c\mu\}} (S_{\mu,d}u\, 
        \nabla S_{\lambda,\bullet < \min\{d,c\mu\}} v) \ , \notag
\end{equation}
One has the following angular decomposition:

\begin{multline}
	s^+_{\lambda,\bullet < \min\{c\mu,d\}}
	( s^\pm_{\mu,d} * 
        s^+_{\lambda,\bullet < \min\{c\mu,d\}}) \\ 
	= \ \sum_{\substack{ \omega_1,\omega_2,\omega_3 \ : \\
        |\omega_1\mp\omega_2|\sim (\frac{d}{\mu})^\frac{1}{2} \\
        |\omega_1 - \omega_3| \sim (\frac{d}{\mu})^\frac{1}{2} } }
        b^{\omega_1}_{(\frac{d}{\mu})^\frac{1}{2}}
        s^+_{\lambda,\bullet < \min\{c\mu,d\}}\left(
        {s^{\omega_2}_{\mu,d}}^\pm \ * \ 
        b^{\omega_3}_{(\frac{d}{\mu})^\frac{1}{2}}
        s^+_{\lambda,
        \bullet < \min\{c\mu,d\}}\right) \ . \notag
\end{multline}\\

\noindent
for the convolution of the associated cutoff functions in Fourier
space. There is a similar decomposition for the terms
$S_{\lambda,\bullet \leqslant d} (S_{\mu,\bullet \leqslant d}u\, 
\nabla S_{\lambda,d} v)$ and \\
$S_{\lambda , d} (S_{\mu,\bullet \leqslant d}u\, 
\nabla S_{\lambda,\bullet < d} v)$ in the range $d < c\mu$, 
where $c\ll 1$ is a fixed small number.
\end{lem}\ret 

\noindent
We will also need a decomposition similar to that of Lemma 
\ref{LH_wide_angle_decomp_lem} for the case where the high frequency
term controls the angle. Again, for a proof see 
\cite{STR_generic1}. This is:\\

\begin{lem}[$Low\times High$ small angle 
decomposition]\label{LH_small_angle_decomp_lem}
For the expression:
\begin{equation}
	S^\omega_{\lambda , d} (S_{\mu,\bullet \leqslant d}u\, 
	\nabla S_{\lambda,\bullet < d} v) \ , \notag
\end{equation}
one has the following angular restriction:\\

\begin{equation}
	{s^{\omega_1}_{\lambda,d}}^+(s^\pm_{\mu,\bullet \leqslant d} *
        s^+_{\lambda, \bullet  <  d} ) \ = \
        {s^{\omega_1}_{\lambda,d}}^+(
        {s^{\omega_2}_{\mu,\bullet \leqslant d}}^\pm \ * \ 
       { s^{\omega_3}_{\lambda,\bullet  < d}}^+) \ , \notag
\end{equation}\\

\noindent
for the convolution of the associated cutoff functions in Fourier space.
Here the angles are restricted to the range
$\ \ |\omega_1 - \omega_3|\sim \sqrt{\frac{d}{\lambda}}\ \ $, and
$\ \ |\omega_1 - \pm\omega_2|\sim \sqrt{\frac{d}{\mu}}\ \ $.
\end{lem}\ret

\noindent
Finally we record here the thickened version of Lemma 
\ref{LH_small_angle_decomp_lem}:\\

\begin{lem}[Thickened $Low\times High$ small angle 
decomposition]\label{th_LH_small_angle_decomp_lem}
For the expression:
\begin{equation}
	S^\omega_{\lambda , d} (S_{\mu}u\, 
	\nabla S_{\lambda,\bullet < c\mu} v) \ , \notag
\end{equation}
one has the following angular restriction:\\

\begin{equation}
	{s^{\omega_1}_{\lambda,d}}^+(s^\pm_\mu *
        s^+_{\lambda, \bullet  <  c\mu} ) \ = \
        {s^{\omega_1}_{\lambda,d}}^+(
        s^\pm_\mu \ * \ 
        b^{\omega_3}_{(\frac{d}{c\lambda})^\frac{1}{2}}
	s^+_{\lambda,\bullet  < c\mu}) \ , \notag
\end{equation}\\

\noindent
for the convolution of the associated cutoff functions in Fourier space.
Here the angles are restricted to the range
$\ \ |\omega_1 - \omega_3|\sim \sqrt{\frac{d}{c\lambda}}\ \ $.
\end{lem}\ret

\noindent
As we have done above, in the sequel
we shall write the decompositions in Lemmas
\ref{LH_wide_angle_decomp_lem}--\ref{th_LH_small_angle_decomp_lem} using the
following shorthand:\\

\begin{align}\hline \hline \notag \\
   \begin{split}
	S_{\lambda , \bullet < \min\{d,c\mu\}} (S_{\mu,d}u\cdot 
        \nabla S_{\lambda,\bullet < \min\{d,c\mu\}} v) \ &= \\
	\sum_\omega B^\omega_{(\frac{d}{\mu})^\frac{1}{2}}
	S_{\lambda , \bullet < \min\{d,c\mu\}}&\bigg(
	S^{\pm \omega}_{\mu,d}u\ \cdot \ \nabla 
	B^\omega_{(\frac{d}{\mu})^\frac{1}{2}} 
	S_{\lambda,\bullet < \min\{d,c\mu\}} v\bigg) \ , 
   \end{split}\label{LH_wide_angle_decomp} \\
	\sum_\omega S^\omega_{\lambda , d} (S_{\mu,\bullet\leqslant d}u\cdot 
        \nabla S_{\lambda,d} v) \ &=  
	\sum_{\substack{\omega_1,\omega_2 \ :\\
	|\omega_1 \mp \omega_2|\sim(\frac{d}{\mu})^\frac{1}{2}}} 
	S^{\omega_1}_{\lambda , d}\bigg(
	S^{\pm\omega_2}_{\mu,\bullet\leqslant d}u\cdot
	\nabla S^{\omega_1}_{\lambda,d} v \bigg) \ , 
	\label{LH_small_angle_decomp} \\
	\sum_\omega S^\omega_{\lambda , d} (S_{\mu}u\, 
	\nabla S_{\lambda,\bullet < c\mu} v) \ &= \
	\sum_{\substack{\omega_1,\omega_2 \ :\\
	|\omega_1 - \omega_2|\sim(\frac{d}{c\lambda})^\frac{1}{2}}}
	S^{\omega_1}_{\lambda , d} \bigg(S_{\mu}u\, 
	\nabla S^{\omega_2}_{\lambda,\bullet < c\mu} v\bigg) \ . 
	\label{th_LH_small_angle_decomp}
	\\ \notag \\ \hline \hline \notag 
\end{align}\\

\noindent
It is important to note that the sum on the right hand side of
\eqref{LH_small_angle_decomp} above is \emph{not} diagonal in the two
angles $\omega_1$ and $\omega_2$. This is one of the main reasons
why we need to employ the extra flexibility in the norms \eqref{Z_norm}. 
Also, while the sum on the right hand side of
\eqref{th_LH_small_angle_decomp} is essentially diagonal for a fixed
small $c\ll 1$, we have elected to keep the more precise form because
we will need to pick a $c$ based on the other implicit constants which
appear in various decompositions in the proof.

\ret\ret

\section{Frequency Decomposition of the Nonlinearity}

The remainder of the paper is devoted to the proof of Theorem
\ref{GWP_theorem}. This will be done through Picard iterating
the integral equation:
\begin{equation}\label{generic_NLW_inegral}
	\phi \ = \ W(f,g) + \Box^{-1} (\phi\, \nabla \phi) \ , 
\end{equation}
in the spaces $G_\Omega^1$ and $F_\Omega^{s,2}\cap G_\Omega^1$ 
for $1 < s$. Due to the
quadratic nature of the nonlinearity, it suffices to prove the following:\\

\begin{thm}[Solution of the division problem]\label{division_thm}
For $n=4$, the $F_\Omega$ and $G_\Omega$ 
spaces solve the division problem for the system \eqref{generic_system}
in the sense that one has the following bilinear estimates for functions
$u$ and $v$:
\begin{align}
	\lp{\Box^{-1}(u\, \nabla v)}{G_\Omega^1} \ &\lesssim \
	\lp{u}{G_\Omega^1}
	\lp{v}{G_\Omega^1}\ , \label{division_estimate1} \\
	\lp{\Box^{-1}(u\, \nabla v)}{F_\Omega^{s,2}} \ &\lesssim \
	\lp{u}{G_\Omega^1}\lp{v}{F_\Omega^{s,2}} + 
	\lp{u}{F_\Omega^{s,2}}\lp{v}{G_\Omega^1}\ . 
	\label{division_estimate2}
\end{align}
\end{thm}\ret

In what follows, we shall concentrate solely on the estimate 
\eqref{division_estimate1}. The other, estimate \eqref{division_estimate2},
will follow directly from the dyadic estimates employed in the proof
of \eqref{division_estimate1}.\\

Our first step is the usual Littlewood-Paley
decomposition of the nonlinear term $\Box^{-1} (u \, \nabla v)$ 
with respect to space--time frequencies:
\begin{equation}\label{LP_decomp1}
	\Box^{-1} (u \, \nabla v) \ = \ \sum_{\lambda_i}
	\Box^{-1}(S_{\lambda_1} u \, \nabla S_{\lambda_2} v) \ .
\end{equation}
We now split the sum on the right hand side of \eqref{LP_decomp1} into the
three cases: $\lambda_1\sim\lambda_2$, $\lambda_1 \ll \lambda_2$, and
$\lambda_2\ll\lambda_1$. In the sequel, we only concentrate on the first two
interactions, the last case being similar to the second through some
standard ``weight trading''. Taking into account the $\ell^1$ Besov 
structure of the $G^1_\Omega$ space, to prove \eqref{division_estimate1},
it suffices to show the following two bilinear estimates:\\

\begin{align}
	\sum_{\mu \ : \ \mu\lesssim max\{\lambda_1,\lambda_2\}} \mu\, 
	\lp{\Box^{-1}(S_{\lambda_1} u \, \nabla S_{\lambda_2} v)}
	{G_{\Omega,\mu}}
	\ &\lesssim \ 
	\lambda_1\lambda_2 \,
	\lp{u}{F_{\Omega,\lambda_1}}\lp{v}{F_{\Omega,\lambda_2}} \ \ , \
	\lambda_1 \sim \lambda_2  \ ,
	\label{HH_est} \\ \notag \\ 
	\lp{\Box^{-1}(S_{\mu} u \, \nabla S_{\lambda} v)}{G_{\Omega,\lambda}}
	\ &\lesssim \ \mu\,
	\lp{u}{G_{\Omega,\mu}}\lp{v}{F_{\Omega,\lambda}} \ \ , \
	\mu \ll \lambda \ . \label{LH_est}
\end{align}\\

Due to the fact that both of the above estimates are scale invariant, it 
suffices to prove them for $\lambda_1 = \lambda_2 = 1$, and $\mu=1$ 
respectively. In the next two subsections, we 
shall break these estimates down further into a series
of cases involving the various function spaces that make up the 
$G_{\Omega,\lambda}$ and $F_{\Omega,\lambda}$ spaces. These estimates will be 
placed in highlighted format for the convenience of the reader. Each 
individual estimate will then be proved separately in the remaining two
sections of the paper.\\

\ret

\subsection{$High\times High$ regime: 
List of dyadic estimates corresponding to \eqref{HH_est}}\label{HH_section}

In what follows, we will not explicitly take angular derivatives
of any of the expressions we are to estimate. Instead, we will prove
bilinear estimates where at most one term in the product\footnote{With the
exception of estimate \eqref{HH_L1Linf_no_Omega} 
below which involves angular derivatives of both
terms in the product. Notice that this is acceptable because the norm
\eqref{ang_Z_norm} does not involve any angular derivatives.}
on the
right hand side contains a norm involving angular
derivatives. Thus, using the Leibniz rule, one can safely add an angular
derivative to all estimates that follow below to get estimate 
\eqref{HH_est} above.\\

We begin by further decomposing the expression $S_\mu\Box^{-1}(S_1 u\, \nabla
S_1 v)$, $\mu \lesssim 1$. Using the (approximate) idempotence
of $S_\mu$, we first compute that:
\begin{equation}
	S_\mu\Box^{-1}(S_1 u\, \nabla S_1 v)  \ = \
	\Box^{-1}S_\mu (S_1 u\, \nabla S_1 v) +
	S_\mu [S_\mu,\Box^{-1}](S_1 u\, \nabla S_1 v) \ . \notag
\end{equation}
To compute the commutator term, let $H_{\bullet\lesssim 1}$
be a function with space-time frequency $\lesssim 1$. Then we
have that:
\begin{align}
  \begin{split}
        S_\mu [S_\mu,\Box^{-1}] H_{\bullet\lesssim 1} \ &= \
	S_\mu E*H_{\bullet\lesssim 1} - S_\mu W
	\left(E*H_{\bullet\lesssim 1}\right) \notag\\
	& \ \ \ \ \ \ \ - S_\mu E*S_\mu H_{\bullet\lesssim 1} + S_\mu W
	\left(E*S_\mu  H_{\bullet\lesssim 1} 
	\right) \ , 
   \end{split} \notag \\
	&= \ -\sum_{\substack{\sigma \ : \\ 
        \mu\leqslant \sigma \lesssim 1}}S_\mu W
	\left(E*S_\sigma H_{\bullet\lesssim 1}\right)
	+ S_\mu W \left(E*S_\mu H_{\bullet\lesssim 1}
	\right) \ , \notag\\
	&= \ -\sum_{\substack{\sigma \ : \\ \mu < \sigma \lesssim 1}}S_\mu W
	\left(E*S_\sigma H_{\bullet\lesssim 1}\right)\ , \notag\\
	&= \ -\sum_{\substack{\sigma \ : \\ \mu < \sigma \lesssim 1}} W
	\left(P_\mu E*S_{\sigma,\sigma} H_{\bullet\lesssim 1} 
	\right) \ . \notag
\end{align}
Therefore, we have that:
\begin{equation}
        S_\mu\Box^{-1}(S_1 u\, \nabla S_1 v)
	\ = \ \Box^{-1}S_\mu (S_1 u\, \nabla S_1 v) \ - 
	\sum_{\sigma \ : \
	\mu < \sigma \lesssim 1}
	W\left(P_\mu\, S_{\sigma,\sigma}\, 
	\frac{1}{\varXi}(S_1 u\, \nabla S_1 v)\right) 
	\ . \label{with_out_comm}
\end{equation} 
To handle the second term on line \eqref{with_out_comm}
it is enough to prove the estimate:\\

\begin{align}\hline \hline \notag \\
	\sum_{\mu \lesssim 1} \ \mu\, \lp{\sum_{\sigma \ : \
	\mu < \sigma \lesssim 1}
	\left(P_\mu\, S_{\sigma,\sigma}\, 
	\frac{1}{\varXi}(S_1 u\, \nabla S_1 v)\right)}{L^\infty(L^2)}
	\ \lesssim \ \lp{u}{F_1}\lp{v}{F_1} \ . \label{extra_LinfL2_est}
	\\ \notag \\
	\hline \hline \notag
\end{align}\\

Next, we move on to the first term on the right hand side of 
\eqref{with_out_comm} above. For this we fix a small $c\ll 1$
as explained in Section \ref{bilin_decomp_sect} and write:
\begin{align}
	S_\mu (S_1 u\, \nabla S_1 v) \ &= \
	S_\mu (S_1 u\, \nabla S_{1,c\mu \leqslant \bullet} v)
	\ + \
	S_\mu (S_{1,c\mu \leqslant \bullet}u\, 
	\nabla S_{1,\bullet < c\mu} v)  \label{rough_high_cdist} \\
	&\ \ \ \ \ \ \ \ \ + \ S_\mu (S_{1,\bullet < c\mu  }u\, 
	\nabla S_{1,\bullet < c\mu} v)
	\ .  \notag
\end{align}
In all that follows here, we will only estimate the terms of 
\eqref{rough_high_cdist} in the $X$ and $Y$ spaces defined on lines
\eqref{Xsd_norm} and \eqref{duhamel_norm} respectively. The addition
of the $L^\infty(L^2)$ estimate, included on line \eqref{fixedF_norm}, 
follows from estimates on these first two norms in a standard way.
See \cite{STR_generic1} for details.
To control the fist two terms on the right hand side of 
\eqref{rough_high_cdist}, we will prove that:\\ 

\begin{align}\hline \hline \notag \\ 
	 \sum_\mu \ \lp{S_\mu (S_1 u\, \nabla 
	S_{1,c\mu \leqslant \bullet} v)}{L^1(L^2)} \ \lesssim \
	\lp{u}{F_1}\lp{v}{F_1} \ , \label{HH_L1L2_no_omega1} \\ 
	\sum_\mu \ \lp{ S_\mu (S_{1,c\mu \leqslant \bullet}u\, 
	\nabla S_{1,\bullet < c\mu} v)}{L^1(L^2)} \ \lesssim \
	\lp{u}{F_1}\lp{v}{F_1} \ . \label{HH_L1L2_no_omega2}
	\\ \notag \\ \hline \hline \notag 
\end{align}\\

It remains to estimate the terms $\Box^{-1} S_\mu 
(S_{1,\bullet < c\mu  }u\, \nabla S_{1,\bullet < c\mu} v)$. To do this
we decompose these expressions with respect to all possible dyadic
distances from the light cone in Fourier space. We group this sum
together as follows:
\begin{align}
	S_\mu (S_{1,\bullet < c\mu  }u\, \nabla 
	S_{1,\bullet < c\mu} v) \ &= \ \sum_{\substack{d,\delta_1,
	\delta_2 \ : \\ 
	\delta_1 < c\mu \\
	\delta_2 < c\mu }} S_{\mu,d} (S_{1,\delta_1}u\, \nabla 
	S_{1,\delta_2} v) \ , \notag \\
   \begin{split}
	&= \sum_{d\ : \ d < c\mu} 
	S_{\mu,\bullet \leqslant d} (S_{1,\bullet \leqslant d}u\, 
	\nabla S_{1,d} v)  \\
	&\ \ \ \ \ \ \ \ + \ \sum_{d\ : \ d < c\mu} 
	S_{\mu,\bullet \leqslant d} (S_{1,d}u\, 
	\nabla S_{1,\bullet < d} v) \\ 
	&\ \ \ \ \ \ \ \ \ \ \ \ \ \ \ \ \ + \
	\sum_d S_{\mu,d} (S_{1,\bullet < \min\{d,c\mu\}}u\, 
	\nabla S_{1,\bullet < \min\{d,c\mu\}} v) \ .
   \end{split}\label{fine_high_cdist}
\end{align}
For the first term on the right hand side of \eqref{fine_high_cdist},
we will prove the two estimates:\\ 

\begin{align}\hline \hline \notag \\
	\sum_\mu \ \lp{\sum_{d\ : \ d < c\mu} 
	S_{\mu,\bullet \leqslant d} (S_{1,\bullet \leqslant d}u\, 
	\nabla S_{1,d} v)}{L^1(L^2)} \ \lesssim \
	\lp{u}{F_{\Omega,1}}\lp{v}{F_1} \ , 
	\label{HH_L1L2_Omega_not_on_d} \\
	\sum_\mu \ \lp{\sum_{d\ : \ d < c\mu} 
	S_{\mu,\bullet \leqslant d} (S_{1,\bullet \leqslant d}u\, 
	\nabla S_{1,d} v)}{L^1(L^2)} \ \lesssim \
	\lp{u}{F_1}\lp{v}{F_{\Omega,1}} \ . \label{HH_L1L2_Omega_on_d}
	\\ \notag \\ \hline \hline \notag 
\end{align}\\

Notice that 
the proof of estimate \eqref{HH_est} for the second 
term in line \eqref{fine_high_cdist}
follows from the proof of 
\eqref{HH_L1L2_Omega_not_on_d}--\eqref{HH_L1L2_Omega_on_d} 
and some weight trading. It remains to
deal with the last term of line \eqref{fine_high_cdist}. 
In this case we'll need 
to rely on the $L^2$ based norms. Using some weight trading, 
this can be reduced to the single estimate:\\

\begin{align}\hline \hline \notag \\
	\sum_{\mu,d}\ d^{-\frac{1}{2}}\, \lp{
	S_{\mu,d} (S_{1,\bullet < \min\{d,c\mu\}}u\, 
	\nabla S_{1,\bullet < \min\{d,c\mu\}} v)}{L^2(L^2)} \ \lesssim \
	\lp{u}{F_{\Omega,1}}\lp{v}{F_1} \ . \label{HH_L2_est}
	\\ \notag \\
	\hline \hline \notag 
\end{align}\\

Because of the limitations of the $L^2$ based spaces (i.e. that they do not
have an analog of \eqref{YZ_incl}--\eqref{ang_YZ_incl}), we also need to
recover the $Z$ norms by hand for this frequency interaction. It suffices to
be able to bound the norms \eqref{Z'_norm} and \eqref{ang_Z_norm}. 
Using the multiplier boundedness
lemma \eqref{spec_mult_lemma} and some weight trading where necessary, 
this reduces to showing the two estimates:\\

\begin{align}\hline \hline \notag \\
   \begin{split}
	\sum_{\mu,d}\ \frac{1}{\mu^\frac{3}{2} 
	d^\frac{1}{2}}\, &\int
	\left( \sum_\omega \lp{S^\omega_{\mu,d} 
	(S_{1,\bullet < \min\{d,c\mu\}}u\, 
	\nabla S_{1,\bullet < \min\{d,c\mu\}} v)(t)}
	{L^\infty_x}^2\right)^\frac{1}{2} \, dt \\
	&\hspace{3in}
	\lesssim \  \lp{u}{F_{\Omega,1}}\lp{v}{F_{1}} \ , 
   \end{split}\label{HH_L1Linf_with_Omega} \\
   \begin{split}
	\sum_{\mu,d}\ \frac{1}{\mu d}\, &\int
	\sup_\omega \lp{S^\omega_{\mu,d} (S_{1,\bullet < \min\{d,c\mu\}}u\, 
	\nabla S_{1,\bullet < \min\{d,c\mu\}} v)(t)}{L^\infty_x}
	\, dt \\
	&\hspace{3in} \lesssim \  \lp{u}{F_{\Omega,1}}
	\lp{v}{F_{\Omega,1}} \ .
   \end{split}\label{HH_L1Linf_no_Omega}
	\\ \notag \\ \hline \hline \notag 
\end{align}\\

\ret

\subsection{$Low\times High$ regime: 
List of dyadic estimates corresponding to \eqref{LH_est}}

Here we follow the same procedure as in the previous section, proving
bilinear estimates where one term in the product can safely absorb
an angular derivative. Of course there is an exception for
estimates \eqref{spec_LH_L1Linf_no_Omega} and  \eqref{LH_L1Linf_no_Omega} 
below, which do not need to absorb any angular derivatives.\\

Our first order of business here is to note that, due to the fact
$1 \ll \lambda$, one has the identity:
\begin{equation}
	S_\lambda \Box^{-1}(S_{1} u \, \nabla S_{\lambda} v) \ = \
	S_\lambda\Box^{-1}S_\lambda(
	S_{1} u \, \nabla S_{\lambda} v) \ . \notag
\end{equation} 
Because of this, we can avoid proving an extra estimate of the form
\eqref{extra_LinfL2_est} here. We now use 
a rough decomposition to isolate things so they
are sufficiently close to the light cone in Fourier space. For this 
purpose we fix a small $c\ll 1$ as described in Section 
\ref{bilin_decomp_sect} and write:
\begin{align}
	S_\lambda(S_{1} u \, \nabla S_{\lambda} v)\ 
	&= \ S_\lambda(S_{1} u \, \nabla S_{\lambda,
	c \leqslant \bullet } v) \ + \ S_{\lambda,c \leqslant \bullet}\,
	(S_{1} u \, \nabla S_{\lambda,\bullet < c } v) 
	\label{rough_low_high_cdist}\\
	&\ \ \ \ \ \ \ \ \ + \ S_{\lambda,\bullet < c }\,
	(S_{1} u \, \nabla S_{\lambda,\bullet < c } v) \ . \notag
\end{align}
Notice that
this decomposition is dual to that of \eqref{rough_high_cdist}.
Again, as we mentioned in the previous subsection, we will
only concern ourselves with estimating the $X$ and $Y$ space
portion of the norm \eqref{fixedF_norm}, as the $L^\infty(L^2)$ type estimate
follows once we have done this.
In order to estimate $\Box^{-1}$ of
the first two terms on the right hand side of \eqref{rough_low_high_cdist}
in the $F_{\Omega,\lambda}$ spaces, it is enough to prove the following two 
estimates:\\

\begin{align}\hline \hline \notag \\
	\lp{S_\lambda(S_{1} u \, \nabla S_{\lambda,
	c \leqslant \bullet } v)}{L^1(L^2)} \ &\lesssim \
	\lp{u}{F_1}\cdot  \ \lambda \ 
	\lp{v}{F_\lambda} \ , \label{LH_L1L2_no_omega1} \\
	\sum_{c \leqslant d}\ d^{-\frac{1}{2}}\ 
	\lp{S_{\lambda,d}\,
	(S_{1} u \, \nabla S_{\lambda,\bullet < c } v)}{L^2(L^2)}
	\ &\lesssim \ \lp{u}{F_1}\cdot  \ \lambda \ \lp{v}{F_\lambda} \ . 
	\label{LH_L1L2_no_omega2}
	\\ \notag \\ \hline \hline \notag 
\end{align}\\

We also need to recover the norms \eqref{Z'_norm} and \eqref{ang_Z_norm}
by hand for the second term on
the right hand side of \eqref{rough_low_high_cdist}. To do this we will
show that:\\

\begin{align}\hline \hline \notag \\
	\sum_{c \leqslant d}\ \frac{1}{\lambda^\frac{3}{2}d^\frac{1}{2}}
	\, \int\left( \sum_\omega \lp{S^\omega_{\lambda,d} (S_1 u\,
	\nabla S_{\lambda,\bullet < c } v)(t)}{L^\infty_x}^2
	\right)^\frac{1}{2}
	 dt \ &\lesssim \   \lp{u}{F_1}\cdot  \ \lambda \ \lp{v}{F_\lambda} \ ,
	\label{spec_LH_L1Linf_with_Omega} \\
	\sum_{c \leqslant d}\ \frac{1}{\lambda d}\, \int  \sup_\omega
	\lp{S^\omega_{\lambda,d} (S_1 u\,
	\nabla S_{\lambda,\bullet < c } v)(t)}{L^\infty_x}\, dt
	 \ &\lesssim \   \lp{u}{F_{1}}\cdot  \ \lambda \ \lp{v}
	{F_{\Omega,\lambda}} \ . \label{spec_LH_L1Linf_no_Omega}
	\\ \notag \\ \hline \hline \notag 
\end{align}\\

At this point, we are in a position to deal with the last term on the 
right hand side of \eqref{rough_low_high_cdist}. We further decompose
this in a manner dual to \eqref{fine_high_cdist}:
\begin{align}
	S_{\lambda,\bullet < c }\,
	(S_{1} u \, \nabla S_{\lambda,\bullet < c } v)
	 \ &= \ \sum_{\substack{d,\delta_1,
	\delta_2 \ : \\ 
	d < c \\
	\delta_2 < c }} S_{\lambda,d} (S_{1,\delta_1}u\, \nabla 
	S_{\lambda,\delta_2} v) \ , \notag \\
   \begin{split}
	&= \sum_{d < c} 
	S_{\lambda,\bullet \leqslant d} (S_{1,\bullet \leqslant d}u\, 
	\nabla S_{\lambda,d} v)  \\
	&\ \ \ \ \ \ \ \ + \ \sum_{d < c} 
	S_{\lambda , d} (S_{1,\bullet \leqslant d}u\, 
	\nabla S_{\lambda,\bullet < d} v) \\ 
	&\ \ \ \ \ \ \ \ \ \ \ \ \ \ \ \ \ + \
	\sum_d S_{\lambda , \bullet < \min\{d,c\}} 
	(S_{1,d}u\, 
	\nabla S_{\lambda,\bullet < \min\{d,c\}} v) \ .
   \end{split}\label{fine_low_high_cdist}
\end{align}
To deal with the first term on the right hand side of 
\eqref{fine_low_high_cdist}, we'll prove the two separate estimates:\\

\begin{align}\hline \hline \notag \\
	\lp{\sum_{d < c} \ S_{\lambda,\bullet \leqslant d} 
	(S_{1,\bullet \leqslant d}u\, \nabla S_{\lambda,d} v)}{L^1(L^2)}
	\ &\lesssim \ \lp{u}{F_1}\cdot \ 
	\lambda\ \lp{v}{F_{\Omega,\lambda}} \ ,
	\label{LH_L1L2_Omega_on_low} \\
	\lp{\sum_{d < c}\ S_{\lambda,\bullet \leqslant d} 
	(S_{1,\bullet \leqslant d}u\, \nabla S_{\lambda,d} v)}{L^1(L^2)}
	\ &\lesssim \ \lp{u}{F_{\Omega,1}}\cdot \ \lambda \
	\lp{v}{F_{\lambda}} \ ,
	\label{LH_L1L2_Omega_on_high}
	\\ \notag \\ \hline \hline \notag 
\end{align}\\

To show the second term on the right hand side of \eqref{fine_low_high_cdist}
is in the $F_{\Omega,\lambda}$ spaces, we prove the pair of estimates:\\

\begin{align}\hline \hline \notag \\
	\sum_{d < c}\ d^{-\frac{1}{2}}\
	\lp{S_{\lambda,d} (S_{1,\bullet \leqslant d}u\, 
	\nabla S_{\lambda,\bullet < d} v)}{L^2(L^2)} \ &\lesssim \
	\lp{u}{F_1}\cdot \ \lambda \ \lp{v}{F_{\Omega,\lambda}} \ ,
	\label{LH_L2L2_Omega_on_low} \\
	\sum_{d < c}\ d^{-\frac{1}{2}}\
	\lp{S_{\lambda,d} (S_{1,\bullet \leqslant d}u\, 
	\nabla S_{\lambda,\bullet < d} v)}{L^2(L^2)} \ &\lesssim \
	\lp{u}{F_{\Omega,1}}\cdot  \ \lambda \ \lp{v}{F_{\lambda}} \ ,
	\label{LH_L2L2_Omega_on_high} 
	\\ \notag \\ \hline \hline \notag 
\end{align}\\

Next, we will need to recover the special $L^1(L^\infty)$ norms for the
second term on the right hand side of \eqref{fine_low_high_cdist}.
Our first estimate will used for controlling the expression
\eqref{fine_low_high_cdist} when the low frequency term in the product
need to absorb an angular derivative.
This is where we will need to use the more complicated version, \eqref{Z_norm},
of the $Z$ norms. Here, the $u^\alpha$ in the sum
\eqref{u_sum} will come from the terms in the following sum:
\begin{equation}
        S_\lambda(S_{\bullet \ll \lambda} u\, \nabla S_\lambda v)
	\ = \ \sum_{\mu \ :  \ \mu \ll \lambda} 
	\ S_{\mu} u\, \nabla S_\lambda v \ . \notag
\end{equation}
Due to the Besov structure in the $G$ space, we only need to deal with 
a single term of this. By rescaling, as we have done throughout
this section, we may assume that $\mu = 1$. For this fixed piece, we shall
use the angle set where $|\theta_d| = d^\frac{1}{2}$. Therefore, it will
suffice to prove:\\

\begin{align}\hline \hline \notag \\ 
   \begin{split}
	\sum_{d}\ \frac{1}{\lambda 
	d^\frac{1}{2}}\, &\int
	\left( \sum_{\substack{\theta \ : \\
	|\theta|\sim \, d^\frac{1}{2}}} \ \sup_{\omega \subseteq \theta}
	\lp{ S^\omega_{\lambda , \bullet < \min\{d,c\}} 
	(S_{1,d}u\, \nabla S_{\lambda,
	\bullet < \min\{d,c\}} v)(t)}{L^\infty_x}^2
	\right)^\frac{1}{2} dt \\ 
	&\hspace{2.5in} \lesssim \
	\lp{u}{F_{1}}\cdot \ \lambda \ \lp{v}{F_{\Omega,\lambda}} \ ,
   \end{split}\label{LH_L1Linf_Omega_on_low}
   \\ \notag \\ \hline \hline \notag 
\end{align}\\

When estimating the second term of 
the expression \eqref{fine_low_high_cdist} where
an angular derivative falls on the high frequency term in the product,
we shall use the simpler norm \eqref{Z'_norm}. Finally, we need to recover
the norm \eqref{ang_Z_norm}. Therefore, we shall
prove the following two estimates:\\

\begin{align}\hline \hline \notag \\
   \begin{split}
	\sum_{d}\ \frac{1}{\lambda^\frac{3}{2} 
	d^\frac{1}{2}}\, &\int
	\left( \sum_\omega \ \lp{ S^\omega_{\lambda , \bullet < \min\{d,c\}} 
	(S_{1,d}u\, \nabla S_{\lambda,
	\bullet < \min\{d,c\}} v)(t)}{L^\infty_x}^2
	\right)^\frac{1}{2} dt \\ 
	&\hspace{2.5in} \lesssim \
	\lp{u}{F_{\Omega,1}}\cdot \ \lambda \ \lp{v}{F_{\lambda}} \ ,
   \end{split}\label{LH_L1Linf_Omega_on_high}\\
  \begin{split}
	\sum_{d}\ \frac{1}{\lambda d}\, &\int
	\sup_\omega \ \lp{ S^\omega_{\lambda , \bullet < \min\{d,c\}} 
	(S_{1,d}u\, \nabla S_{\lambda,
	\bullet < \min\{d,c\}} v)(t)}{L^\infty_x}\, 
	dt \\ 
	&\hspace{2.5in} \lesssim \
	\lp{u}{F_{\Omega,1}}\cdot \ \lambda \ \lp{v}{F_{\Omega,\lambda}} \ ,
   \end{split}\label{LH_L1Linf_no_Omega}
	\\ \notag \\ \hline \hline \notag 
\end{align}\\

It remains to deal with the last term on the right hand side of 
\eqref{fine_low_high_cdist}. For this we'll prove the following two estimates 
involving the $Z_1$ and $Z_{\Omega,1}$ spaces:\\

\begin{align}\hline \hline \notag \\
	\lp{\sum_d \ S_{\lambda , \bullet < \min\{d,c\}} 
	(S_{1,d}u\, 
	\nabla S_{\lambda,\bullet < \min\{d,c\}} v)}{L^1(L^2)} \ &\lesssim \
	\lp{u}{Z_1}\cdot \ \lambda \ \lp{v}{F_{\Omega,\lambda}} \ ,
	\label{special_L1L2_omega_on_low} \\
	\lp{\sum_d \ S_{\lambda , \bullet < \min\{d,c\}} 
	(S_{1,d}u\, 
	\nabla S_{\lambda,\bullet < \min\{d,c\}} v)}{L^1(L^2)} \ &\lesssim \
	\lp{u}{Z_{\Omega,1}}\cdot \ \lambda \ \lp{v}{F_{\lambda}} \ .
	\label{special_L1L2_omega_on_high}
	\\ \notag \\ \hline \hline \notag 
\end{align}

\ret\ret

\section{Inductive Estimates \textbf{I}. $High\times High$ 
Frequency Interactions}

We now begin the proof of the boxed estimates in section 
\ref{HH_section}. To streamline the process, we will list out the 
individual steps in the proofs all in a row with very little text
explanations in between. At the beginning of each block of estimates,
we list in order of use
the various Sobolev, Strichartz, and multiplier boundedness
estimates that are used in the
lines that follow. It would have been more convenient to include these
references on the lines where they are used, but there does not seem to be 
enough room in the typesetting to allow this.\\

\begin{proof}[proof of estimate \eqref{extra_LinfL2_est}]
For a fixed $\sigma$ and $\mu$ we use the multiplier boundedness Lemma
\ref{generic_mult_lemma}, the localized Sobolev estimate \eqref{L1_loc_sob},
H\"olders inequality, and the energy estimate \eqref{F_energy_str_est} 
to compute that:
\begin{align}
	&\lp{S_{\sigma,\sigma}\, 
	\frac{1}{\varXi}\, P_\mu
	(S_1 u\, \nabla S_1 v)}{L^\infty(L^2)} \ , \notag \\
	\lesssim \ \ \ \sigma^{-2}\ &\lp{P_\mu
	(S_1 u\, \nabla S_1 v)}{L^\infty(L^2)} \ , \notag \\
	\lesssim \ \ \left(\frac{\mu}{\sigma}\right)^2
	&\lp{S_1 u\, \nabla S_1 v}{L^\infty(L^1)} \ , \notag \\
	\lesssim \  \left(\frac{\mu}{\sigma}\right)^2\
	&\lp{S_1 u}{L^\infty(L^2)}\cdot
	\lp{\nabla S_1 v}{L^\infty(L^2)} \ , \notag \\
	\lesssim \ \left(\frac{\mu}{\sigma}\right)^2\ &\lp{u}{F_1}\cdot
	\lp{v}{F_1} \ . \notag
\end{align}
Multiplying this last expression by $\mu$ and then summing yields:
\begin{align}
	\text{(L.H.S.)}\eqref{extra_LinfL2_est} \ &\lesssim \
	\sum_{\mu \lesssim 1}\ \mu \ \sum_{\mu < \sigma} \
	\left(\frac{\mu}{\sigma}\right)^2\ \lp{u}{F_1}\cdot
	\lp{v}{F_1} \ , \notag \\
	&\lesssim \ \lp{u}{F_1}\cdot \lp{v}{F_1} \ . \notag
\end{align}
\end{proof}


\begin{proof}[proof of estimates 
\eqref{HH_L1L2_no_omega1}--\eqref{HH_L1L2_no_omega2}]
Due to the boundedness of the multiplier $S_{1,\bullet < c\mu }$ on the 
$F$ spaces, and the ability to trade the $\nabla$ weight between either
term of the product, it suffices to prove the first estimate 
\eqref{HH_L1L2_no_omega1}. We begin by noting that:
\begin{equation}
	\text{(L.H.S.)}\eqref{HH_L1L2_no_omega1} \ \lesssim \
	\sum_{\substack{\sigma,\mu \ : \\
	\sigma \leqslant \mu \lesssim 1}} \ \lp{P_\sigma S_\mu (S_1 u\, \nabla 
	S_{1,c\mu \leqslant \bullet} v)}{L^1(L^2)} \ . \notag
\end{equation}
We now fix both $\sigma$ and $\mu$ for the moment and use the local Sobolev
estimate \eqref{L3/2_loc_sob}, H\"olders inequality, the Strichartz
estimate \eqref{F_L2L6_str_est}, and the $L^2$ estimate \eqref{sum_L2_F_est} 
to compute that:
\begin{align}
	& & &\lp{P_\sigma S_\mu (S_1 u\, \nabla 
	S_{1,c\mu \leqslant \bullet} v)}{L^1(L^2)} \ , \notag \\ 
	&\lesssim&  
	\sigma^\frac{2}{3} \ &\lp{S_1 u\, \nabla 
	S_{1,c\mu \leqslant \bullet} v}{L^1(L^\frac{3}{2})} \ , \notag \\
	&\lesssim& c^{-\frac{1}{2}}\sigma^\frac{1}{6}
	\left(\frac{\sigma}{\mu}\right)^\frac{1}{2}\ &\lp{S_1 u}{L^2(L^6)}
	\cdot\, (c\mu)^\frac{1}{2}
	\lp{\nabla S_{1,c\mu \leqslant \bullet} v}{L^2(L^2)} \ , \notag \\
	&\lesssim&  c^{-\frac{1}{2}}\sigma^\frac{1}{6}
	\left(\frac{\sigma}{\mu}\right)^\frac{1}{2}\ &\lp{u}{F_1}\cdot
	\lp{v}{F_1} \ . \notag
\end{align} 
This last expression can now be summed over $\sigma$ and $\mu$ to yield:
\begin{align}
	\text{(L.H.S.)}\eqref{HH_L1L2_no_omega1} \ &\lesssim \
	\sum_{\sigma \lesssim 1}\ \sum_{\sigma \leqslant \mu}\ 
	c^{-\frac{1}{2}}\sigma^\frac{1}{6}
	\left(\frac{\sigma}{\mu}\right)^\frac{1}{2}\ \lp{u}{F_1}\cdot
	\lp{v}{F_1} \ , \notag \\
	&\lesssim \ c^{-\frac{1}{2}}\ \lp{u}{F_1}\cdot\lp{v}{F_1} \ . \notag 
\end{align}
For a fixed small $c$, this yields the desired result.
\end{proof}\ret


\begin{proof}[proof of estimate \eqref{HH_L1L2_Omega_not_on_d}]
We begin by fixing both $\mu$ and $d$. Next, we use in order the angular 
decomposition \eqref{HH_angle_decomp}, the multiplier Lemma 
\ref{generic_mult_lemma} and orthogonality, the local Sobolev estimate
\eqref{L6/5_loc_sob}, the angular concentration estimate \eqref{L3_ang_const},
and finally the Strichartz estimate \eqref{F_L2L3_str_est} 
as well as the $L^2$ estimate \eqref{L2_F_est} to compute that:
\begin{align}
	& & &\int \lp{S_{\mu,\bullet \leqslant d} 
	(S_{1,\bullet \leqslant d}u\, 
	\nabla S_{1,d} v)(t)}{L^2_x}\, dt \ , \notag \\
	&=& 
	&\int   
	\lp{\sum_\omega S^{\pm\omega}_{\mu,\bullet \leqslant d} 
	(B^{-\omega}_{(\frac{d}{\mu})^\frac{1}{2}}
	S_{1,\bullet \leqslant d}u\, 
	\nabla B^{\omega}_{(\frac{d}{\mu})^\frac{1}{2}}
	S_{1,d} v)(t)}{L^2_x} \, dt \ , \notag \\
	&\lesssim& &\int \left(\sum_\omega
	\lp{B^{\pm\omega}_{(\frac{d}{\mu})^\frac{1}{2}}P_\mu 
	(B^{-\omega}_{(\frac{d}{\mu})^\frac{1}{2}}
	S_{1,\bullet \leqslant d}u\, 
	\nabla B^{\omega}_{(\frac{d}{\mu})^\frac{1}{2}}
	S_{1,d} v)(t)}{L^2_x}^2\right)^\frac{1}{2} \, dt \ , \notag \\
	&\lesssim& \mu^{\frac{5}{6}-} d^{\frac{1}{2}-}\ &\int
	\left(\sum_\omega  
	\lp{ (B^{-\omega}_{(\frac{d}{\mu})^\frac{1}{2}}
	S_{1,\bullet \leqslant d}u\, 
	\nabla B^{\omega}_{(\frac{d}{\mu})^\frac{1}{2}}
	S_{1,d} v)(t)}{L^{\frac{6}{5}+}_x}^2\right)^\frac{1}{2}\, dt 
	\ , \notag \\
	&\lesssim& \mu^{\frac{5}{6}-} d^{\frac{1}{2}-}\ &\int
	\sup_\omega\lp{ B^{-\omega}_{(\frac{d}{\mu})^\frac{1}{2}}
	S_{1,\bullet \leqslant d}u\, (t)}{L^{3+}_x} \cdot
	\left(\sum_\omega \lp{\nabla 
	B^{\omega}_{(\frac{d}{\mu})^\frac{1}{2}}
	S_{1,d} v\, (t)}{L^2_x}^2\right)^\frac{1}{2}\, dt \ , \notag \\
	&\lesssim&  \mu^{\frac{7}{12}-} d^{\frac{3}{4}-}\ &\int
	\lp{ S_{1,\bullet
	\leqslant d}\langle\Omega\rangle^\frac{1}{2} u\, (t)}{L^{3+}_x} \cdot
	\lp{\nabla S_{1,d} v\, (t)}{L^2_x}\, dt \ , \notag \\
	&\lesssim& \mu^{\frac{5}{6}-}\left(\frac{d}{\mu}\right)
	^{\frac{1}{4}-}\
	&\lp{S_{1,\bullet
	\leqslant d} \langle\Omega\rangle^\frac{1}{2}  u}{L^2(L^{3+})}\cdot\
	d^\frac{1}{2}\lp{\nabla S_{1,d} v}{L^2(L^2)} \ , \notag \\
	&\lesssim&  \mu^{\frac{5}{6}-}\left(\frac{d}{\mu}\right)
	^{\frac{1}{4}-}\
	&\lp{u}{F_{\Omega,1}}\cdot\lp{v}{F_{1}} \ . \notag
\end{align}
This last expression can now be safely summed over both $\mu$ and $d$ to
yield:
\begin{align}
	\text{(L.H.S.)}\eqref{HH_L1L2_Omega_not_on_d} \ &\lesssim \ 
	\sum_{\substack{\mu,d \ : \\
	\mu \lesssim 1 \, , \, d < c\mu }}
	\mu^{\frac{5}{6}-}\left(\frac{d}{\mu}\right)
	^{\frac{1}{4}-}\
	\lp{u}{F_{1}}\cdot \lp{v}{F_{\Omega,1}} \ , \notag \\
	&\lesssim \ \ \ \ \lp{u}{F_{\Omega,1}}\cdot \lp{v}{F_{1}} \ . \notag
\end{align}
\end{proof}\ret


\begin{proof}[proof of estimate \eqref{HH_L1L2_Omega_on_d}]
We begin by fixing both $\mu$ and $d$, we use in order the angular 
decomposition \eqref{HH_angle_decomp}, the multiplier Lemma 
\ref{generic_mult_lemma} and orthogonality, the local Sobolev estimate
\eqref{L3/2_loc_sob}, the angular concentration estimate \eqref{L2_ang_const},
and finally the Strichartz estimate \eqref{F_L2L6_str_est} 
as well as the $L^2$ estimate \eqref{L2_F_est} to compute that:
\begin{align}
	& & &\int \lp{S_{\mu,\bullet \leqslant d} 
	(S_{1,\bullet \leqslant d}u\, 
	\nabla S_{1,d} v)(t)}{L^2_x}\, dt \ , \notag \\
	&=& 
	&\int   
	\lp{\sum_\omega S^{\pm\omega}_{\mu,\bullet \leqslant d} 
	(B^{-\omega}_{(\frac{d}{\mu})^\frac{1}{2}}
	S_{1,\bullet \leqslant d}u\, 
	\nabla B^{\omega}_{(\frac{d}{\mu})^\frac{1}{2}}
	S_{1,d} v)(t)}{L^2_x} \, dt \ , \notag \\
	&\lesssim& &\int \left(\sum_\omega
	\lp{B^{\pm\omega}_{(\frac{d}{\mu})^\frac{1}{2}}P_\mu 
	(B^{-\omega}_{(\frac{d}{\mu})^\frac{1}{2}}
	S_{1,\bullet \leqslant d}u\, 
	\nabla B^{\omega}_{(\frac{d}{\mu})^\frac{1}{2}}
	S_{1,d} v)(t)}{L^2_x}^2\right)^\frac{1}{2} \, dt \ , \notag \\
	&\lesssim&  \mu^\frac{5}{12}d^\frac{1}{4}\ &\int
	\left(\sum_\omega   
	\lp{ (B^{-\omega}_{(\frac{d}{\mu})^\frac{1}{2}}
	S_{1,\bullet \leqslant d}u\, 
	\nabla B^{\omega}_{(\frac{d}{\mu})^\frac{1}{2}}
	S_{1,d} v)(t)}{L^\frac{3}{2}_x}^2\right)^\frac{1}{2}\, dt 
	\ , \notag \\
	&\lesssim&  \mu^\frac{5}{12}d^\frac{1}{4}\ &\int
	\left(\sum_\omega\lp{ B^{-\omega}_{(\frac{d}{\mu})^\frac{1}{2}}
	S_{1,\bullet \leqslant d}u\, (t)}{L^6_x}^2\right)^\frac{1}{2}\cdot
	\sup_\omega \lp{\nabla B^{\omega}_{(\frac{d}{\mu})^\frac{1}{2}}
	S_{1,d} v\, (t)}{L^2_x}\, dt \ , \notag \\
	&\lesssim& \ \mu^{-\frac{1}{12}}d^\frac{3}{4}\ &\int
	\left(\sum_\omega\lp{ B^{-\omega}_{(\frac{d}{\mu})^\frac{1}{2}}
	S_{1,\bullet\leqslant d}u\, (t)}{L^6_x}^2\right)^\frac{1}{2}\cdot
	\lp{\nabla S_{1,d} \langle\Omega\rangle 
	v\, (t)}{L^2_x}\, dt \ , \notag \\
	&\lesssim& \mu^\frac{1}{6}\left(\frac{d}{\mu}\right)^\frac{1}{4}\
	&\left(\sum_\omega\lp{ B^{-\omega}_{(\frac{d}{\mu})^\frac{1}{2}}
	S_{1,\bullet\leqslant d}u}{L^2(L^6)}^2\right)^\frac{1}{2}\cdot\
	d^\frac{1}{2}\lp{\nabla S_{1,d} \langle\Omega\rangle 
	v}{L^2(L^2)} \ , \notag \\
	&\lesssim&   \mu^\frac{1}{6}\left(\frac{d}{\mu}\right)^\frac{1}{4}\
	&\lp{u}{F_1}\cdot\lp{v}{F_{\Omega,1}} \ . \notag
\end{align}
This last expression can now be safely summed over both $\mu$ and $d$ to
yield:
\begin{align}
	\text{(L.H.S.)}\eqref{HH_L1L2_Omega_on_d} \ &\lesssim \ 
	\sum_{\substack{\mu,d \ : \\
	\mu \lesssim 1 \, , \, d < c\mu }}
	\mu^\frac{1}{6} \left(\frac{d}{\mu}\right)^\frac{1}{4}
	\lp{u}{F_{1}}\cdot \lp{v}{F_{\Omega,1}} \ , \notag \\
	&\lesssim \ \ \ \ \lp{u}{F_{1}}\cdot \lp{v}{F_{\Omega,1}} \ . \notag
\end{align}
\end{proof}\ret


\begin{proof}[proof of estimate \eqref{HH_L2_est}]
We begin by fixing both $\mu$ and $d$, we use in order the angular 
decomposition \eqref{HH_angle_decomp} and orthogonality, the local Sobolev 
estimate \eqref{L3/2_loc_sob}, the Strichartz estimates 
\eqref{F_L2L6_str_est} and \eqref{F_energy_str_est},
and finally the angular concentration estimate \eqref{F_ang_const}
to compute that:
\begin{align}
	& &  
	&\lp{S_{\mu,d} (S_{1,\bullet < \min\{d,c\mu\}}u\, 
	\nabla S_{1,\bullet < \min\{d,c\mu\}} v)}{L^2(L^2)} \ , \notag \\
	&\lesssim& 
	&\left(\sum_\omega \lp{S^{\pm\omega}_{\mu,d} (
	B^{-\omega}_{(\frac{d}{\mu})^\frac{1}{2}}
	S_{1,\bullet < \min\{d,c\mu\}}u\, 
	\nabla B^\omega_{(\frac{d}{\mu})^\frac{1}{2}}
	S_{1,\bullet < \min\{d,c\mu\}} v)}{L^2(L^2)}^2\right)^\frac{1}{2}
	\ , \notag \\
	&\lesssim&  \mu^\frac{5}{12}d^\frac{1}{4}\
	&\left(\sum_\omega \lp{
	B^{-\omega}_{(\frac{d}{\mu})^\frac{1}{2}}
	S_{1,\bullet < \min\{d,c\mu\}}u\, 
	\nabla B^\omega_{(\frac{d}{\mu})^\frac{1}{2}}
	S_{1,\bullet < \min\{d,c\mu\}} v}{L^2(L^\frac{3}{2})}^2
	\right)^\frac{1}{2}\ , \notag \\
	&\lesssim& \mu^\frac{5}{12}d^\frac{1}{4}\ 
	&\sup_\omega
	\lp{B^{-\omega}_{(\frac{d}{\mu})^\frac{1}{2}}
	S_{1,\bullet < \min\{d,c\mu\}}u}{L^2(L^6)}\ \notag \\
	& & &\hspace{1in}
	\cdot \ \left(\sum_\omega \lp{\nabla 
	B^\omega_{(\frac{d}{\mu})^\frac{1}{2}}
	S_{1,\bullet < \min\{d,c\mu\}} v}{L^\infty(L^2)}^2\right)^\frac{1}{2}
	\ , \notag \\
	&\lesssim& \mu^\frac{5}{12}d^\frac{1}{4}\ &\sup_\omega
	\lp{B^{-\omega}_{(\frac{d}{\mu})^\frac{1}{2}} u}{F_1}\cdot
	\left(\sum_\omega \lp{\nabla 
	B^\omega_{(\frac{d}{\mu})^\frac{1}{2}} v}{F_1}^2
	\right)^\frac{1}{2} \ , \notag \\
	&\lesssim&  \mu^{-\frac{1}{12}}d^\frac{3}{4}\ 
	&\lp{u}{F_{\Omega,1}}\cdot \lp{v}{F_1} \ . \notag
\end{align}
Multiplying this last expression through by $d^{-\frac{1}{2}}$ and
summing over $d$ and $\mu$ yields:
\begin{align}
	\text{(L.H.S.)}\eqref{HH_L2_est} \ &\lesssim \ 
	\sum_{\substack{\mu,d \ : \\
	\mu \lesssim 1 \, , \, d \leqslant \mu }}
	\mu^\frac{1}{6} \left(\frac{d}{\mu}\right)^\frac{1}{4}
	\lp{u}{F_{\Omega,1}}\cdot \lp{v}{F_1} \ , \notag \\
	&\lesssim \ \ \ \ \lp{u}{F_{\Omega,1}}\cdot \lp{v}{F_1} \ . \notag
\end{align}
\end{proof}\ret


\begin{proof}[proof of estimate \eqref{HH_L1Linf_with_Omega}]
This estimate is truly bilinear in nature, essentially due to the failure
of the $L^2(L^3)$ endpoint version of \eqref{ang_str_est}. 
Therefore we must proceed
in a more detailed fashion than the rest of the estimates in this section.
For a fixed $\mu$ and $d$, we use the partial angular decomposition 
\eqref{HH_angle_partial_decomp} 
in conjunction with the special multiplier bound 
\eqref{sqsum_L1Linf_bound_no_box}, and
the local Sobolev estimate \eqref{LinftL2_loc_sob}
to compute that:
\begin{align}
	& & &\int
	\left( \sum_\omega \lp{S^\omega_{\mu,d} 
	(S_{1,\bullet < \min\{d,c\mu\}}u\, 
	\nabla S_{1,\bullet < \min\{d,c\mu\}} v)(t)}
	{L^\infty_x}^2\right)^\frac{1}{2} \, dt 
	\ , \label{HH_L1Linf_with_Omega_start} \\
	&\lesssim& &\int
	\left( \sum_\omega \lp{B^\omega_{(\frac{d}{\mu})^\frac{1}{2}}P_\mu 
	(S_{1,\bullet < \min\{d,c\mu\}}u\, 
	\nabla B^\omega_{(\frac{d}{\mu})^\frac{1}{2}}
	S_{1,\bullet < \min\{d,c\mu\}} v)(t)}
	{L^\infty_x}^2\right)^\frac{1}{2} \, dt 
	\ , \notag \\
	&\lesssim&  \mu^{\frac{5}{4}-}d^{\frac{3}{4}-}
	\, &\int
	\left( \sum_\omega \lp{ 
	(S_{1,\bullet < \min\{d,c\mu\}}u\, 
	\nabla B^\omega_{(\frac{d}{\mu})^\frac{1}{2}}
	S_{1,\bullet < \min\{d,c\mu\}} v)(t)}
	{L^{2+}_x}^2\right)^\frac{1}{2} \, dt \ . \notag 
\end{align}
We now need to incorporate the ``improved'' Strichartz estimate
\eqref{F_imp_L2L3_str_est}. To do this, we begin by fixing
$t$ and use the dual--scale Sobolev estimate
from  \cite{TKR} along with H\"olders inequality to compute that:
\begin{align}
	&\ & &\lp{ (S_{1,\bullet < \min\{d,c\mu\}}u\, 
	\nabla B^\omega_{(\frac{d}{\mu})^\frac{1}{2}}
	S_{1,\bullet < \min\{d,c\mu\}} v)(t)}
	{L^{2+}_x} \ , \notag \\
	&\lesssim& 
	\mu^{2+} \ &\left( \sum_{Q_\alpha} \lp{ 
	S_{1,\bullet < \min\{d,c\mu\}}u\, 
	\nabla B^\omega_{(\frac{d}{\mu})^\frac{1}{2}}
	S_{1,\bullet < \min\{d,c\mu\}} v)(t)}
	{L^1(Q_\alpha)}^{2+}\right)^\frac{1}{2+} \ , \notag \\
	&\lesssim& \mu^{2+}\ &\left( \sum_{Q_\alpha}\lp{ 
	S_{1,\bullet < \min\{d,c\mu\}}u\, (t)}
	{L^2(Q_\alpha)}^{3+}\right)^\frac{1}{3+}\ \cdot \notag \\
	&\ & &\ \ \ \ \ \ \ \ \ \ \ \ \ \ \ \ \ \ \  
	\left( \sum_{Q_\alpha}\lp{
	\nabla B^\omega_{(\frac{d}{\mu})^\frac{1}{2}}
	S_{1,\bullet < \min\{d,c\mu\}} v\, (t)}
	{L^2(Q_\alpha)}^{6}\right)^\frac{1}{6} \ . \notag
\end{align}
We now square sum this last expression with respect to $\omega$ and
integrate with respect to time, followed by a use of H\"olders
inequality and the Strichartz estimates 
\eqref{F_imp_L2L3_str_est} and \eqref{imp_L2L6_str_est} 
to compute that:

\begin{align}
	&\ & \hbox{(L.H.S.)}\eqref{HH_L1Linf_with_Omega_start} \ , \notag \\
	&\lesssim& 
	\mu^{\frac{13}{4}-}d^{\frac{3}{4}-} \, &\int
	\left( \sum_{Q_\alpha}\lp{ 
	S_{1,\bullet < \min\{d,c\mu\}}u\, (t)}
	{L^2(Q_\alpha)}^{3+}\right)^\frac{1}{3+}\ \cdot \notag \\
	&\ & & \ \ \ \ \ \ \ \ 
	\left(\sum_\omega
	\left( \sum_{Q_\alpha}\lp{
	\nabla B^\omega_{(\frac{d}{\mu})^\frac{1}{2}}
	S_{1,\bullet < \min\{d,c\mu\}} v\, (t)}
	{L^2(Q_\alpha)}^{6}\right)^\frac{1}{3}\right)^\frac{1}{2}\, dt
	\ , \notag \\
	&\lesssim& \mu^{\frac{13}{4}-}d^{\frac{3}{4}-}\ 
	&\lp{ \left( \sum_{Q_\alpha}\lp{ 
	S_{1,\bullet < \min\{d,c\mu\}}u\, (t)}
	{L^2(Q_\alpha)}^{3+}\right)^\frac{1}{3+} }{L_t^2} \ \cdot \notag \\
	&\ & &\ \ \ \ \ \ \ \ 
	\left(\sum_\omega\
	\lp{ \left( \sum_{Q_\alpha}\lp{
	\nabla B^\omega_{(\frac{d}{\mu})^\frac{1}{2}}
	S_{1,\bullet < \min\{d,c\mu\}} v\, (t)}
	{L^2(Q_\alpha)}^{6}\right)^\frac{1}{6} }{L^2_t}^2 
	\right)^\frac{1}{2} \ , \notag \\
	&\lesssim&  \mu^{\frac{7}{4}-}d^{\frac{3}{4}-}\
	&\lp{u}{F_{\Omega,1}} \left( \sum_\omega\lp{
	B^\omega_{(\frac{d}{\mu})^\frac{1}{2}}
	v}{F_1}^2\right)^\frac{1}{2} \ ,  \notag \\
	&\lesssim&  \mu^{\frac{7}{4}-}d^{\frac{3}{4}-}\
	&\lp{u}{F_{\Omega,1}} \lp{v}{F_1} \ .  \notag
\end{align}
Multiplying this last expression by $\mu^{-\frac{3}{2}}d^{-\frac{1}{2}}$ 
and summing yields: 
\begin{align}
	\text{(L.H.S.)}\eqref{HH_L1Linf_with_Omega} \ &\lesssim \
	\sum_{\substack{\mu,d \ : \\
	\mu \lesssim 1 \, , \,
	d \leqslant \mu }} \mu^{\frac{1}{2}-} 
	\left(\frac{d}{\mu}\right)^{\frac{1}{4}-}\
	\lp{u}{F_{\Omega,1}} \lp{v}{F_1} \ , \notag \\
	&\lesssim \ \ \ \  \lp{u}{F_{\Omega,1}} \lp{v}{F_1} \ . \notag 
\end{align}
\end{proof}\ret


\begin{proof}[proof of estimates \eqref{HH_L1Linf_no_Omega}]
We begin by fixing both $\mu$ and $d$, and using the angular decomposition
\eqref{HH_angle_decomp} in conjunction with the special multiplier bound 
\eqref{sup_L1Linf_bound_no_box}, the local Sobolev estimate 
\eqref{LinftL3/2_loc_sob}, the angular concentration estimate
\eqref{L3_ang_const}, and finally the Strichartz estimate 
\eqref{high_F_L2L3_str_est} to compute that:
\begin{align}
	&\ & &\int
	\sup_\omega \lp{S^\omega_{\mu,d} (S_{1,\bullet < \min\{d,c\mu\}}u\, 
	\nabla S_{1,\bullet < \min\{d,c\mu\}} v)(t)}{L^\infty_x}
	\, dt \ , \notag \\
	&\lesssim&  &\int
	\sup_\omega \lp{B^\omega_{(\frac{d}{\mu})^\frac{1}{2}}P_\mu 
	(B^{-\omega}_{(\frac{d}{\mu})^\frac{1}{2}}
	S_{1,\bullet < \min\{d,c\mu\}}u\, 
	\nabla B^\omega_{(\frac{d}{\mu})^\frac{1}{2}} 
	S_{1,\bullet < \min\{d,c\mu\}} v)(t)}{L^\infty_x}
	\, dt \ , \notag \\
	&\lesssim&  \mu^{\frac{5}{3}-}d^{1-} \ &\int
	\sup_\omega \lp{(B^{-\omega}_{(\frac{d}{\mu})^\frac{1}{2}}
	S_{1,\bullet < \min\{d,c\mu\}}u\, 
	\nabla B^\omega_{(\frac{d}{\mu})^\frac{1}{2}}
	S_{1,\bullet < \min\{d,c\mu\}} v)(t)}{L^{\frac{3}{2}+}_x}
	\, dt \ , \notag \\
	&\lesssim&  \mu^{\frac{5}{3}-}d^{1-} \ &\int
	\sup_\omega \Big(\lp{B^{-\omega}_{(\frac{d}{\mu})^\frac{1}{2}}
	S_{1,\bullet < \min\{d,c\mu\}}u\, (t)}{L_x^{3+}} \notag \\
	&\ & & \ \ \ \ \ \ \ \ \ \ \ \ \ \ \ \ \ \ \ \ \
	\cdot \ \lp{\nabla B^{\omega}_{(\frac{d}{\mu})^\frac{1}{2}}
	S_{1,\bullet < \min\{d,c\mu\}} v\, (t)}{L^{3+}_x}\Big) \, dt
	\ , \notag \\
	&\lesssim&  \mu^{\frac{7}{6}-}d^{\frac{3}{2}-} \ &\int
	\lp{S_{1,\bullet < \min\{d,c\mu\}}\langle\Omega\rangle^\frac{1}{2}
	u\, (t)}{L_x^{3+}}\cdot  
	\lp{\nabla 
	S_{1,\bullet < \min\{d,c\mu\}} \langle\Omega\rangle^\frac{1}{2}
	v\, (t)}{L^{3+}_x} \, dt
	\ , \notag \\  
	&\lesssim&  \mu^{\frac{7}{6}-}d^{\frac{3}{2}-} \ 
	&\lp{S_{1,\bullet < \min\{d,c\mu\}}\langle\Omega\rangle^\frac{1}{2} u}
	{L^2(L^{3+})}\cdot
	\lp{\nabla S_{1,\bullet < \min\{d,c\mu\}} 
	\langle\Omega\rangle^\frac{1}{2}v}
	{L^2(L^{3+})} \ , \notag \\
	&\lesssim& \mu^{\frac{7}{6}-}d^{\frac{3}{2}-} \
	&\lp{u}{F_{\Omega,1}} \cdot \lp{v}{F_{\Omega,1}} \ . \notag
\end{align}
Multiplying this last expression through by $(\mu d)^{-1}$ and summing,
we get that:
\begin{align}
	\text{(L.H.S.)}\eqref{HH_L1Linf_no_Omega} \ &\lesssim \
	\sum_{\substack{\mu,d \ : \\
	\mu \lesssim 1 \, , \,
	d \leqslant \mu }}
	\mu^{\frac{2}{3}-}\left(\frac{d}{\mu}\right)^{\frac{1}{2}-}\
	\lp{u}{F_{\Omega,1}} \cdot \lp{v}{F_{\Omega,1}} \ , \notag \\
	&\lesssim \ \ \ \ \lp{u}{F_{\Omega,1}} \cdot \lp{v}{F_{\Omega,1}} 
	\ . \notag 
\end{align}
\end{proof}\ret

\ret\ret

\section{Inductive Estimates \textbf{II}. $Low\times High$ 
Frequency Interactions}

\begin{proof}[proof of estimate \eqref{LH_L1L2_no_omega1}]
The proof of this estimate follows directly from 
the boundedness of the $S_\lambda$ multiplier, H\"olders inequality,
and the Strichartz estimate \eqref{F_L2Linf_str_est} as well as the 
$L^2$ estimate \eqref{sum_L2_F_est}:
\begin{align}
	&\lp{S_\lambda(S_{1} u \, \nabla S_{\lambda,
        c \leqslant \bullet } v)}{L^1(L^2)} \ , \notag \\
	\lesssim \ c^{-\frac{1}{2}}\ 
	&\lp{S_{1} u}{L^2(L^\infty)}\cdot \ c^\frac{1}{2}\
	\lp{\nabla S_{\lambda,
        c \leqslant \bullet } v}{L^2(L^2)} \ , \notag \\
	\lesssim \ c^{-\frac{1}{2}}\ &\lp{u}{F_1}\cdot
	\ \lambda \ \lp{v}{F_{\lambda}} \ . \notag
\end{align}
For a fixed small $c$, this yields the desired result.
\end{proof}\ret


\begin{proof}[proof of estimate \eqref{LH_L1L2_no_omega2}]
The proof here is a simple matter of H\"olders inequality and the
Strichartz estimates \eqref{F_L2Linf_str_est} and 
\eqref{high_F_energy_str_est}:
\begin{align}
	&\ & \sum_{c \leqslant d}\ d^{-\frac{1}{2}}\ 
        &\lp{S_{\lambda,d}\, (S_{1} u \, \nabla 
	S_{\lambda,\bullet < c } v)}{L^2(L^2)} \ , \notag \\
	&\lesssim&  \sum_{c \leqslant d}\ d^{-\frac{1}{2}}\
	&\lp{ S_{1} u}{L^2(L^\infty)}\cdot
	\lp{\nabla S_{\lambda,\bullet < c } v}{L^\infty(L^2)} \ , \notag \\
	&\lesssim& \ \sum_{c \leqslant d}\ d^{-\frac{1}{2}}\
	&\lp{u}{F_1}\cdot \ \lambda \ \lp{v}{F_\lambda} \ , \notag \\
	&\lesssim& \text{ln}(c)\ 
	&\lp{u}{F_1}\cdot \ \lambda \ \lp{v}{F_\lambda} \ . \notag
\end{align}
For a fixed small $c$, this yields the desired result.
\end{proof}\ret


\begin{proof}[proof of estimate \eqref{spec_LH_L1Linf_with_Omega}]
For a fixed $d$ we use the decomposition \eqref{th_LH_small_angle_decomp} 
in conjunction with the multiplier bound \eqref{sqsum_L1Linf_bound_no_box}, 
H\"olders inequality, the local Sobolev estimate 
\eqref{low_LinftL6_loc_sob}, and the Strichartz
estimates \eqref{F_L2Linf_str_est} and \eqref{high_F_L2L6_str_est} 
to prove that:
\begin{align}
	& & &\int \left( \sum_\omega \lp{S^\omega_{\lambda,d} (S_1 u\,
        \nabla S_{\lambda,\bullet < c } v)(t)}{L^\infty_x}^2
        \right)^\frac{1}{2}
         dt \ , \notag \\
	&\lesssim& &\int \left( \sum_{\substack{ \omega_1,\omega_2\ : \\
	|\omega_1 - \omega_2| \sim (\frac{d}{c\lambda})^\frac{1}{2}}} 
	\lp{S^{\omega_1}_{\lambda,d} (S_1 u\,
        \nabla B^{\omega_2}_{(\frac{d}{c\lambda})^\frac{1}{2}}
	S_{\lambda,\bullet < c } v)(t)}{L^\infty_x}^2
        \right)^\frac{1}{2} dt \ , \notag \\
	&\lesssim& c^{-\frac{3}{4}} \ &\int
	\left( \sum_\omega 
	\lp{(S_1 u\,
        \nabla B^{\omega}_{(\frac{d}{c\lambda})^\frac{1}{2}}
	S_{\lambda,\bullet < c } v)(t)}{L^\infty_x}^2
        \right)^\frac{1}{2} dt \ , \notag \\
	&\lesssim& c^{-\frac{3}{4}} \ &\int \lp{S_1 u\, (t)}{L^\infty_x}
	\cdot \left( \sum_\omega \lp{\nabla 
	B^{\omega}_{(\frac{d}{c\lambda})^\frac{1}{2}}
	S_{\lambda,\bullet < c } v)(t)}{L^\infty_x}^2\right)^\frac{1}{2}
	dt \ , \notag \\
	&\lesssim& c^{-1} \lambda^\frac{5}{12}d^\frac{1}{4}
	\ &\lp{S_1 u}{L^2(L^\infty)}\cdot
	\left( \sum_\omega \lp{\nabla 
	B^{\omega}_{(\frac{d}{c\lambda})^\frac{1}{2}}
	S_{\lambda,\bullet < c } v}{L^2(L^6)}^2\right)^\frac{1}{2}
	\ , \notag \\
	&\lesssim& c^{-1} \lambda^\frac{15}{12}d^\frac{1}{4}
	\ &\lp{u}{F_1}\cdot \ \lambda\,
	\left( \sum_\omega \lp{ B^{\omega}_{(\frac{d}{c\lambda})^\frac{1}{2}}
	u}{F_\lambda} \right)^\frac{1}{2} \ , \notag \\
	&\lesssim& c^{-1}\lambda^\frac{15}{12}d^\frac{1}{4}\
	&\lp{u}{F_1}\cdot \ \lambda \ \lp{v}{F_\lambda} \ . \notag
\end{align}
Multiplying the last line above by the quantity $\lambda^{-\frac{3}{2}}
d^{-\frac{1}{2}}$, using the fact that $1\ll\lambda$,  and summing 
over $d$ yields:
\begin{align}
	\text{(L.H.S.)}\eqref{spec_LH_L1Linf_with_Omega} \ &\lesssim \
	\sum_{c\leqslant d} \ c^{-1}\lambda^{-\frac{1}{4}}d^{-\frac{1}{2}}\
	\lp{u}{F_1}\cdot \ \lambda \ \lp{v}{F_\lambda} \ , \notag \\
	&\lesssim \ c^{-1}\, \text{ln}(c)\
	\lp{u}{F_1}\cdot \ \lambda \ \lp{v}{F_\lambda} \ . \notag
\end{align}
For a fixed small $c$, this yields the desired result.
\end{proof}\ret


\begin{proof}[proof of estimate \eqref{spec_LH_L1Linf_no_Omega}]
For a fixed $d$ we use in order
the decomposition \eqref{th_LH_small_angle_decomp},
the multiplier bound \eqref{sup_L1Linf_bound_no_box}, the local Sobolev
estimate \eqref{low_LinftL3_loc_sob}, the angular concentration
estimate \eqref{L3_ang_const}, and the Strichartz estimates 
\eqref{F_L2Linf_str_est} and \eqref{high_F_L2L3_str_est}
to compute that:
\begin{align}
	& & &\int \sup_\omega \lp{S^\omega_{\lambda,d} (S_1 u\,
        \nabla S_{\lambda,\bullet < c } v)(t)}{L^\infty_x}\, dt \ , \notag \\
	&\lesssim& &\int \sup_{\substack{ \omega_1,\omega_2\ : \\
	|\omega_1 - \omega_2| \sim (\frac{d}{c\lambda})^\frac{1}{2}}} 
	\lp{S^{\omega_1}_{\lambda,d} (S_1 u\,
        \nabla B^{\omega_2}_{(\frac{d}{c\lambda})^\frac{1}{2}}
	S_{\lambda,\bullet < c } v)(t)}{L^\infty_x}\, dt \ , \notag \\
	&\lesssim& &\int \sup_\omega \lp{(S_1 u\,
        \nabla B^{\omega}_{(\frac{d}{c\lambda})^\frac{1}{2}}
	S_{\lambda,\bullet < c } v)(t)}{L^\infty_x}\, dt \ , \notag \\
	&\lesssim& &\int \lp{S_1 u\, (t)}{L^\infty_x}\cdot
        \ \sup_\omega\lp{\nabla B^{\omega}_{(\frac{d}{c\lambda})^\frac{1}{2}}
	S_{\lambda,\bullet < c } v\, (t)}{L^\infty_x}\, dt \ , \notag \\
	&\lesssim& c^{(-\frac{1}{2})+}\lambda^{\frac{5}{6}-}d^{\frac{1}{2}-}\
	&\int \lp{S_1 u\, (t)}{L^\infty_x}\cdot
        \ \sup_\omega\lp{\nabla B^{\omega}_{(\frac{d}{c\lambda})^\frac{1}{2}}
	S_{\lambda,\bullet < c } v\, (t)}{L^{3+}_x}\, dt \ , \notag \\
	&\lesssim& c^{(-\frac{3}{4})+} \lambda^{\frac{7}{12}-}d^{\frac{3}{4}-}\
	&\int \lp{S_1 u\, (t)}{L^\infty_x}\cdot
        \lp{\nabla \langle\Omega\rangle^\frac{1}{2}
	S_{\lambda,\bullet < c } v\, (t)}{L^{3+}_x}\, dt \ , \notag \\
	&\lesssim& c^{(-\frac{3}{4})+} \lambda^{\frac{7}{12}-}d^{\frac{3}{4}-}\
	&\lp{S_1 u}{L^2(L^\infty)}\cdot
        \lp{\nabla \langle\Omega\rangle^\frac{1}{2}
	S_{\lambda,\bullet < c } v}{L^2(L^{3+})} \ , \notag \\
	&\lesssim& c^{(-\frac{3}{4})+} \lambda^{\frac{3}{4}+}d^{\frac{3}{4}-}\
	&\lp{u}{F_1}\cdot\lambda\ \lp{v}{F_{\Omega,\lambda}} \ . \notag
\end{align}
Multiplying this last line through by $(\lambda d)^{-1}$ and 
summing, while using the fact that $1 \ll \lambda$ yields:
\begin{align}
        \text{(L.H.S.)}\eqref{spec_LH_L1Linf_no_Omega} \ &\lesssim \
	\sum_{c\leqslant d} \ c^{(-\frac{3}{4})+} \lambda^{(-\frac{1}{4})+}
	d^{-(\frac{1}{4}+)}\ \lp{u}{F_1}\cdot\lambda\ 
	\lp{v}{F_{\Omega,\lambda}} \ , \notag \\
	&\lesssim \ c^{-\frac{3}{4}}\
	\lp{u}{F_1}\cdot\lambda\ \lp{v}{F_{\Omega,\lambda}} \ . \notag
\end{align}
For a fixed $c$, this yields the desired result.
\end{proof}\ret


\begin{proof}[proof of estimate \eqref{LH_L1L2_Omega_on_low}]
For a fixed $d$, we use the angular decomposition \eqref{LH_wide_angle_decomp},
H\"olders inequality, the local Sobolev estimate \eqref{LinftL6_loc_sob} 
and the concentration estimate \eqref{L2_ang_const}, and finally the 
Strichartz estimate \eqref{F_L2L6_str_est} and the $L^2$ estimate
\eqref{L2_F_est} to compute that:
\begin{align}
	&\ & &\int \lp{S_{\lambda,\bullet \leqslant d} 
        (S_{1,\bullet \leqslant d}u\, \nabla S_{\lambda,d} v)(t)}{L^2_x}
	\ \, \notag \\
	&=& &\int  \left(\sum_\omega  
	\lp{B^\omega_{d^\frac{1}{2}}P_\lambda
        (S^{\pm\omega}_{1,\bullet \leqslant d}u\cdot
	\nabla B^\omega_{d^\frac{1}{2}}S_{\lambda,d} v)(t)}{L^2_x}^2
	\right)^\frac{1}{2} dt \ , \notag \\
	&\lesssim& &\int  \left(\sum_\omega  
	\lp{S^{\pm\omega}_{1,\bullet \leqslant d}u\, (t)}{L^\infty}^2
	\right)^\frac{1}{2} \cdot \ \sup_\omega \lp{
	\nabla B^\omega_{d^\frac{1}{2}}S_{\lambda,d} v\, (t)}{L^2_x}
	\, dt \ , \notag \\
	&\lesssim& d^\frac{3}{4} \ &\int  \left(\sum_\omega  
	\lp{S^{\pm\omega}_{1,\bullet \leqslant d}u\, (t)}{L^6}^2
	\right)^\frac{1}{2} \cdot \lp{
	\nabla S_{\lambda,d} \langle\Omega\rangle v\, (t)}{L^2_x}
	\, dt \ , \notag \\
	&\lesssim& d^\frac{1}{4}\ &\left(\sum_\omega\
	\lp{S^{\pm\omega}_{1,\bullet \leqslant d}u}{L^2(L^6)}^2
	\right)^\frac{1}{2}\cdot \ d^\frac{1}{2}\ 
	\lp{\nabla S_{\lambda,d} \langle\Omega\rangle v}{L^2(L^2)} \ , \notag \\
	&\lesssim&  d^\frac{1}{4}\ &\lp{u}{F_1}\cdot \ \lambda \ 
	\lp{v}{F_{\Omega,\lambda}} \ . \notag
\end{align}
This last expression can now be safely summed over $d$ to yield:
\begin{align}
	\text{(L.H.S.)}\eqref{LH_L1L2_Omega_on_low} \ &\lesssim \
	\sum_{d \leqslant 1}\ d^\frac{1}{4}\ \lp{u}{F_1}\cdot \ \lambda \ 
	\lp{v}{F_{\Omega,\lambda}} \ , \notag \\
	&\lesssim \ \lp{u}{F_1}\cdot \ \lambda \ 
	\lp{v}{F_{\Omega,\lambda}} \ . \notag
\end{align}
\end{proof}


\begin{proof}[proof of estimate \eqref{LH_L1L2_Omega_on_high}]
Here, for a fixed $d$, we use in order the angular decomposition 
\eqref{LH_wide_angle_decomp},
H\"olders inequality, the local Sobolev estimate \eqref{LinftL3_loc_sob},
the concentration estimate \eqref{L3_ang_const}, and the 
Strichartz estimate \eqref{F_L2L3_str_est} as well as the $L^2$ estimate
\eqref{L2_F_est} to compute that: 
\begin{align}
	&\ & &\int \lp{S_{\lambda,\bullet \leqslant d} 
        (S_{1,\bullet \leqslant d}u\cdot \nabla S_{\lambda,d} v)(t)}{L^2_x}
	\ \, \notag \\
	&=& &\int  \left(\sum_\omega  
	\lp{B^\omega_{d^\frac{1}{2}}P_\lambda
        (S^{\pm\omega}_{1,\bullet \leqslant d}u\cdot
	\nabla B^\omega_{d^\frac{1}{2}}S_{\lambda,d} v)(t)}{L^2_x}^2
	\right)^\frac{1}{2} dt \ , \notag \\
	&\lesssim& &\int  \sup_\omega  
	\lp{S^{\pm\omega}_{1,\bullet \leqslant d}u\, (t)}{L^\infty}
	\cdot \left( \sum_\omega \lp{
	\nabla B^\omega_{d^\frac{1}{2}}S_{\lambda,d} v\, (t)}{L^2_x}^2
	\right)^\frac{1}{2} \, dt \ , \notag \\
	&\lesssim& d^{\frac{1}{2}-} \ &\int  \sup_\omega  
	\lp{S^{\pm\omega}_{1,\bullet \leqslant d}u\, (t)}{L^{3+}}
	\cdot \lp{ \nabla S_{\lambda,d} v\, (t)}{L^2_x}
	\, dt \ , \notag \\
	&\lesssim&  d^{\frac{3}{4}-} \ &\int    
	\lp{S_{1,\bullet \leqslant d}\langle\Omega\rangle
	^\frac{1}{2}u\, (t)}{L^{3+}}
	\cdot \lp{ \nabla S_{\lambda,d} v\, (t)}{L^2_x}
	\, dt \ , \notag \\
	&\lesssim& d^{\frac{1}{4}-}\
	&\lp{S_{1,\bullet \leqslant d}\langle\Omega\rangle^\frac{1}{2}u 
	}{L^2(L^{3+})}\cdot \ d^\frac{1}{2} \ 
	\lp{\nabla S_{\lambda,d} v}{L^2(L^2)} \ , \notag\\
	&\lesssim&  d^{\frac{1}{4}-}\ &\lp{u}{F_{\Omega,1}}\cdot \ \lambda \ 
	\lp{v}{F_{\lambda}} \ . \notag
\end{align}
This last expression can now be safely summed over $d$ to yield:
\begin{align}
	\text{(L.H.S.)}\eqref{LH_L1L2_Omega_on_high} \ &\lesssim \
	\sum_{d \leqslant 1}\ d^{\frac{1}{4}-}\ \lp{u}{F_{\Omega,1}}
	\cdot \ \lambda \ \lp{v}{F_{\lambda}} \ , \notag \\
	&\lesssim \ \lp{u}{F_{\Omega,1}}\cdot \ \lambda \ 
	\lp{v}{F_\lambda} \ . \notag
\end{align}
\end{proof}\ret


\begin{proof}[proof of estimate \eqref{LH_L2L2_Omega_on_low}]
The proof here uses essentially the same steps as 
\eqref{LH_L1L2_Omega_on_low}
above. We begin by fixing $d$ and use the angular decomposition 
\eqref{LH_wide_angle_decomp},
H\"olders inequality and the local Sobolev estimate \eqref{LinftL6_loc_sob},
the Strichartz estimates  \eqref{F_L2L6_str_est} and 
\eqref{high_F_energy_str_est},
and the concentration estimate \eqref{F_ang_const} to compute that:
\begin{align}
	&\ & &\lp{S_{\lambda,d} (S_{1,\bullet \leqslant d}u\cdot 
	\nabla S_{\lambda,\bullet < d} v)}{L^2(L^2)} \ , \notag \\
	&\lesssim& &\bigg(\sum_\omega 
	\lp{ B^\omega_{d^\frac{1}{2}}
	S_{\lambda,d} (S^{\pm\omega}_{1,\bullet \leqslant d}u\cdot 
	\nabla B^\omega_{d^\frac{1}{2}} 
	S_{\lambda,\bullet < d} v)}{L^2(L^2)}^2\bigg)^\frac{1}{2} \ , \notag\\
	&\lesssim& d^{\frac{1}{4}}\ &\bigg(\sum_\omega 
	\lp{S^{\pm\omega}_{1,\bullet \leqslant d}u}{L^2(L^6)}^2
	\bigg)^\frac{1}{2}\cdot 
	\ \sup_\omega\lp{\nabla B^\omega_{d^\frac{1}{2}} 
	S_{\lambda,\bullet < d} v}{L^\infty(L^2)} \ , \notag \\
	&\lesssim&  d^{\frac{1}{4}}\ &\bigg(\sum_\omega 
	\lp{ B^\omega_{d^\frac{1}{2}} u}{F_1}^2\bigg)^\frac{1}{2}\cdot 
	\ \lambda \ \sup_\omega\lp{B^\omega_{d^\frac{1}{2}} 
	v}{F_\lambda} \ , \notag \\
	&\lesssim& 
	d^\frac{3}{4} \ &\lp{u}{F_1}\cdot \ \lambda \ 
	\lp{v}{F_{\Omega,\lambda}} \ . \notag
\end{align}
We can now multiply this last expression by $d^{-\frac{1}{2}}$ and sum
to yield:
\begin{align}
	\text{(L.H.S.)}\eqref{LH_L2L2_Omega_on_low} \ &\lesssim \
	\sum_{d \leqslant 1}\ d^\frac{1}{4}\ \lp{u}{F_1}\cdot \ \lambda \ 
	\lp{v}{F_{\Omega,\lambda}} \ , \notag \\
	&\lesssim \ \lp{u}{F_1}\cdot \ \lambda \ 
	\lp{v}{F_{\Omega,\lambda}} \ . \notag
\end{align}
\end{proof}\ret


\begin{proof}[proof of estimate \eqref{LH_L2L2_Omega_on_high}]
The proof here uses essentially the same steps as 
\eqref{LH_L1L2_Omega_on_high} above. For fixed  $d$ we use in order the 
angular decomposition \eqref{LH_wide_angle_decomp},
H\"olders inequality and the local Sobolev estimate \eqref{LinftL3_loc_sob},
the concentration estimate \eqref{L3_ang_const}, and the 
Strichartz estimates \eqref{F_L2L3_str_est} and 
\eqref{high_F_energy_str_est} to compute that:
\begin{align}
	&\ & &\lp{S_{\lambda,d} (S_{1,\bullet \leqslant d}u\cdot 
	\nabla S_{\lambda,\bullet < d} v)}{L^2(L^2)} \ , \notag \\
	&\lesssim& &\bigg(\sum_\omega 
	\lp{ B^\omega_{d^\frac{1}{2}}
	S_{\lambda,d} (S^{\pm\omega}_{1,\bullet \leqslant d}u\cdot 
	\nabla B^\omega_{d^\frac{1}{2}} 
	S_{\lambda,\bullet < d} v)}{L^2(L^2)}^2\bigg)^\frac{1}{2} \ , \notag\\
	&\lesssim& d^{\frac{1}{2}-}\ 
	&\sup_\omega\ \lp{S^{\pm\omega}_{1,\bullet \leqslant d}u}{L^2(L^{3+})}
	\cdot \bigg(\sum_\omega
	\ \lp{\nabla B^\omega_{d^\frac{1}{2}} 
	S_{\lambda,\bullet < d} v}{L^\infty(L^2)}^2\bigg)^\frac{1}{2} 
	\ , \notag \\
	&\lesssim& d^{\frac{3}{4}-}\ 
	&\lp{S^{\pm\omega}_{1,\bullet \leqslant d}\langle\Omega\rangle
	^\frac{1}{2}u}
	{L^2(L^{3+})}\cdot \bigg(\sum_\omega
	\ \lp{\nabla B^\omega_{d^\frac{1}{2}} 
	S_{\lambda,\bullet < d} v}{L^\infty(L^2)}^2\bigg)^\frac{1}{2} 
	\ , \notag \\
	&\lesssim& 
	d^{\frac{3}{4}-} \ &\lp{u}{F_{\Omega,1}}\cdot \ \lambda \ 
	\lp{v}{F_{\lambda}} \ . \notag
\end{align}
We can now multiply this last expression by $d^{-\frac{1}{2}}$ and sum
to yield:
\begin{align}
	\text{(L.H.S.)}\eqref{LH_L2L2_Omega_on_high} \ &\lesssim \
	\sum_{d \leqslant 1}\ d^{\frac{1}{4}-}\ 
	\lp{u}{F_{\Omega,1}}\cdot \ \lambda \ 
	\lp{v}{F_{\lambda}} \ , \notag \\
	&\lesssim \ \lp{u}{F_{\Omega,1}}\cdot \ \lambda \ 
	\lp{v}{F_{\lambda}} \ . \notag
\end{align}
\end{proof}\ret


\begin{proof}[proof of estimate \eqref{LH_L1Linf_Omega_on_low}]
For a fixed $d$ we use the angular decomposition 
\eqref{LH_small_angle_decomp} in conjunction
with the multiplier bound \eqref{spec_sqsum_L1Linf_bound_no_box},
the multiplier Lemma \ref{generic_mult_lemma} and H\"olders inequality,
the local Sobolev lemmas \eqref{LinftL6_loc_sob} and 
\eqref{low_LinftL3_loc_sob},
the angular concentration estimate \eqref{L3_ang_const}, 
and finally the Strichartz estimates \eqref{F_L2L6_str_est} and 
\eqref{high_F_L2L3_str_est} to compute that:
\begin{align}
	&\ & &\int
        \left( \sum_{\substack{\theta \ : \\
        |\theta|\sim \, d^\frac{1}{2}}} \ \sup_{\omega \subseteq \theta}
        \lp{ S^\omega_{\lambda , \bullet < \min\{d,c\}} 
        (S_{1,d}u\, \nabla S_{\lambda,
        \bullet < \min\{d,c\}} v)(t)}{L^\infty_x}^2
        \right)^\frac{1}{2} dt \ , \notag \\
	&=& &\int
        \left( \sum_{\substack{\theta \ : \\
        |\theta|\sim \, d^\frac{1}{2}}} \ \sup_{\omega \subseteq \theta}
        \lp{ B^\omega_{(\frac{d}{\lambda})^\frac{1}{2}}P_\lambda 
        (S^{\pm\theta}_{1,d}u\, \nabla S^\omega_{\lambda,
        \bullet < \min\{d,c\}} v)(t)}{L^\infty_x}^2
        \right)^\frac{1}{2} dt \ , \notag \\
	&\lesssim& &\int \left( \sum_\theta  \ 
        \lp{ S^{\pm\theta}_{1,d}u\,(t)}{L^\infty_x}^2\cdot\  
	\sup_{\omega \subseteq \theta}\lp{\nabla S^\omega_{\lambda,
        \bullet < \min\{d,c\}} v)(t)}{L^\infty_x}^2
        \right)^\frac{1}{2} dt \ , \notag \\
	&\lesssim& \lambda^{\frac{5}{6}-}d^{\frac{3}{4}-}\
	&\int \left( \sum_\theta  \ 
        \lp{ S^{\pm\theta}_{1,d}u\,(t)}{L^6_x}^2\right)^\frac{1}{2}
	\cdot\  
	\sup_{\omega }\lp{\nabla 
	S^\omega_{\lambda,
        \bullet < \min\{d,c\}} v\, (t)}{L^{3+}_x}\, dt \ , \notag \\
	&\lesssim& \lambda^{\frac{7}{12}-}d^{1-}\
	&\int \left( \sum_\theta  \ 
        \lp{ S^{\pm\theta}_{1,d}u\,(t)}{L^6_x}^2\right)^\frac{1}{2}
	\cdot \lp{\nabla 
	S_{\lambda,
        \bullet < \min\{d,c\}} \langle\Omega\rangle^\frac{1}{2}
	v\, (t)}{L^{3+}_x}\, dt \ , \notag \\
	&\lesssim& \lambda^{\frac{3}{4}+}d^{1-}\
	&\lp{u}{F_1}\cdot \ \lambda \ \lp{v}{F_{\Omega,\lambda}} \ . \notag
\end{align}
Multiplying this last expression through by $(\lambda d^\frac{1}{2})^{-1}$
and summing over $d$ yields:
\begin{align}
	\text{(L.H.S.)}\eqref{LH_L1Linf_Omega_on_low} \ &\lesssim \
	\sum_{d \leqslant 1}\ d^{\frac{1}{4}-}\left(\frac{d}{\lambda}
	\right)^{\frac{1}{4}-}\ 
	\lp{u}{F_{1}}\cdot \ \lambda \ 
	\lp{v}{F_{\Omega,\lambda}} \ , \notag \\
	&\lesssim \ \lp{u}{F_{1}}\cdot \ \lambda \ 
	\lp{v}{F_{\Omega,\lambda}} \ . \notag
\end{align}
\end{proof}\ret


\begin{proof}[proof of estimate \eqref{LH_L1Linf_Omega_on_high}]
We begin by fixing $d$ and use the angular decomposition 
\eqref{LH_small_angle_decomp} in conjunction with the multiplier
bound \eqref{sqsum_L1Linf_bound_no_box} and H\"olders inequality, 
the local Sobolev estimates \eqref{LinftL6_loc_sob} and
\eqref{low_LinftL3_loc_sob},
the angular concentration estimate \eqref{L3_ang_const}, and the Strichartz
estimates \eqref{F_L2L6_str_est} and \eqref{high_F_L2L3_str_est}
to compute that:
\begin{align}
	&\ & &\int
        \left( \sum_\omega \ \lp{ S^\omega_{\lambda , \bullet < \min\{d,c\}} 
        (S_{1,d}u\, \nabla S_{\lambda,
        \bullet < \min\{d,c\}} v)(t)}{L^\infty_x}^2
        \right)^\frac{1}{2} dt \ , \notag \\
	&\lesssim& 
	&\int \left(\sum_{\substack{\omega_1,\omega_2 \ :\\
        |\omega_1 \mp \omega_2|\sim d^\frac{1}{2}}}
	\lp{ S^{\pm\omega_2}_{1,d}u\,(t)}{L^\infty_x}^2\cdot 
	\lp{\nabla S^{\omega_1}_{\lambda,
        \bullet < \min\{d,c\}} v\,(t)}{L^\infty_x}^2
        \right)^\frac{1}{2} dt \ , \notag \\
	&\lesssim&  \lambda^\frac{5}{12}d^{\frac{3}{4}-}\
	&\int \sup_\omega
	\lp{ S^{\pm\omega}_{1,d}u\,(t)}{L^{3+}_x}\cdot 
	\left(\sum_\omega \lp{\nabla S^{\omega}_{\lambda,
        \bullet < \min\{d,c\}} v\,(t)}{L^6_x}^2
        \right)^\frac{1}{2} dt \ , \notag \\
	&\lesssim& \lambda^\frac{5}{12}d^{1-}\
	&\lp{S_{1,d}\langle\Omega\rangle^\frac{1}{2} u}{L^2(L^{3+})}\cdot 
	\left(\sum_\omega \lp{\nabla S^{\omega}_{\lambda,
        \bullet < \min\{d,c\}} v}{L^2(L^6)}^2
        \right)^\frac{1}{2} \ , \notag \\
	&\lesssim&  \lambda^\frac{15}{12}d^{1-}\
	&\lp{u}{F_{\Omega,1}}\cdot \ \lambda \ \lp{v}{F_\lambda} \ . \notag
\end{align}
Multiplying this last expression through by $(\lambda^\frac{3}{2} 
d^\frac{1}{2})^{-1}$ and summing over $d$ yields:
\begin{align}
	\text{(L.H.S.)}\eqref{LH_L1Linf_Omega_on_low} \ &\lesssim \
	\sum_{d \leqslant 1}\ d^{\frac{1}{4}-}\left(\frac{d}{\lambda}
	\right)^{\frac{1}{4}}\ 
	\lp{u}{F_{\Omega,1}}\cdot \ \lambda \ 
	\lp{v}{F_{\lambda}} \ , \notag \\
	&\lesssim \ \lp{u}{F_{\Omega,1}}\cdot \ \lambda \ 
	\lp{v}{F_{\lambda}} \ . \notag
\end{align}
\end{proof}\ret


\begin{proof}[proof of estimate \eqref{LH_L1Linf_no_Omega}]
Fixing $d$ we use in order the angular decomposition 
\eqref{LH_small_angle_decomp} in conjunction with the multiplier
bound \eqref{sup_L1Linf_bound_no_box} and H\"olders inequality,
the local Sobolev estimates \eqref{LinftL3_loc_sob} and 
\eqref{low_LinftL3_loc_sob}, the angular concentration
estimate \eqref{L3_ang_const}, and the Strichartz
estimates \eqref{F_L2L3_str_est} and \eqref{high_F_L2L3_str_est}
to compute that:
\begin{align}
	&\ & &\int \sup_\omega \ 
	\lp{ S^\omega_{\lambda , \bullet < \min\{d,c\}} 
        (S_{1,d}u\, \nabla S_{\lambda,
        \bullet < \min\{d,c\}} v)(t)}{L^\infty_x}\, dt \ , \notag \\
	&\lesssim& &\int \sup_\omega \ 
	\lp{ S^{\pm\omega}_{1,d}u\,(t)}{L^\infty_x}\cdot \
        \sup_\omega\ 
	\lp{ \nabla S^\omega_{\lambda,
        \bullet < \min\{d,c\}} v\, (t)}{L^\infty_x}\, dt \ , \notag \\
	&\lesssim& 
	\lambda^{\frac{5}{6}-} d^{1-}\
	&\int \sup_\omega \ 
	\lp{ S^{\pm\omega}_{1,d}u\,(t)}{L^{3+}_x}\cdot \
        \sup_\omega\ 
	\lp{ \nabla S^\omega_{\lambda,
        \bullet < \min\{d,c\}} v\, (t)}{L^{3+}_x}\, dt \ , \notag \\
	&\lesssim& \lambda^{\frac{7}{12}-}d^{\frac{3}{2}-}\
	&\lp{ S_{1,d}\langle\Omega\rangle^\frac{1}{2} u}{L^2(L^{3+})}\ \cdot \
	\lp{\nabla S_{\lambda,
        \bullet < \min\{d,c\}} \langle\Omega\rangle^\frac{1}{2}
	v}{L^2(L^{3+})} \ , \notag \\
	&\lesssim&  \lambda^{\frac{3}{4} + } d^{\frac{3}{2}-}\
	&\lp{u}{F_{\Omega,1}}\cdot \ \lambda \
        \lp{v}{F_{\Omega,\lambda}} \ . \notag
\end{align}
Multiplying this last expression through by $(\lambda d)^{-1}$ and
summing over yields: 
\begin{align}
	\text{(L.H.S.)}\eqref{LH_L1Linf_no_Omega} \ &\lesssim \
	\sum_{d \leqslant 1}\ d^{\frac{1}{4}-}\left(\frac{d}{\lambda}
	\right)^{\frac{1}{4}-}\ 
	\lp{u}{F_{\Omega,1}}\cdot \ \lambda \ 
	\lp{v}{F_{\Omega,\lambda}} \ , \notag \\
	&\lesssim \ \lp{u}{F_{\Omega,1}}\cdot \ \lambda \ 
	\lp{v}{F_{\Omega,\lambda}} \ . \notag
\end{align}
\end{proof}\ret


\begin{proof}[proof of estimate \eqref{special_L1L2_omega_on_low}]
For any decomposition $u=\sum_\alpha u^\alpha$, we
begin by fixing $d$ and using the angular decomposition 
\eqref{LH_wide_angle_decomp}
in conjunction with the multiplier Lemma \ref{generic_mult_lemma} 
and H\"olders inequality to compute that:
\begin{align}
	& & &\int \lp{S_{\lambda , \bullet < \min\{d,c\}} 
	(S_{1,d}u\cdot
	\nabla S_{\lambda,\bullet < \min\{d,c\}} v)(t)}{L^2_x}\, dt 
	\ , \label{special_L1L2_omega_on_low_start}\\
	&\lesssim& \sum_\alpha\ &\int \lp{S_{\lambda , \bullet < \min\{d,c\}} 
	(S_{1,d}u^\alpha\cdot
	\nabla S_{\lambda,\bullet < \min\{d,c\}} v)(t)}{L^2_x}\, dt 
	\ , \notag \\
	&\lesssim& \sum_\alpha &\int \left(\sum_\omega
	\lp{ B^\omega_{d^\frac{1}{2}} P_\lambda  
	(S^{\pm\omega}_{1,d}u^\alpha\cdot
	\nabla B^\omega_{d^\frac{1}{2}}
	S_{\lambda,\bullet < \min\{d,c\}} v)(t)}{L^2_x}^2\right)^\frac{1}{2}
	\, dt \ , \notag \\
	&\lesssim& \sum_\alpha &\int \left(\sum_\omega\
	\lp{S^{\pm\omega}_{1,d}u^\alpha\, (t)}{L^\infty_x}^2\cdot
	\lp{\nabla B^\omega_{d^\frac{1}{2}}
	S_{\lambda,\bullet < \min\{d,c\}} v\, (t)}{L^2_x}^2
	\right)^\frac{1}{2} dt \ . \notag
\end{align}
Now, for each fixed $\alpha$, we let 
$\{\theta_{\alpha,d}\}$ be a collection of angles such that 
$|\omega| \leqslant |\theta_{\alpha,d}|$. Then with a repeated use of H\"olders
inequality, the concentration estimate \eqref{L2_ang_const}, and the 
Strichartz estimate \eqref{high_F_energy_str_est}, we see that:
\begin{align}
	& & &\text{(L.H.S.)}\eqref{special_L1L2_omega_on_low_start}\text{\ for
        a fixed $\alpha$} \ , \notag\\
	&\lesssim& &\int \left(\sum_{\theta_{\alpha,d}}\ 
	\sup_{\omega\subseteq\theta_{\alpha,d}}\
	\lp{S^{\pm\omega}_{1,d}u^\alpha\, (t)}{L^\infty_x}^2
	\cdot\ \sum_{\omega\subseteq\theta_{\alpha,d}}\ 
	\lp{\nabla B^\omega_{d^\frac{1}{2}}
	S_{\lambda,\bullet < \min\{d,c\}} v\, (t)}{L^2_x}^2
	\right)^\frac{1}{2} dt \ , \notag \\
	&\lesssim& &\int
	\left(\sum_{\theta_{\alpha,d}}\ 
	\sup_{\omega\subseteq\theta_{\alpha,d}}\
	\lp{S^{\pm\omega}_{1,d}u^\alpha\, (t)}{L^\infty_x}^2\right)^\frac{1}{2}
	\cdot\  \sup_{\theta_{\alpha,d}} \lp{\nabla 
	B^{\theta_{\alpha,d}}_{|\theta_{\alpha,d}|}
	S_{\lambda,\bullet < \min\{d,c\}} v\, (t)}{L^2_x}\, dt \ , \notag\\
	&\lesssim&  |\theta_{\alpha,d}| \ &\int
	\left(\sum_{\theta_{\alpha,d}}\ 
	\sup_{\omega\subseteq\theta_{\alpha,d}}\
	\lp{S^{\pm\omega}_{1,d}u^\alpha\, (t)}{L^\infty_x}^2\right)^\frac{1}{2}
	\cdot\  \lp{\nabla S_{\lambda,\bullet < \min\{d,c\}} 
	\langle\Omega\rangle v\, (t)}{L^2_x}\, dt \ , \notag\\
	&\lesssim&  |\theta_{\alpha,d}| \ &\int
	\left(\sum_{\theta_{\alpha,d}}\ 
	\sup_{\omega\subseteq\theta_{\alpha,d}}\
	\lp{S^{\pm\omega}_{1,d}u^\alpha\, (t)}
	{L^\infty_x}^2\right)^\frac{1}{2} dt
	\ \cdot\  \lp{\nabla S_{\lambda,\bullet < \min\{d,c\}} 
	\langle\Omega\rangle v}{L^\infty(L^2)} \ , \notag\\
	&\lesssim&  |\theta_{\alpha,d}| \ &\int
	\left(\sum_{\theta_{\alpha,d}}\ 
	\sup_{\omega\subseteq\theta_{\alpha,d}}\
	\lp{S^{\pm\omega}_{1,d}u^\alpha\, (t)}
	{L^\infty_x}^2\right)^\frac{1}{2} dt
	\ \cdot\  \lambda \ \lp{v}{F_{\Omega,\lambda}} \ . \notag
\end{align}
Finally, assuming that $\{\theta_{\alpha,d}\}$ and the $u^\alpha$
minimize the $Z_1$ norm for this particular $u$, we get that:
\begin{align}
	\text{(L.H.S.)}\eqref{special_L1L2_omega_on_low} 
	\ &\lesssim \ \sum_{\alpha , d} \ |\theta_{\alpha,d}|\ \int
	\left(\sum_{\theta_{\alpha,d}}\ 
	\sup_{\omega\subseteq\theta_{\alpha,d}}\
	\lp{S^{\pm\omega}_{1,d}u^\alpha\, (t)}
	{L^\infty_x}^2\right)^\frac{1}{2} dt
	\ \cdot \ \lambda \lp{v}{F_{\Omega,\lambda}}\ , \notag\\
	&\lesssim \ \lp{v}{Z_1}\cdot \ \lambda \ 
	\lp{v}{F_{\Omega,\lambda}}\ .\notag
\end{align}
\end{proof}\ret


\begin{proof}[proof of estimate \eqref{special_L1L2_omega_on_high}]
Fixing $d$ we use the angular decomposition 
\eqref{LH_wide_angle_decomp}
in conjunction with the multiplier Lemma \ref{generic_mult_lemma}, several
rounds of H\"olders inequality, and the Strichartz estimate 
\eqref{high_F_energy_str_est} to compute that:
\begin{align}
	& & &\int \lp{S_{\lambda , \bullet < \min\{d,c\}} (S_{1,d}u\, 
        \nabla S_{\lambda,\bullet < \min\{d,c\}} v)(t)}{L^2_x}\, dt 
	\ , \notag \\
	&\lesssim& &\int \left(\sum_\omega\
	\lp{S^{\pm\omega}_{1,d}u\, (t)}{L^\infty_x}^2\cdot
	\lp{\nabla B^\omega_{d^\frac{1}{2}}
	S_{\lambda,\bullet < \min\{d,c\}} v\, (t)}{L^2_x}^2
	\right)^\frac{1}{2} dt \ , \notag \\
	&\lesssim& &\int \sup_\omega 
	\lp{S^{\pm\omega}_{1,d}u\, (t)}{L^\infty_x}\cdot
	\left(\sum_\omega\ \lp{\nabla B^\omega_{d^\frac{1}{2}}
	S_{\lambda,\bullet < \min\{d,c\}} v\, (t)}{L^2_x}^2
	\right)^\frac{1}{2} dt \ , \notag \\
	&\lesssim& &\int \sup_\omega 
	\lp{S^{\pm\omega}_{1,d}u\, (t)}{L^\infty_x}\, dt\ \cdot \
	\lp{\nabla S_{\lambda,\bullet < \min\{d,c\}} v}
	{L^\infty(L^2)} \ , \notag\\
	&\lesssim& &\int \sup_\omega 
	\lp{S^{\pm\omega}_{1,d}u\, (t)}{L^\infty_x}\, dt\ \cdot \
	\lambda \ \lp{v}{F_\lambda} \ . \notag
\end{align}
Summing this last expression over $d$ yields:
\begin{align}
	\text{(L.H.S.)}\eqref{special_L1L2_omega_on_high} 
	\ &\lesssim \ \sum_d \ \int \sup_\omega 
	\lp{S^{\pm\omega}_{1,d}u\, (t)}{L^\infty_x}\, dt\ \cdot \
	\lambda \ \lp{v}{F_\lambda}\ , \notag\\
	&\lesssim \ \lp{v}{Z_{\Omega,1}}\cdot \ \lambda \ 
	\lp{v}{F_{\lambda}}\ .\notag
\end{align}

\end{proof}

\ret\ret


\end{document}